\documentclass [twoside,reqno, 12pt] {amsart}
  
\usepackage{graphicx}
\usepackage{amsfonts}
\usepackage{amssymb}
\usepackage{a4} 

\newcommand{\mattiHOX}[1]{}

\newcommand{\HOX}[1]{}

\def\prln{{\it Phys. Rev. Lett.}\,}

\def\lambdapoistettu{}

\def\uusif{f}
\def\principalPhi{\Phi}
\def\uusiPhi{\Phi}
\def\uusiPsi{\Psi}
\def\uusivarphi{\varphi}
\def\uusieta{\eta}
\newtheorem{thm}{Theorem}[section]
\newtheorem{theorem}[thm]{Theorem}
\newtheorem{corollary}[thm]{Corollary}
\newtheorem{lemma}[thm]{Lemma}
\newtheorem{proposition}[thm]{Proposition}

\newtheorem{definition}[thm]{Definition}

\theoremstyle{definition}

\numberwithin{equation}{section}

\newcommand{\mlltext}{}
\newcommand{\lltext}{}
\newcommand{\ltext}{}
\newcommand{\mmmtext}{}

\def \ol {\overline}
\def \b {{\beta}}

\def \K {{\mathcal K}}

\def\re{\hbox{Re}\,}
\def\im{\hbox{Im}\,}
\def\Re{\hbox{Re}\,}         
\def\Im{\hbox{Im}\,}
\def\A{{\mathcal A}}

\renewcommand{\O}{{\mathcal O}}
\def\Om{\Omega}
\renewcommand{\H}{{\mathcal H}}

\newcommand{\R}{{\mathbb R}}
\newcommand{\D}{{\mathbb D}}
\newcommand{\Sbb}{\mathbb S}

\newcommand{\C}{{\mathbb C}}
\newcommand{\N}{{\mathbb N}}
\newcommand{\proofbox}{{\hfill $\Box$}}

\newcommand{\B}{{\mathcal N}}

\def\hat{\widehat}
\def\tilde{\widetilde}
\def\bfo{\begin {eqnarray*} }
\def\efo{\end {eqnarray*} }
\def\ba{\begin {eqnarray*} }
\def\ea{\end {eqnarray*} }
\def\beq{\begin {eqnarray}}
\def\eeq{\end {eqnarray}}
\def\supp{\hbox{supp}\,}

\def\diag{\hbox{diag }}
\def\det{\hbox{det}\,}
\def\tr{\hbox{tr}\,}
\def\bra{\langle}

\def\cet{\rangle}

\def\e{\varepsilon}
\def\p{\partial}
\def\dbar{\overline \partial}
\def\a{\lambda}
\def\la{\lambda}

\def\Z{{\mathbb Z}}

\def\S{\Sbb}
\begin{document}

\title[Invisibility and Visibility]{The borderlines of the invisibility and  visibility
for Calder\'on's inverse
problem 
}

\author[Astala]{Kari Astala}
\address
        {Department of Mathematics and Statistics \\
         University of Helsinki\\
         P.O. Box 68 \\
         FI-00014   Helsinki\\
         Finland}

\author[Lassas]{Matti Lassas}
\address
        {{\it         E-mail: Kari.Astala@helsinki.fi, Matti.Lassas@helsinki.fi, ljp@rni.helsinki.fi}}


\author[P\"aiv\"arinta]{Lassi P\"aiv\"arinta}

\date{}
\maketitle

\begin{abstract}
We  consider the determination
of a conductivity function in a two-dimensional domain from the Cauchy data of the solutions of the conductivity equation on the boundary. 
We prove  uniqueness results for this inverse problem, posed by Calder\'on,  for  conductivities that
are degenerate, that is, they may
 not be bounded
from above or below. 
Elliptic equations with such coefficient functions are essential for physical models
used in
transformation optics and the metamaterial constructions.
In particular, for scalar conductivities 
we solve the inverse problem in a class which is larger than $L^\infty$.
Also, we give new counterexamples for the uniqueness of the  inverse conductivity problem. 

We say that a conductivity is visible if the inverse problem is solvable so that the inside of the domain can be uniquely determined, up to a change of
coordinates, using  the boundary measurements.
The present counterexamples for the inverse problem have been related to  the invisibility cloaking. 
This means that there are  conductivities for which a part of the domain is   shielded from
detection via boundary measurements. Such conductivities are called  invisibility cloaks.

In the present paper we identify 
  the borderline of the visible conductivities 
and the  borderline of invisibility 
cloaking conductivities.
Surprisingly, these borderlines are not the same. We show that between the visible and the cloaking conductivities
there are the electric holograms, conductivities which create an illusion of a non-existing  body.
The electric holograms give counterexamples for the uniqueness of the inverse problem  which are less
degenerate than the previously known ones. 
These examples are constructed using transformation optics and the inverse maps 
of the Iwaniec-Martin mappings. The uniqueness results are based on
combining the complex geometrical optics, the properties of the mappings with subexponentially integrable
distortion, and the Orlicz space techniques.

\end{abstract}

\section{Introduction and main results}

Invisibility cloaking has been a very topical 
subject in recent studies in mathematics, physics, and 
material science  \cite{AE,GKLU,GLU3,MN,Le,MOV,PSS1,Weder}.
By invisibility cloaking we mean the possibility, 
both theoretical and practical, of shielding a
region or object from
detection via electromagnetic fields. 

{\mlltext The counterexamples
for inverse problems and the proposals for invisibility cloaking are closely related.
In 2003,
before the appearance of practical possibilities for cloaking,
 it was shown in  \cite{GLU2,GLU3} that passive
objects  can be coated with a layer of material with a degenerate
conductivity which makes the object undetectable by 
the electrostatic boundary
measurements.
These constructions were 
based on the blow up maps
and gave counterexamples for the uniqueness of  inverse conductivity
problem in the three and
higher dimensional cases.
In two dimensional case, the mathematical theory of the cloaking examples for conductivity
equation have been studied in \cite{KSVW,KOVW,LZ,Ng}. 


The interest in cloaking was
raised in particular in 2006 when it was realized 
that practical
cloaking constructions are possible using so-called metamaterials
which allow fairly arbitrary specification of electromagnetic material
parameters. The construction of Leonhardt
\cite{Le} was based on conformal mapping on a non-trivial Riemannian surface. 
At the same time, Pendry et al
\cite{PSS1}  proposed a cloaking construction
for Maxwell's equations using a blow up map 
 and the idea was demonstrated in laboratory experiments
\cite{SMJCPSS}. There are also other suggestions for cloaking  based on active sources \cite{MOV} or
negative material parameters \cite{AE,MN}.


In this paper we consider both new counterexamples and uniqueness  
results for inverse problems.}
We study  Calder\'on's inverse problem
in the two dimensional case, that is, 
 the
question  whether an  unknown conductivity distribution inside a domain
can be determined  from the voltage and current
measurements made on the boundary.
Mathematically the measurements correspond to the knowledge of the
Dirichlet-to-Neumann map  $ \Lambda_\sigma $ (or the quadratic form)
associated to $\sigma$, i.e., the map taking the Dirichlet boundary
values of the solution of
the conductivity equation
\begin{equation}\label{johty}
\nabla\cdot \sigma(x)\nabla u(x) = 0
\end{equation}
 to the
corresponding Neumann boundary values,
\beq\label{eq: DN map}
 \Lambda_\sigma : \  u|_{\p \Omega}\mapsto \nu \, \,\cdotp \sigma \nabla u|_{\p \Omega}.
\eeq
In the classical theory of the problem, the conductivity $\sigma$
is bounded uniformly from above and below. 
The problem was originally proposed by Calder\'on \cite{C} in 1980.
Sylvester and Uhlmann \cite{SU} proved the unique identifiability of the
conductivity in dimensions three and higher for isotropic conductivities which
are $C^\infty$-smooth, and Nachman gave a reconstruction method
\cite{N}. In three dimensions or higher unique identifiability of the
conductivity is known for conductivities with $3/2$ derivatives
\cite{PPU}, \cite{BT} and $C^{1,\alpha}$-smooth conductivities which
are  $C^\infty$ smooth outside surfaces on which they have
conormal singularities \cite{GLU1}. The problems has also been solved
with measurements only on a part of the boundary \cite{KSU}.

   In two dimensions the first global solution of the inverse conductivity problem
  is due to Nachman \cite{N1}   for conductivities with two derivatives. In this seminal paper
   the $\dbar$ technique was first time  used in the study of Calderon's inverse problem.
The smoothness requirements were reduced in \cite{BU} to Lipschitz
conductivities.
Finally, in  \cite{AP} the uniqueness of the inverse 
problem  was proven {\mlltext in the form that the problem was
originally formulated in \cite{C}, i.e., for}
 general
isotropic conductivities in $L^\infty$ which are bounded from below and above by positive constants.

$$
\begin{array}{ll}
\hspace{-2cm}\includegraphics[width=.4\linewidth]{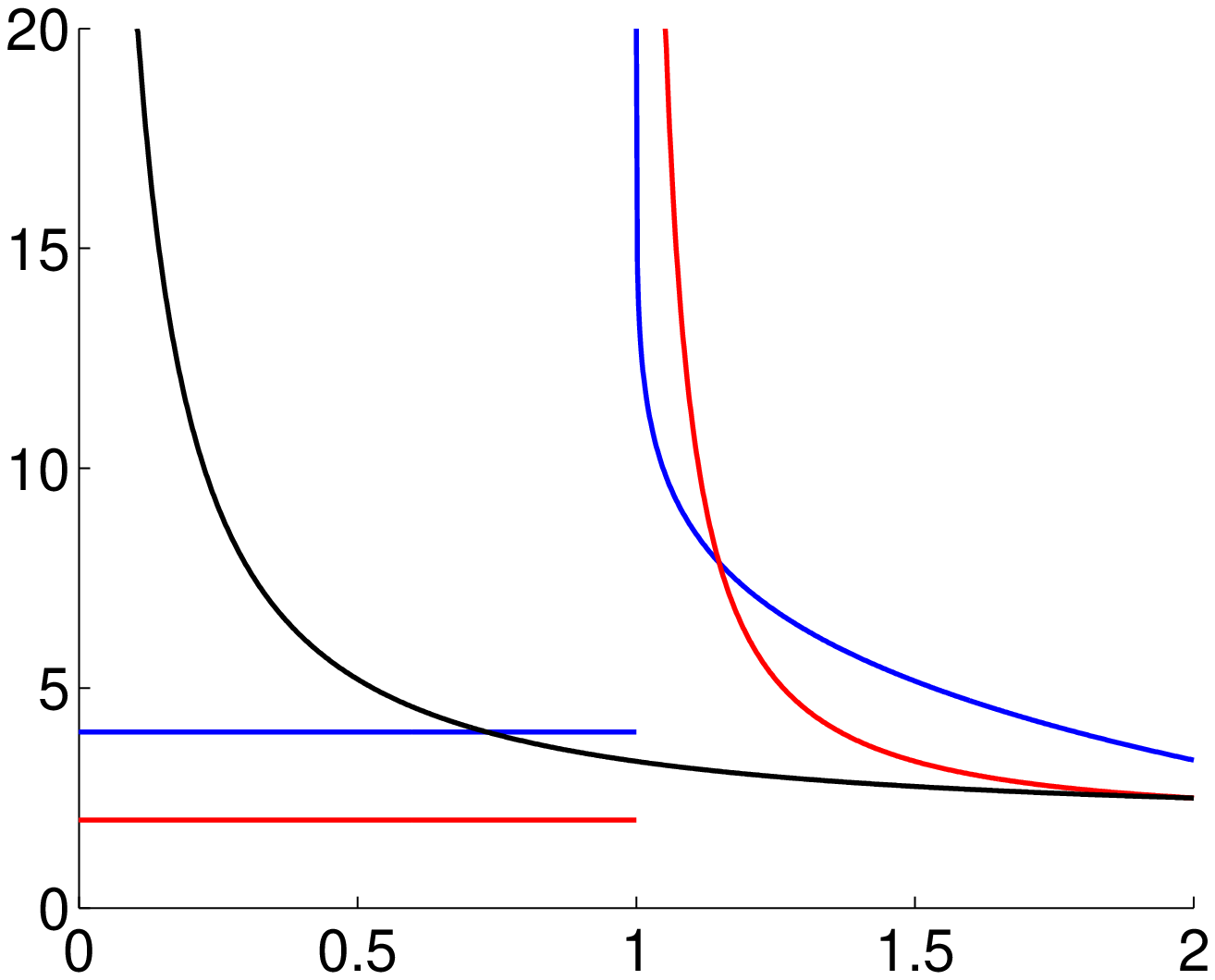}\hspace{-.7cm}&\includegraphics[width=.4\linewidth]{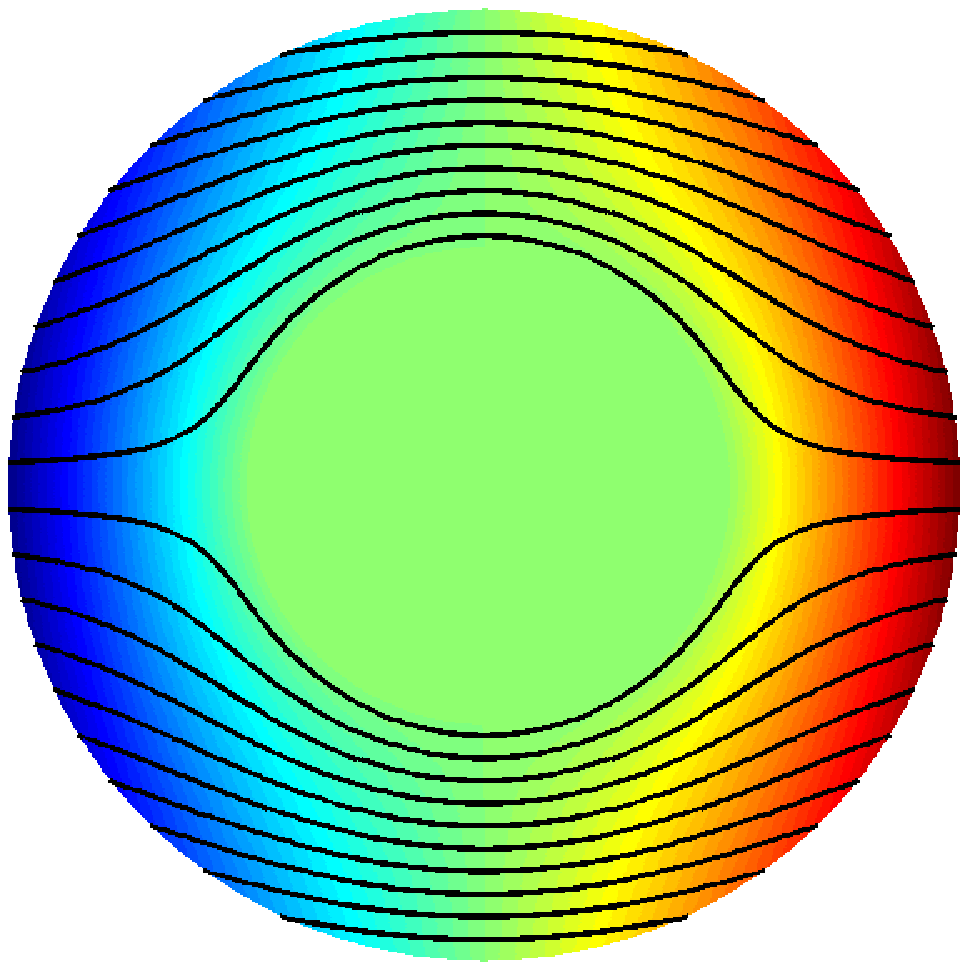}\hspace{-1cm}
\includegraphics[width=.4\linewidth]{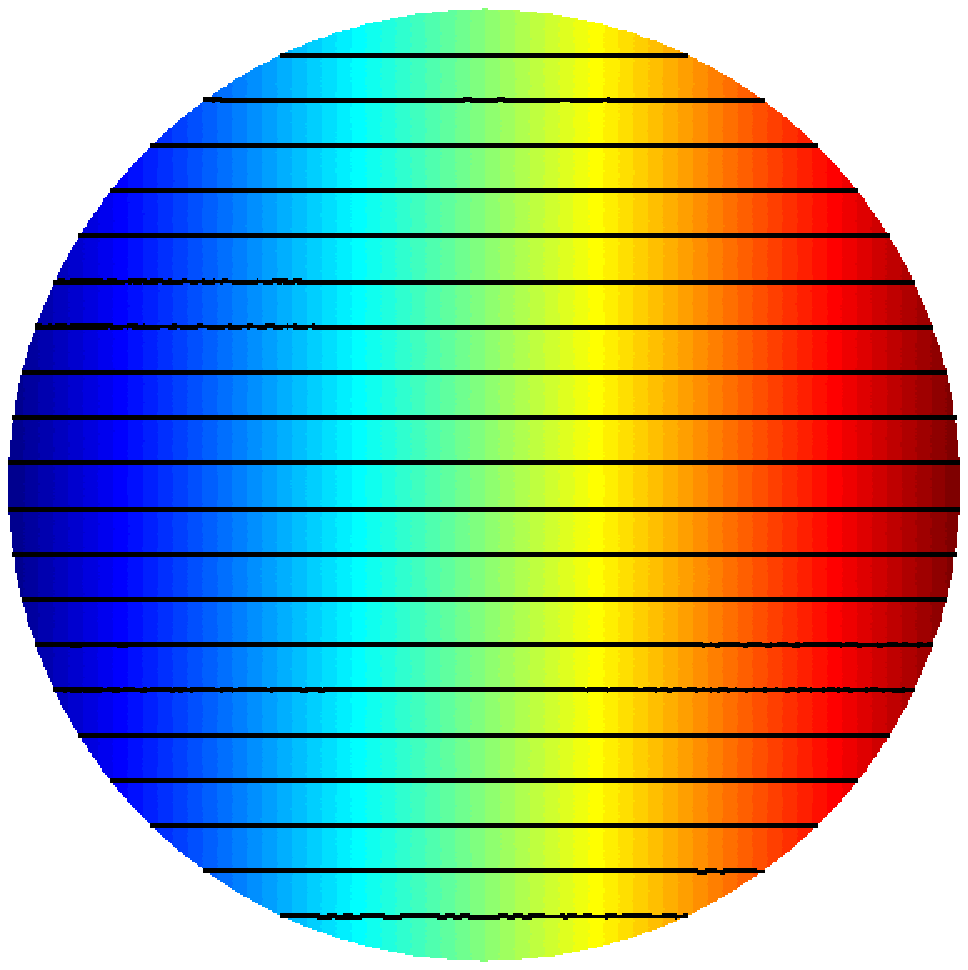}\hspace{-3cm}\\
\hspace{-2cm}\includegraphics[width=.4\linewidth]{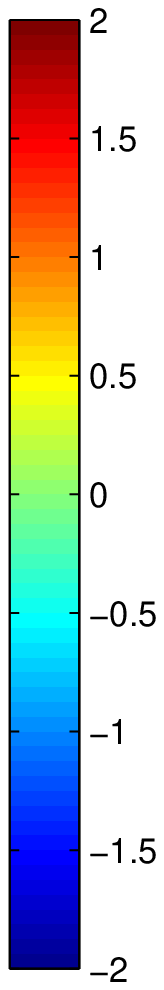}\hspace{-.7cm} &\includegraphics[width=.4\linewidth]{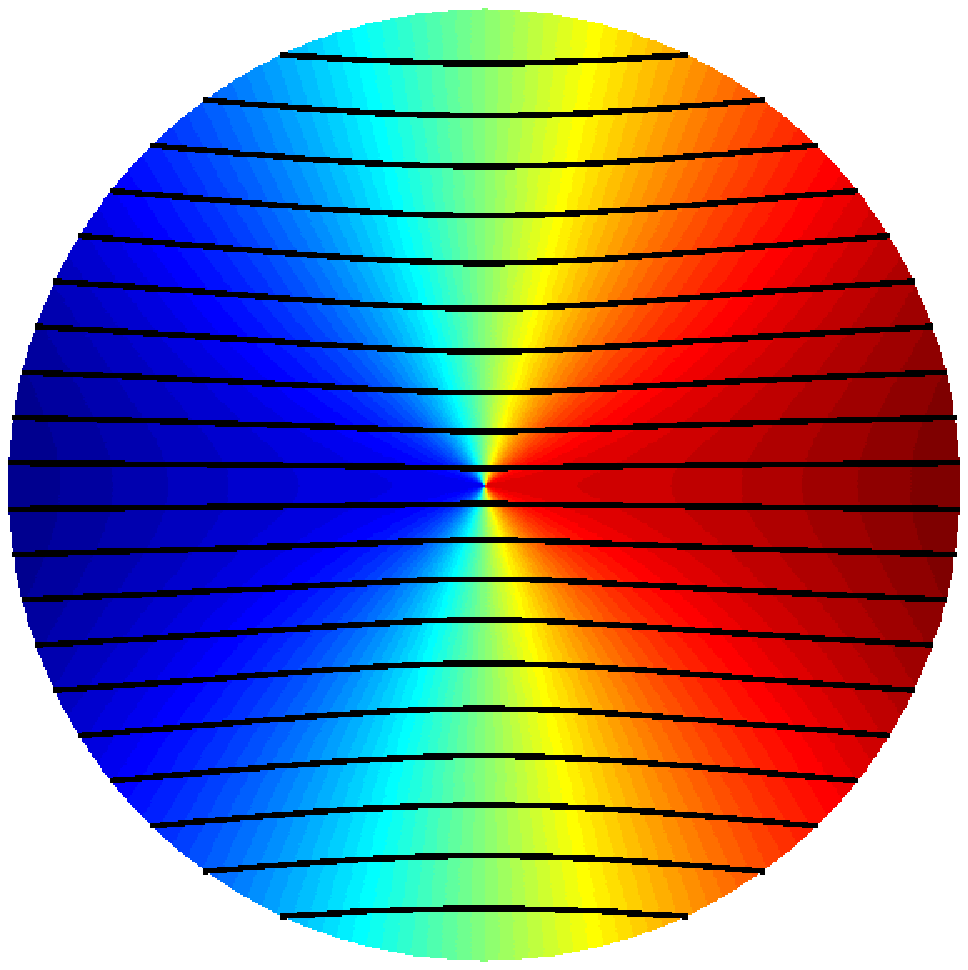}\hspace{-1cm}
\includegraphics[width=.4\linewidth]{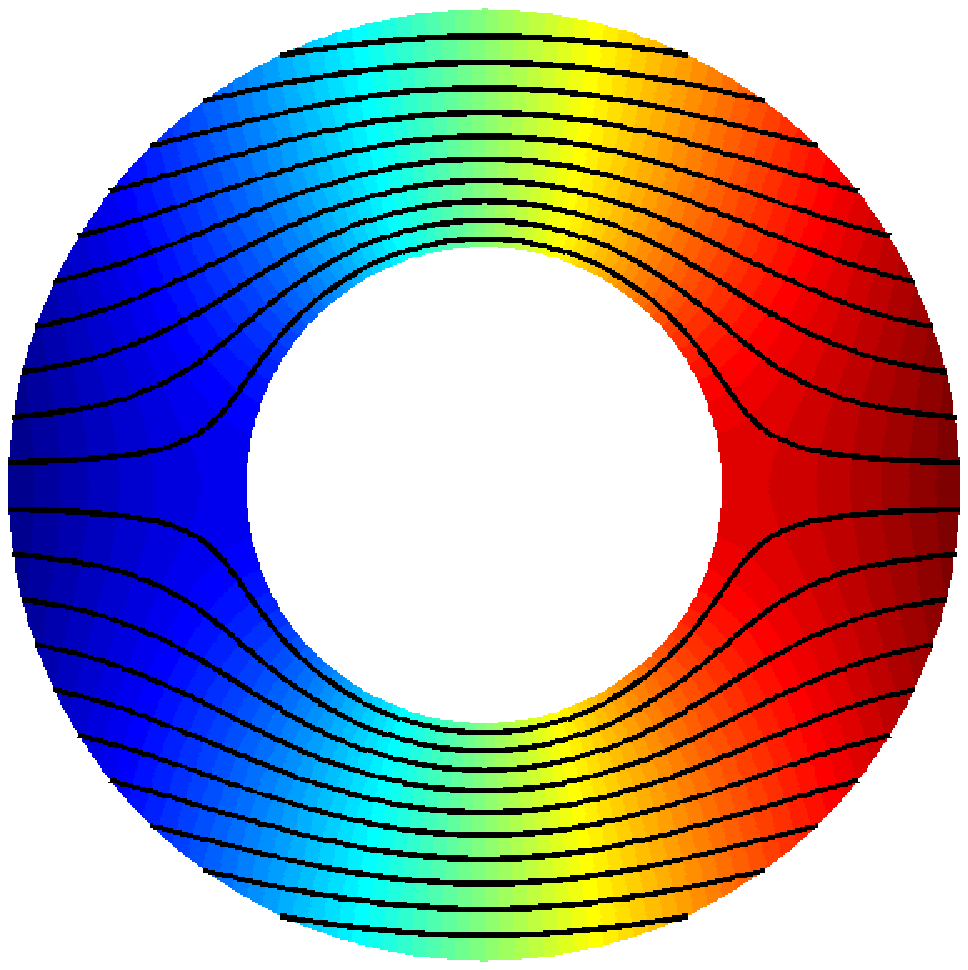}\hspace{-3cm}\ 
\end{array}$$
{{\bf Figure 1. Left}.\  $\tr (\sigma)$ \it of three radial and singular conductivities on the positive $x$ axis. The curves
correspond to the invisibility cloaking conductivity   (red), with the singularity $\sigma^{22}(x,0)\sim (|x|-1)^{-1}$ for
$|x|>1$,
 a visible conductivity (blue) with a $\log \log$ type singularity at $|x|=1$, 
 and an electric hologram  (black) with the conductivity   having the singularity $\sigma^{11}(x,0)\sim |x|^{-1}$. 
 {\bf Right, Top.}\ 
 All measurements on the boundary  of  
 the invisibility cloak (left) coincide with the measurements for the homogeneous disc (right).
 The color shows the value of the solution 
 $u$ with the boundary value $u(x,y)|_{\p B(2)}=x$
 and the black curves are the integral curves
 of the current $-\sigma\nabla u$.  
 {\bf Right, Bottom.}\ All measurements on the boundary of  
 the electric hologram (left) coincide with the measurements for  an isolating disc covered with the homogeneous medium (right). The solutions and the current lines corresponding to
 the boundary value $u|_{\p B(2)}=x$ are shown.
}
\bigskip

The Calder\'on problem with an anisotropic, i.e., matrix-valued,
conductivity that
is uniformly bounded from above and below
has been studied in  two dimensions \cite{S, N1,LU,  ALP,IUY} and
   in dimensions $n\geq 3$
\cite{LeU,LU, Salo}.   For example, for the
anisotropic inverse  conductivity problem in the two dimensional case it is
known that
the Dirichlet-to-Neumann map determines a regular conductivity tensor up to a
diffeomorphism
$ F:\overline\Omega\to
\overline\Omega$, i.e.,\ {\lltext one can obtain an image of the interior of $\Omega$ in
deformed coordinates}. This implies that the inverse problem is not uniquely
solvable, but the non-uniqueness of the problem can be characterized.
We note that the problem in higher dimensions is presently
solved only {\mlltext in special cases, like when} the conductivity
is real analytic.




In this work we will study the inverse conductivity problem
in the two dimensional case with  
{\ltext degenerate conductivities}. 
{\lltext Such conductivities appear in physical models where the medium
varies continuously from a perfect conductor to a perfect insulator. As an 
example, we may consider a case where the conductivity goes to
zero or to infinity near $\p D$ where $D\subset\Omega$ is a smooth open set.}
 We ask what kind of degeneracy prevents solving the
inverse problem, that is, we study what is the border of visibility.
We also ask what kind of degeneracy makes it even possible 
to coat of an arbitrary object so that it appears the same
as a homogeneous  body in all static measurements, that is, 
we study what is  the border of the invisibility cloaking. {\ltext Surprisingly,
these borders are not the same;}
We identify these borderlines and
show that between them 
there are the electric holograms,  that is, the conductivities  creating an illusion of a non-existing  body
(see Fig.\ 1).
These conductivities are the counterexamples for the unique solvability of inverse problems 
for which 
even the topology of the domain can not be determined using boundary measurements.
Our main result for {\lltext the uniqueness of the inverse problem are given in Theorems \ref{main 2b both}, \ref{main 2b isotropic}, and} \ref{main 2b trace}
and the counterexamples are formulated in 
Theorems \ref{prop: invisibility}  and \ref{thm: hologram}.

{\mlltext The cloaking constructions have given rise for the design technique called the transformation optics. 
The metamaterials build to operate
at microwave frequencies \cite{SMJCPSS} and  near the optical frequencies \cite{Er,Val} are inherently
prone to dispersion, so that realistic cloaking must currently be considered
as occurring at  very narrow range of wavelengths.
Fortunately, in many physical applications the materials need
to operate only near a single frequency. 
In particular, the cloaking type constructions have 
inspired suggestions for possible
devices producing extreme effects on wave propagation, including invisibility cloaks for magnetostatics  \cite{WP},
 acoustics \cite{CummerSchurig,ChenChan} and quantum mechanics  \cite{Zhang,GKLU6}; field rotators \cite{ChenRotate};  electromagnetic wormholes \cite{GKLUWorm}; 
  invisible sensors \cite {AE2,GKLUSensor}; superantennas \cite{LeT}; perfect absorbers \cite{Lan};
  {\lltext and wave amplifiers \cite{GKLU7}}. 
It has turned out that the designs that are based on well posed mathematical models, {\ltext e.g.\ approximate cloaks,} 
 have excellent properties {\ltext when compared to ad hoc constructions}. Due to this, it is important to know what is the exact degree of non-regularity which
 is needed for invisibility cloaking or solving the inverse problems. 

 Finally, we note that the differential equations with degenerate coefficients modeling
 cloaking devices have 
 turned out to have interesting properties, such as non-existence results for solutions with non-zero sources
  \cite{GKLU,LiZ} and the local and non-local hidden boundary conditions \cite{LZ,Ng}.
}

The structure of the paper is the following. The main results and the formulation
of the boundary measurements  are presented in
the first section. The proofs for the existence of the solutions of
the direct problem as well as for the new counterexamples and
the invisibility cloaking examples with a non-smooth background are given
in Section 2.  The uniqueness result for the isotropic conductivities 
is proven in Sections 3-4 and the reduction of the general problem to the isotropic case is
shown in Section 5. In Sections 3-5, the degeneracy of the conductivity 
causes that the exponentially growing solutions, the standard
tools used to study Calderon's inverse problem,  can not be
constructed using purely microlocal or  functional analytic methods.
Because of this we will  extensively need the topological properties of the solutions:
By Stoilow's theorem the solutions are compositions of analytic functions and
homeomorphisms. Using this, the continuity properties of the weakly
 monotone maps, and the Orlicz-estimates holding for homeomorphisms
  we prove the existence of the solutions
 in the Sobolev-Orlicz spaces.  These  spaces are chosen so
that we can obtain subexponential asymptotics for the families of exponentially
growing solutions needed in the $\dbar$ technique used to solve the inverse
problem.

\subsection{Definition of measurements and solvability}\label{solva}

Let $\Omega\subset \R^2$ be a bounded {\ltext simply connected} domain with a smooth boundary.
Let $\Sigma=\Sigma(\Omega)$ be the class of measurable matrix valued functions 
$\sigma: \Omega \to M$, where $M$ is the set of 
generalized matrices $m$ of the  form
\ba
m=U \left(\begin{array}{cc} \lambda_1 & 0 \\
 0&  \lambda_2\\ \end{array}\right) U^t
\ea
 where $U\in \R^{2\times 2}$ is an orthogonal matrix, $U^{-1}=U^t$
and $\lambda_1,\lambda_2\in [0,\infty) $
{\lltext We denote by
$W^{s,p}(\Omega)$ and $H^s(\Omega)=W^{s,2}(\Omega)$ the standard Sobolev spaces.}

In the following, let  $dm(z)$ denote the Lebesgue measure in $\C$ and $|E|$ be
the   Lebesgue measure of the set $E\subset \C$.
Instead of defining the   Dirichlet-to-Neumann operator 
for the above  conductivities, we consider the corresponding quadratic forms.

\begin{definition}\label{def: 1A} {\rm Let $h\in H^{1/2}(\partial \Omega)$.
The Dirichlet-to-Neumann quadratic form corresponding to the conductivity
$\sigma\in \Sigma(\Omega)$ is
given by
\beq\label{def: DN map}
Q_\sigma[h]=\inf A_\sigma[u]\quad\hbox{where, }
A_\sigma[u]=\int_\Omega \sigma(z)\nabla u(z)\, \,\cdotp\nabla u(z)\,\,dm(z),
\eeq
and the infimum is taken over real valued
$u\in L^1(\Omega)$ such that $ \nabla u\in L^1(\Omega)^3$
and $u|_{\partial \Omega}=h$. In the case where $Q_\sigma[h] < \infty$ and $A_\sigma[u]$ reaches
its minimum at some $u$, we say that $u$ is a $W^{1,1}(\Omega)$
solution of the conductivity problem.}
\end{definition}

In the case when $\sigma$ is smooth, bounded from below and above by positive
constants,
$Q_\sigma[h]$ is the quadratic form corresponding the Dirichlet-to-Neumann
map (\ref{eq: DN map}),
\beq\label{eq: Q forms}
Q_\sigma[h] =  \int_{\partial \Omega} h \, \Lambda_{\sigma} h\,  dS ,
\eeq
where $dS$ is the length measure on $\p \Omega$.
Physically, $Q_\sigma[h]$ corresponds to the power needed to keep voltage
$h$ at the boundary. {\lltext For smooth conductivities bounded from below}, for every $h \in H^{1/2}(\partial \Omega)$ the integral $A_\sigma[u]$ always has a unique minimizer $u \in H^1(\Omega)$
with $u|_{\p \Omega}=h$, which is also a distributional solution to (\ref{johty}).  Conversely, for functions $u \in H^1(\Omega)$ their traces lie in $H^{1/2}(\partial \Omega)$. It is for this reason that we chose to consider   the $H^{1/2}$-boundary functions also in the most general case.
We interpret that the Dirichlet-to-Neumann form
corresponds to the idealization of the boundary measurements 
 for $\sigma\in \Sigma(\Omega)$.

We note that the conductivities  studied in the context of cloaking
are not even in $L^1_{loc}$. 
As $\sigma$ is unbounded it is possible that $Q_\sigma[h]=\infty$. Even
if $Q_\sigma[h]$ is finite, the minimization problem in (\ref{def: DN map}) may generally have no 
minimizer and even if they exist the minimizers
need not be distributional solutions to (\ref{johty}). 
However, if the singularities of  $\sigma$ {\ltext are not too strong}, minimizers satisfying (\ref{johty}) do always  exist. 
To show this, we need to define a suitable subclasses 
of  degenerate conductivities.

Let $\sigma\in \Sigma(\Omega)$. We start with precise quantities describing the possible degeneracy or loss of uniform ellipticity. First, a natural measure of
the anisotropy of the conductivity $\sigma$ at $z\in \Omega$ is
\ba
k_\sigma(z)={\sqrt{\frac{\lambda_1(z)}{\lambda_2(z)}}},
\ea 
where $\lambda_1(z)$ and  $\lambda_2(z)$ are the eigenvalues
of the matrix $\sigma(z)$, $\lambda_1(z)\geq \lambda_2(z)$. If we want to simultaneously control both the size and the anisotropy, this is measured by the ellipticity $K(z)=K_\sigma(z)$  of $\sigma(z)$, i.e. the smallest number $1 \leq K(z) \leq \infty$ such that 
\beq\label{uusi1}
 \frac 1{K(z)} |\xi|^2\leq \xi \,\cdotp \sigma(z)\xi\leq K(z) |\xi|^2,\quad \hbox{for all }\xi \in \R^2.
\eeq
For a general, positive matrix valued function $\sigma(z)$ we have at $z\in \Omega$
\beq\label{uusi2}
K(z) = k_\sigma(z) \max \{ [\det \sigma(z)]^{1/2}, [\det \sigma(z)]^{-1/2} \}.
\eeq
Consequently, we always have the following simple estimates.

\begin{lemma}\label{lem: Inequality} For any measurable matrix function $\sigma(z)$ we have
\ba \quad  \frac{1}{4} \left[ \tr \sigma(z)+\tr (\sigma(z)^{-1}) \right] \leq K(z) \leq  \tr \sigma(z)+\tr (\sigma(z)^{-1}).
\ea
\end{lemma}

{\bf Proof.} 
Let $\lambda_{max}$ and $\lambda_{min}$ be the  eigenvalues
of $\sigma=\sigma(z)$. Then
$K(z)=\max(\lambda_{max},\lambda_{min}^{-1})$. Since $ \tr \sigma(z) = \lambda_{max} +\lambda_{min}$ and $\tr (\sigma(z)^{-1}) =\lambda_{max}^{-1}+\lambda_{min}^{-1}$, the claim  follows.
\hfill \proofbox
\medskip

{\lltext Due to Lemma \ref{lem: Inequality} we use}  the  quantity $\;  \tr \sigma(z)+\tr (\sigma(z)^{-1}) \;$  as a  measure of size and anisotropy of  $\sigma(z)$.
\medskip

For the degenerate elliptic equations it may be that 
the optimization problem (\ref{def: DN map})  has a minimizer which satisfies the conductivity equation
but this solution may not have  the standard $W^{1,2}_{loc}$ regularity.
Therefore more subtle smoothness estimates are required. We start with the  exponentially integrable conductivities, and the natural energy estimates they require. As an important  consequence we will see  the correct  Orlicz-Sobolev regularity to work with. 
These observations are based on the following elementary inequality.

\begin{lemma}\label{lem: Kari 1}
 Let $K\geq 1$ and $A\in \R^{2\times 2}$ be a symmetric matrix satisfying
 \ba
 \frac 1K |\xi|^2\leq \xi \,\cdotp A\xi\leq K |\xi|^2,\quad \xi\in \R^2.
 \ea
 Then for every $p>0$
 \ba
 \frac {|\xi|^2} {\log(e + |\xi|^2)}+ \frac {|A\xi|^2} {\log(e + |A\xi|^2)}\leq
 \frac 2p\left( \xi \,\cdotp A\xi+e^{pK}\right).
 \ea  
\end{lemma}
  
\noindent {\bf Proof.}
Since  $K\geq 1$ and $t\mapsto t/\log(e+t)$ is an increasing function, we have
\ba
\frac {|\xi|^2} {\log(e + |\xi|^2)}&\leq& \frac {K \xi \,\cdotp A\xi} {\log(e + K \xi \,\cdotp A\xi)}\\
&\leq& \frac 1p \left(\frac {\xi \,\cdotp A\xi} {\log(e +  \xi \,\cdotp A\xi)}\right) pK\\
&\leq&
 \frac 1p\left( \xi \,\cdotp A\xi+e^{pK}\right),
\ea 
where the last estimate follows from the inequality
\ba
ab\leq a\log(e+a)+e^b,\quad a,b\geq 0.
\ea 
Moreover, as $K$ is at least as large as the maximal eigenvalue of $A$,
we have $|A\xi|^2\leq K\xi \,\cdotp A\xi$. Thus we see as above 
that
\ba
\frac {|A\xi|^2} {\log(e + |A\xi|^2)}&\leq& \frac {K \xi \,\cdotp A\xi} {\log(e + K \xi \,\cdotp A\xi)}
\leq
 \frac 1p\left( \xi \,\cdotp A\xi+e^{pK}\right).
\ea 
Adding the above estimates together proves the claim.\hfill \proofbox 
\medskip

Lemma \ref{lem: Kari 1} implies in particular that if $\sigma(z)$ is symmetric matrix valued function
satisfying (\ref{uusi1}) for a.e.\ $z\in \Omega$
and $u\in W^{1,1}(\Omega)$, then always
\beq\label{eq: Kari 2}\\ \nonumber
 p\int_\Omega \frac {|\nabla u(z)|^2} {\log(e + |\nabla u(z)|^2)}dm(z)
 \leq \int_\Omega  \nabla u(z) \,\cdotp \sigma(z)\nabla u(z)dm(z)+
 \int_\Omega  e^{pK(z)}dm(z),\hspace{-1cm}\\
 \nonumber
\hspace{-.1cm}
 p\int_\Omega \frac {|\sigma(z)\nabla u(z)|^2} {\log(e + |\sigma(z)\nabla u(z)|^2)}dm(z)
 \hspace{-.1cm}\leq\hspace{-.1cm} \int_\Omega  \nabla u(z) \,\cdotp \sigma(z)\nabla u(z)dm(z)+
 \int_\Omega  e^{pK(z)}dm(z).\hspace{-1cm}
\eeq
Note that these inequalities are valid whether $u$ is a solution of the conductivity
equation or not!

{\lltext Due to (\ref{eq: Kari 2}),} we see that to analyze finite energy solutions corresponding to a
singular conductivity of exponentially integrable ellipticity, we are naturally led to consider the regularity gauge
\beq\label{gauge2}
Q(t)=  \frac {t^2}
{\log(e + t)},\quad t\geq 0.
\eeq
We say accordingly that $f$ belongs to the Orlicz space $W^{1,Q}(\Omega)$, cf.\ Appendix, if $f$ and 
 its first distributional derivatives are in $L^1(\Omega)$
and
\ba
 \int_{\Omega}  \frac {|\nabla f(z)|^2}
{\log(e + |\nabla f(z)|)} \; \,dm(z)  < \infty .
\ea   

The first existence result for solutions corresponding to degenerate conductivities is now given as follows. 

\begin{theorem}\label{thm: Kari 3} 
 Let $\sigma(z)$ be a measurable symmetric matrix valued function.
 Suppose further that for some $p>0$,
\beq\label{eq: C1 estimat}
\int_\Omega \exp(p\, [\tr \sigma(z)+\tr (\sigma(z)^{-1})])\,\,dm(z)=C_1<\infty.
\eeq
Then, if $h\in H^{1/2}(\p \Omega)$ 
{\ltext is such that $Q_\sigma[h]<\infty$} and $X=\{v\in W^{1,1}(\Omega);\ v|_{\p\Omega}=h\},$
 there is a unique $w\in X$ such that
\beq\label{eq: Matti-Kari 1A}
A_\sigma[w]=\inf\{A_\sigma[v]\ ;\ v\in X\}.
\eeq
Moreover, $w$ satisfies the conductivity equation
\beq\label{eq: Matti-Kari 1}
\nabla \,\cdotp \sigma \nabla w=0\quad\hbox{in }\Omega
\eeq
in sense of distributions, and it has the regularity $w\in W^{1,Q}(\Omega)\cap C(\Omega)$.
\end{theorem}

{\ltext Note that if $\sigma$ is bounded near $\p \Omega$ then $Q_\sigma[h]<\infty$
for all $h\in H^{1/2}(\p \Omega)$.}
{\mmmtext Theorem \ref{thm: Kari 3}  is proven in Theorem \ref{thm: Kari 3BB} and Corollary \ref{laitos}
in a more general setting. }


Theorem \ref{thm: Kari 3} yields
that for conductivities satisfying (\ref{eq: C1 estimat})
and being 1 near $\p \Omega$ we can define
the Dirichlet-to-Neumann map
\beq\label{Eq: DN}
\Lambda_\sigma: H^{1/2}(\p \Om)\to H^{-1/2}(\p \Om),\quad
\Lambda_\sigma (u|_{\p \Omega})=\nu \cdot  \sigma\nabla u|_{\p \Omega},
\eeq
where $u$ satisfies (\ref{johty}).

The reader should consider  the exponential condition (\ref{eq: C1 estimat}) as being  close to  the optimal one,  still allowing uniqueness in the inverse  problem.
   Indeed, in view of  Theorem \ref{thm: hologram} and Section  \ref{uusilabel1} below,
the most general situation where the Calder\'on inverse problem can be solved
 involves conductivities whose singularities satisfy a physically interesting small relaxation of the condition (\ref{eq: C1 estimat}). 
 Before solving
inverse problems for conductivities satisfying (\ref{eq: C1 estimat}) we discuss some 
counterexamples.

\subsection{Counterexamples for the unique solvability of the inverse problem}
\label{sec: Counterexamples}

Let  $F:\Omega_1\to\Omega_2,$ $ y=F(x)$ be 
an {\lltext orientation preserving} homeomorphism between domains $\Omega_1,\Omega_2\subset  \C$ 
for which $F$ and its inverse $F^{-1}$ are at least $W^{1,1}$-smooth and let
$\sigma(x)=[\sigma^{jk}(x)]_{j,k=1}^2\in \Sigma(\Omega_1)$ be a conductivity on
$\Omega_1$.
Then the map $F$ pushes $\sigma$ forward to a conductivity $(F_*\sigma)(y)$, defined 
on $\Omega_2$ and  given by
\beq\label{transf}
 (F_*\sigma)(y)=
\frac 1{[\det DF(x)]} DF(x)\, \sigma(x)\, DF(x)^t, \quad x = F^{-1}(y).
\eeq
The main methods for constructing counterexamples to  Calder\'on's problem are based on the following principle.

\begin{proposition}\label{prop: change of coord. part 2}
Assume {\ltext  that  $\sigma,\tilde \sigma\in \Sigma(\Omega)$    
   satisfy (\ref{eq: C1 estimat}), and let
$F:\Omega\to \tilde \Omega$ be a homeomorphism so that
$F$ and $F^{-1}$ are $W^{1,Q}$-smooth 
  and  $C^1$-smooth near the boundary, and $F|_{\p \Omega}=id$.
Suppose that $\tilde \sigma=F_*\sigma$. 
Then $Q_\sigma=Q_{\tilde \sigma}$.}
\end{proposition}

This proposition generalizes the previously known results  \cite{KV} to less smooth diffeomorphisms and 
conductivities and it
follows from Lemma \ref{lem: change of coord.} proven later.

\subsection{Counterexample 1: Invisibility cloaking}
%
{\ltext We consider here invisibility cloaking in general background $\sigma$, that is,
we aim to coat an arbitrary  body with a layer of exotic material so that the coated  body
 appears in measurements the same as the background conductivity $\sigma$. Usually one
 is interested in the case when  the background conductivity $\sigma$  is equal to
the constant $\gamma=1$. However, we consider here a more general
case and assume that $\sigma$ is a $L^\infty$-smooth conductivity in $\overline B(2)$,
$\sigma(z)\geq c_0I,$ $c_0>0$.}
{\mmmtext Here, $B(\rho)$ is an open 2-dimensional disc of radius $\rho$ and center zero
and $\overline B(\rho)$ is its closure. Consider 
a homeomorphism
\beq\label{general blow up}
F:\overline B(2)\setminus\{0\}\to
\overline B(2)
\setminus \K
\eeq
where $\K\subset B(2)$ is a compact set which is the closure of a smooth open set
and suppose $F:\overline B(2)\setminus\{0\}\to
\overline B(2)
\setminus \K$ and its inverse $F^{-1}$ are 
$C^1$-smooth in $\overline B(2)\setminus\{0\}$ and $\overline B(2)
\setminus \K$, correspondingly. We also require that $F(z) = z$ for $z \in \partial  B(2)$.
The standard example {\lltext of invisibility cloaking}
 is the case when $\K=\overline B(1)$ and  the map
\beq\label{blow up}
 F_0(z)=(\frac {|z|}2+1)\frac z{|z|}.
\eeq

Using the map (\ref{general blow up}), we  
define a singular conductivity
\beq\label{eq: sing cond}
\tilde \sigma (z)=\left\{\begin{array}{ll}
(F_*\sigma)(z)
   & \hbox{for }z\in B(2)\setminus \K,\\
\eta(z) & \hbox{for }z\in  \K,\end{array}\right.
\eeq
where $\eta(z)=[\eta^{jk}(x)]$ is any symmetric measurable matrix satisfying
$c_1I\leq \eta(z)\leq c_2I$ with $c_1,c_2>0$. 
The conductivity $\tilde \sigma$ is called the cloaking conductivity obtained from
the transformation map $F$ and background conductivity $\sigma$ and $\eta(z)$
is the conductivity of the cloaked (i.e.\ hidden) object.

{\lltext In particular, choosing  $\sigma$ to be  the constant conductivity $\sigma=1$,  $\K=\overline B(1)$,
and $F$ to be the map $F_0$ given in (\ref{blow up})}, we obtain  the standard example} of the  
invisibility cloaking.
In dimensions $n\geq 3$ it shown in 2003 in \cite{GLU2,GLU3}
{\lltext that the Dirichlet-to-Neumann map corresponding to $H^1(\Omega)$ solutions
for the conductivity (\ref{eq: sing cond}) coincide with the
Dirichlet-to-Neumann map for  $\sigma=1$.}
In  2008, the analogous result was proven in   the two-dimensional case in
\cite{KSVW}. For cloaking results for the Helmholtz equation with frequency $k\not =0$
and for Maxwell's system in dimensions $n\geq 3$, see results in \cite{GKLU}.
We note also that John Ball \cite{Ball2} has used the push forward by the
analogous radial blow-up maps 
to study the discontinuity of the solutions of partial differential equations,
in particular the appearance of cavitation in the non-linear elasticity.

In the sequel we   consider cloaking results using measurements
given in Definition \ref{def: 1A}. As we have formulated the boundary measurements
in a new way, that is, in terms of the Dirichlet-to-Neumann forms $Q_\sigma$ associated to the class $W^{1,1}(\Omega)$,
we present the complete proof of the following proposition
in Subsection \ref{subsec: proof for counterexamples}.

\begin{theorem}\label{prop: invisibility} {\mmmtext (i) Let $\sigma\in 
L^\infty(B(2))$ be a scalar conductivity, $\sigma(x)\geq c_0>0$,
$\K\subset B(2)$ be a relatively compact open set with smooth boundary
and $F:\overline B(2)\setminus\{0\}\to
\overline B(2)
\setminus \K$ be a homeomorphism.
Assume that $F$ and $F^{-1}$ are 
$C^1$-smooth in $\overline B(2)\setminus\{0\}$ and $\overline B(2)
\setminus \K$, correspondingly and 
 $F|_{\p B(2)}=id.$
 {\ltext  Moreover, assume there is $C_0>0$ such that 
$\|DF^{-1}(x)\| \leq C_0$ for all $x\in \overline B(2)
\setminus \K$.}
Let  $\tilde \sigma$ be the conductivity defined in
(\ref{eq: sing cond}). 
Then  the boundary measurements for $\tilde \sigma$ and $\sigma$
coincide in the sense that $
Q_{\tilde \sigma}=Q_\sigma.
$
\smallskip

(ii) Let $\tilde \sigma$ be a cloaking conductivity {\lltext of the form (\ref{eq: sing cond})} obtained from
the transformation map $F$ and  the background conductivity $\sigma$ 
{\ltext where $F$ and $\sigma$ satisfy the  conditions in (i).} 
Then
\beq\label{eq: not in L^1}
\tr( {\tilde \sigma}) \not \in L^1(B(2)\setminus \K).
\eeq
}
\end{theorem}

The result (\ref{eq: not in L^1}) us optimal in the following sense.
{\mmmtext When $F$ is the map $F_0$ in (\ref{blow up}) and $\sigma=1$, the eigenvalues of the cloaking conductivity ${\tilde \sigma}$ in  $B(2)\setminus \overline B(1)$ behaves 
asymptotically as
$(|z|-1)$ and $(|z|-1)^{-1}$ as $|z|\to 1$. This cloaking conductivity has so strong degeneracy that (\ref{eq: not in L^1})
holds. On the other hand,
\beq\label{eq: in L^1_w}
\tr( {\tilde \sigma})  \in L^1_{weak}(B(2)).
\eeq
where $L^1_{weak}$ is the weak-$L^1$ space.
We note that in the case when $\sigma=1$,  $\det ({\tilde \sigma})$ is identically 1 in $B(2)\setminus \overline  B(1)$.}

The formula (\ref{eq: in L^1_w}) for the blow up map $F_0$ in (\ref{blow up}) 
and Theorem \ref{prop: invisibility} identify the {\it borderline of the invisibility}
for the trace of the conductivity:
Any cloaking conductivity $\tilde \sigma$ satisfies $\tr( {\tilde \sigma}) \not \in L^1(B(2))$
and there is an example of a cloaking conductivity for which $\tr( {\tilde \sigma})  \in L^1_{weak}(B(2)).$
Thus the borderline of invisibility is  the same as
the border between the space $L^1$ and the weak-$L^1$ space. 

\subsection{Counterexample 2: Illusion of a non-existent obstacle}\label{conter}

Next we consider new   counterexamples for the inverse problem which could be considered as creating
an illusion of a non-existing obstacle. 
The example is based on a radial shrinking map, that is, a mapping $ B(2)\setminus \overline B(1)\to B(2)\setminus \{0\}$.
The suitable maps are the inverse maps of the blow-up maps
$F_1:B(2)\setminus\{0\}\to B(2) \setminus \overline B(1)$ which
are constructed by  Iwaniec and Martin \cite{IM} and have 
the optimal smoothness. Alternative constructions for such blow up maps have also been proposed
by Kauhanen et al, see \cite{KKMOZ}.
Using the properties of these maps and
defining a conductivity $\sigma_1=(F_1^{-1})_* 1$ on $B(2)\setminus \{0\}$ 
we will later prove the following result.

\begin{theorem} \label{thm: hologram}
Let  $\gamma_1$ be a conductivity in $B(2)$ which
 is identically 1 in $B(2) \setminus \overline B(1)$ and zero in $B(1)$ 
and ${\mathcal A}:[1,\infty]\to [0,\infty]$ 
be any strictly increasing positive smooth function with $\A(1)=0$ which is
 sub-linear in the sense that 
\beq\label{IM cond 1}
\int_1^\infty \; \frac{{\mathcal A}(t)}{t^2}dt  < \infty.
\eeq
 Then there is a conductivity
 $\sigma_1\in \Sigma(B_2)$ satisfying $\det(\sigma_1) = 1$ and
 \beq
& &
 \int_{B(2)} \, \exp(
{\mathcal A}(\tr(\sigma_1)+\tr(\sigma_1^{-1}))) \,dm(z)  < \infty,
\label{IM3}
\eeq
such that $Q_{\sigma_1}=Q_{\gamma_1}$, i.e.,
the boundary measurements corresponding to $\sigma_1$ and $\gamma_1$ coincide.
\end{theorem}

We observe that for instance the function ${\mathcal A}_0(t)=t/(1+\log t)^{1+ \varepsilon}$
satisfies   (\ref {IM cond 1}) and for such weight function
 $\sigma_1\in L^1(B_2)$.
The proof of Theorem \ref{thm: hologram} is given in Subsection  \ref{subsec: proof for counterexamples}.

Note that $\gamma_1$ corresponds to the case when $B(1)$ is a perfect  insulator which
is surrounded with constant conductivity 1. 
Thus Theorem \ref{thm: hologram} can be interpreted by saying that there
is a relatively weakly degenerated conductivity satisfying integrability condition (\ref{IM3})
that creates in the boundary observations an illusion of an obstacle that does
not exists. Thus the conductivity can be considered as "electric hologram''.
As the obstacle can be considered as a "hole" in the domain, we
can say also that even the topology of the domain can not be detected. In other words,  Calder\'on's program 
to image the conductivity inside a domain using the boundary measurements
cannot work  within the class of degenerate  conductivities satisfying  (\ref{IM cond 1}) and (\ref{IM3}).

\subsection{Positive results for  Calder\'on's inverse problem} \label{uusilabel1}

{\lltext Let us formulate our first main  theorem which deals on 
inverse problems for an\-is\-o\-trop\-ic conductivities 
where both the trace and the determinant of the conductivity can be degenerate.} 

\begin{theorem}\label{main 2b both}
{\lltext Let $\Omega \subset \C$ be a bounded simply connected domain with smooth boundary. 
 Let $\sigma_1, \sigma_2\in \Sigma(\Omega)$ be matrix valued conductivities in  $\Omega$
 which 
 satisfy the integrability condition
 \ba
\int_\Omega \, \exp(p(\tr \sigma(z)+\tr(\sigma(z)^{-1})))\,dm(z)  < \infty
\ea}
for some $p>1$.
 Moreover, assume that 
\beq\label{eq: exp^3 integrability}
\hspace{-7mm}\int_\Omega {\mathcal E}(q \;\det \sigma_j(z))\,
dm(z)<\infty,\ \mbox{ for some } q>0,
\eeq 
where $\mathcal E(t)=\exp(\exp(\exp(t^{1/2}+t^{-1/2})))$ and
 $
Q_{\sigma_1} = Q_{\sigma_2}. 
$
Then there is a $W^{1,1}_{loc}$-homeomorphism $F:\Omega\to \Omega$ satisfying
 $F|_{\partial \Omega} = \mathrm{id}$ such that
\beq\label{push}
 \sigma_1
= F_*\, \sigma_2.
\eeq
 \end{theorem}

%


Equation (\ref{push}) can be stated as saying that $\sigma_1$ and $\sigma_2$ are the
same up to a change of coordinates, that is, the invariant
manifold structures corresponding to these conductivities are the same,
cf.\ \cite{LeU,LU}. 

{\lltext In the case when the conductivities are isotropic 
we can improve the result of Theorem \ref{main 2b both}.
The following theorem is our second main result for uniqueness of the inverse problem.}

\begin{theorem}\label{main 2b isotropic}
{\lltext Let $\Omega \subset \C$ be a bounded simply connected domain with smooth boundary. 
If $\sigma_1,\sigma_2\in \Sigma(\Omega)$ are isotropic conductivities, i.e.,
$\sigma_j(z)=\gamma_j(z)I$, $\gamma_j(z)\in   [0,\infty]$ satisfying 
\beq\label{eq: expexp integrability}
\int_\Omega \exp\bigl(\exp\bigl[ q(\gamma_j(z)+ \frac 1{\gamma_j(z)})\bigr]\bigr)\,dm(z)<\infty,\quad \mbox{ for some } q>0,
\eeq 
and $Q_{\sigma_1} = Q_{\sigma_2},$ then
 $\sigma_1=\sigma_2$.} \end{theorem}


Let us next consider {\lltext anisotropic} 
conductivities with bounded determinant but more degenerate ellipticity function $K_\sigma(z)$
defined in (\ref{uusi1}), and ask how far can we then generalize Theorem \ref{main 2b both}.
Motivated by the counterexample given in Theorem \ref{thm: hologram} we
 consider the following class: We say that $\sigma\in \Sigma(\Omega)$ has an exponentially degenerated 
anisotropy with a weight ${\mathcal A}$ 
and denote $\sigma\in \Sigma_{\mathcal A}=\Sigma_{\mathcal A}(\Omega)$ 
if $\sigma(z)\in \R^{2\times 2}$ for a.e.\ $z\in \Omega$ and 
 \begin{equation} \label{hotelli55}
 \int_\Omega \, \exp({\mathcal A}(\tr \sigma+\tr(\sigma^{-1})))\,dm(z)  < \infty.
\end{equation}
In view of Theorem \ref{thm: hologram}, for obtaining uniqueness for the inverse problem we need to consider  weights that are strictly increasing 
positive  {\ltext smooth} functions ${\mathcal A}:[1,\infty] \to [0,\infty]$, ${\mathcal A}(1) = 0$, with 
\beq\label{cond for A function}
\int_1^\infty \frac{{\mathcal A}(t)}{t^2}\,dt  = \infty \quad \quad  \mbox{ and }  \quad \quad  t {\mathcal A}'(t) \to \infty, \; \mbox{ as $t\to \infty$.}
\eeq
We say that $\A$ has almost linear growth if (\ref{cond for A function}) holds.
The point here is the first condition, that is, the divergence of the integral. The second condition is a technicality, which is  satisfied by all weights one encounters in practice \mmmtext {(which do not oscillate too much)}; the condition guarantees  that the Sobolev-gauge function $P(t)$ defined below in (\ref{hotelli57}) {\ltext is equivalent
to a convex function for large} $t$,
{\mmmtext see \cite[Lem.\ 20.5.4]{AIM}.} 

Note in particular that affine weights  $\A(t)=pt-p,$ $p>0$ satisfy the condition (\ref{cond for A function}). To develop uniqueness results for inverse problems within the class $\Sigma_{\mathcal A}$, the first questions we face are  to establish the right Sobolev-Orlicz regularity for the solutions $u$ of finite energy, $ A_\sigma[u] < \infty$, and solving the Dirichlet problem with given boundary values.

To start with this, we need the counterpart of the gauge $Q(t)$ defined at (\ref{gauge2}). In the case of a general weight  ${\mathcal A}$ we  {\lltext define} 
 \begin{equation} \label{hotelli57}
P(t) = \left\{\begin{array}{ll}
  t^2, & \quad \hbox{for}\quad   0   \leqslant t   <1,\\
\frac{t^2} {{\mathcal A}^{-1}(\log (t^2))} \,,  &\quad \hbox{for}\quad    t   \geqslant 1\end{array} \right.
\end{equation}
where ${\mathcal A}^{-1}$ is the inverse function of ${\mathcal A}$.
As an example, note that  if ${\mathcal A}$ is affine, ${\mathcal A}(t) = p\,t- p $ for some number $p >0$, then the condition (\ref{hotelli55}) takes us back to  the exponentially integrable distortion  
of Theorem \ref{main 2b both}, while  
$P(t) = t^2( 1 + \frac{1}{p} \log ^{+} (t^2)  )^{-1}$ is  equivalent to the gauge function $Q(t)$ used at (\ref{gauge2}).

 {\ltext The inequalities
(\ref{eq: Kari 2})  corresponding to the case
when   ${\mathcal A}$ is affine can be generalized for the following result
holding for general gauge  ${\mathcal A}$ satisfying  (\ref{cond for A function}).}

\begin{lemma}\label{calA11}   Suppose $u \in W^{1,1}_{loc}(\Omega)$ {\lltext and 
$\A$ satisfies  the almost linear growth condition (\ref{cond for A function}).}
Then
\ba
 \int_{\Omega}(P\bigl(|\nabla u|\bigr)  +P\bigl(|\sigma\nabla u|\bigr)) \,dm  
\leq 2  \int_{\Omega} e^{{\mathcal A}\left(\text{tr}\, \sigma+\text{tr}\,(\sigma^{-1})\right)} 
 \,dm(z)+ 2  \int_{\Omega} \nabla u \,\cdotp \sigma\nabla u \,dm\nonumber 
\ea
for every measurable function of symmetric matrices  $\sigma(z)\in \R^{2\times 2}$.
\end{lemma}

\noindent {\bf Proof.} We have in fact pointwise estimates. For these, 
note first that the conditions for $\A(t)$ imply that 
$P(t) \leq t^2$ for every $t\geq 0$. Hence, if $|\nabla u(z)|^2 \leq  
\exp {\mathcal A}\left(\tr \sigma(z)+\tr(\sigma^{-1}(z))\right) $ then
\beq\label{P exp equat pre1}
P\bigl(|\nabla u(z)|\bigr)   \leqslant \exp {\mathcal A}\left(\tr \sigma(z)+\tr(\sigma^{-1}(z))\right) .
\eeq
If, however, $|\nabla u(z)|^2 > \exp {\mathcal A}\left(\tr \sigma(z)+\tr(\sigma^{-1}(z))\right)$, then 
\beq\label{P exp equat pre2}
P\bigl(|\nabla u(z)|\bigr) = \frac{|\nabla u(z)|^2}{ {\mathcal A}^{-1}( \log |\nabla u(z)|^2)} \leq \frac{|\nabla u(z)|^2}{\tr(\sigma^{-1}(z))} \leq \nabla u(z) \,\cdotp \sigma(z)\nabla u(z).
\eeq
Thus at a.e.\ $z\in \Omega$ we have 
\beq\label{P exp equat}
P\bigl(|\nabla u(z)|\bigr)   \leq \exp {\mathcal A}\left(\tr \sigma(z)+\tr(\sigma^{-1}(z))\right) +  \nabla u(z) \,\cdotp \sigma(z)\nabla u(z).
\eeq
Similar arguments give pointwise bounds to $P\bigl(|\sigma(z)\nabla u(z)|\bigr)$.
Summing these estimates and integrating these pointwise estimates over $\Omega$ proves the claim.
\hfill \proofbox 
\medskip

In following,  we say that $u \in W^{1,1}_{loc}(\Omega)$  is in the Orlicz space
$W^{1,P}(\Omega)$ if
\ba
 \int_{\Omega}P\bigl(|\nabla u(z)|\bigr) \,dm(z)<\infty.  
\ea
\mattiHOX{Kerro, mitka ovat tekniset paatokalut ja uutuudet}

There are further important reasons  that make  the gauge $P(t)$  a natural and useful choice.  For instance, in   constructing  a minimizer for the energy $A_\sigma[u]$  we are faced with the problem of possible  equicontinuity  of Sobolev functions with uniformly bounded $A_\sigma[u]$. In view of Lemma \ref{calA11} this is reduced to describing those weight functions $\A(t)$ for which the condition $P\bigl(|\nabla u(z)|\bigr) \in L^1(\Omega) $ implies that the continuity modulus of $u$ can be estimated.  As we will
see later in (\ref{modulus of cont.}), this follows for 
weakly monotone functions $u$ (in particular for homeomorphisms),  as soon as  the divergence condition
  \begin{equation} \label{Ap171}
  \int_1^\infty \frac{P(t)}{ t^{3}} \, dt = \infty
  \end{equation}
 is satisfied, that is, $P(t)$ has almost quadratic growth.  
In fact, note that the divergence of the integral $\int_1^\infty \frac{ {\mathcal A}(t)}{ t^{2}} \, dt $ is equivalent to
\begin{equation} \label{hotelli56}
\int_1^\infty \frac{P(t)}{ t^{3}} \, dt = \frac{1}{2} \int_1^\infty \frac{ {\mathcal A}'(t)}{ t} \, dt = \frac{1}{2} \int_1^\infty \frac{ {\mathcal A}(t)}{ t^{2}} \, dt  = \infty
\end{equation}
where we have  used the  substitution $ {\mathcal A}(s) = \log (t^2)$. Thus the condition (\ref{cond for A function}) is directly connected to the smoothness properties of solutions of finite energy,  for  conductivities  
satisfying  (\ref{hotelli55}). 


 {\lltext 
We are now ready to formulate our third main theorem 
for uniqueness of inverse problem, which gives a sharp result for singular an\-is\-o\-trop\-ic conductivities
with a determinant bounded from above and below by positive constants.}

\begin{theorem}\label{main 2b trace}
{\lltext Let $\Omega \subset \C$ be a bounded simply connected domain with smooth boundary and
 ${\mathcal A}:[1,\infty)\to [0,\infty)$ be a strictly increasing smooth function satisfying  the almost linear growth condition (\ref{cond for A function}).  Let $\sigma_1, \sigma_2\in \Sigma(\Omega)$ be matrix valued conductivities in  $\Omega$
 which 
 satisfy the integrability condition
 \begin{equation} \label{hotelli55kopioitu}
\int_\Omega \, \exp({\mathcal A}(\tr \sigma(z)+\tr(\sigma(z)^{-1})))\,dm(z)  < \infty.
\end{equation}
 Moreover, suppose that $c_1\leq \det(\sigma_j(z))\leq c_2$, $z\in \Omega$, $j=1,2$ for
 some $c_1,c_2>0$
and
 $
Q_{\sigma_1} = Q_{\sigma_2}. 
$
Then there is a $W^{1,1}_{loc}$-homeomorphism $F:\Omega\to \Omega$ satisfying
 $F|_{\partial \Omega} = \mathrm{id}$ such that
\ba
 \sigma_1
= F_*\, \sigma_2.
\ea }
 \end{theorem}

We note that the determination of $\sigma$ from $Q_\sigma$ in {\lltext
Theorems \ref{main 2b both}, \ref{main 2b isotropic}, and \ref{main 2b trace} is  constructive
in the sense that one can write an algorithm which constructs $\sigma$ from $\Lambda_\sigma$.
For example,
for the non-degenerate scalar conductivities such a construction has been  numerically implemented in
\cite{AMPPS}.

Let us next discuss the borderline of the visibility somewhat formally.  Below we say  
that a conductivity is visible if there is  an algorithm 
which reconstructs the conductivity $\sigma$ from 
the boundary measurements $Q_\sigma$, possibly 
 up to  a change of coordinates. In other words,
 for visible conductivities one can use the boundary measurements to produce an image of the conductivity
in the interior of $\Omega$ in some deformed coordinates}.
For simplicity, let us consider conductivities  with $\det \sigma$ bounded from above and below.
Then, Theorems \ref{thm: hologram}  and \ref{main 2b trace} can be interpreted by saying that
the almost linear growth condition (\ref{cond for A function}) for the weight function $\A$ gives  
the {\it borderline of visibility} for the trace of the conductivity matrix:  
 If
$\A$
satisfies (\ref{cond for A function}), the  conductivities satisfying the integrability condition
(\ref{hotelli55kopioitu}) are visible.
However, if $\A$  does not satisfy 
 (\ref{cond for A function})   
 we can construct a conductivity in $\Omega$ satisfying the integrability condition (\ref{hotelli55kopioitu}) which appears as if an obstacle (which does not exist
in reality) would have included in the domain.

 Thus the borderline of the visibility  is between any spaces $\Sigma_{\A_{1}}$ and
 $\Sigma_{\A_{2}}$ where $\A_1$ satisfies condition (\ref{cond for A function}) and $\A_2$ does not
 satisfy it. Example of such gauge functions are $\A_1(t)=t(1+\log t)^{-1}$ and $\A_2(t)=t(1+\log t)^{-1-\e}$ with $\e>0$.

Summarizing the results, in terms of the trace of the conductivity,
we have identified the borderline of visible conductivities and the  borderline of  {\lltext 
invisibility} cloaking conductivities.
Moreover, these borderlines are not the same and between the visible and the {\lltext invisibility 
cloaking} conductivities
there are conductivities creating electric holograms.


\section{Proofs for the existence and uniqueness of the  solution of the direct problem and for the counterexamples.}

First we  show that under the conditions (\ref{hotelli55}) and (\ref{cond for A function}) 
  the Dirichlet problem for the conductivity equation admits 
a unique solution $u$ with finite energy $A_\sigma[u]$. 

  \subsection{The Dirichlet  problem} \label{sec: direct}

{\mmmtext In this section we prove  Theorem \ref{thm: Kari 3}. In fact, we prove
it in a more general
setting than it was stated.}

\begin{theorem}\label{thm: Kari 3BB}
{\mmmtext {\ltext  Let $\sigma\in \Sigma_{\mathcal A}(\Omega)$
where $\A$ satisfies the almost linear growth condition (\ref{IM cond 1}).
 Then, if $\, h\in H^{1/2}(\p \Omega)$ is such that $Q_\sigma[h]<\infty$} and $X=\{v\in W^{1,1}(\Omega);\ v|_{\p\Omega}=h\},$
 there is a unique $w\in X$ satisfies (\ref{eq: Matti-Kari 1A}).
Moreover, $w$ satisfies the conductivity equation
\beq\label{eq: Matti-Kari 1C}
\nabla \,\cdotp \sigma \nabla w=0\quad\hbox{in }\Omega
\eeq
in sense of distributions, and it has the regularity $w\in W^{1,P}(\Omega)$.} 
\end{theorem}

\noindent {\bf  Proof.}  For $N>0$, denote
$
\Omega_N=\{x\in \Omega;\ \|\sigma(x)\|+\|\sigma(x)^{-1}\|\leq N\}.
$
Let $w_n\in X$ be such {\ltext that
\ba
\lim_{n\to \infty} A_\sigma[w_n]=C_0=\inf\{A_\sigma[v]\ ;\ v\in X\}=Q_\sigma[h]<\infty
\ea
and $A_\sigma[w_n]<C_0+1$. 
{\mmmtext Then by Lemma \ref{calA11},
\beq\label{Ma 1}\hspace{-2mm}
 \int_\Omega P(|\nabla w_n(x)|)\,dm(x)+ \hspace{-1mm}\int_\Omega P(|\sigma (x)\nabla w_n(x)|)\,dm(x)\hspace{-.5mm}\leq\hspace{-.5mm} 
2(C_1+C_0+1)=C_2,
\eeq
where $C_1=\int_\Omega e^{\A(K(z))}dm(z)$.}
{\ltext By \cite[Lem.\ 20.5.3, 20.5.4]{AIM}, there is a convex, 
and unbounded
function $\Phi:[0,\infty)\to \R$ such that  $\Phi(t)\leq P(t)+c_0\leq 2\Phi(t)$ with some $c_0>0$ and
moreover, the function $t\mapsto \Phi(t^{5/8})$ is convex and increasing. This implies  that 
$P(t)\geq c_1t^{8/5}-c_2$ for some $c_1>0$, $c_2\in \R$. 
Thus (\ref{Ma 1}) yields
that for all $1<q\leq 8/5$}
\ba
\|\nabla w_n\|_{L^q(\Omega)}\leq C_3=C_3(q,C_0,C_1),\quad \hbox{for }n\in \Z_+.
\ea
Using the Poincare inequality in $L^q(\Omega)$ and that
 $(w_n-w_1)|_{\p \Omega}=0$, we see that $\|w_n-w_1\|_{L^q(\Omega)}\leq  C_4C_3$.
  Thus, there is $C_5$ such that $\|\nabla w_n\|_{W^{1,q}(\Omega)}<C_5$
for all $n$.}
By Banach-Alaloglu theorem this implies that  by restricting to a subsequence of
$(w_n)_{n=1}^\infty$, which we denote in sequel also by $w_n$,
such that $w_n\to w$ in $W^{1,q}(\Omega)$ as $n\to \infty$.
 {\ltext As $W^{1,q}(\Omega)$ embeds compactly
to $H^s(\Omega)$ for $s<2(1-q^{-1})$ we see that $\| w_n-w\|_{H^s(\Omega)}\to 0$
as $n\to \infty$ for all $s\in (\frac 12,\frac 34)$. Thus  $w_n|_{\p\Omega}\to w|_{\p\Omega}$ 
in $H^{s-1/2}(\p\Omega)$ as $n\to \infty$. This implies 
that $w|_{\p\Omega}=h$ and $w\in X$.}
Moreover,
for any $N>0$
\ba
\frac 1N \int_{\Omega_N}  |\nabla w_n(x)|^2\, dm(x)\leq 
\int_{\Omega_N}  \nabla w_n(x) \,\cdotp \sigma (x) \nabla w_n(x) dm(x)\leq C_0+1.
\ea
This implies that $\nabla w_n|_{\Omega_N}$ are uniformly bounded in $L^2(\Omega_N)^2$.
Thus by restricting to a subsequence, we can assume that 
$\nabla w_n|_{\Omega_N}$ converges weakly in $L^2(\Omega_N)^2$ as $n\to \infty$.
Clearly, the weak limit must be $\nabla w|_{\Omega_N}$. 
Since the norm $V\mapsto (\int_{\Omega_N} V \,\cdotp \sigma V\, dm)^{1/2}$ in $L^2(\Omega_N)^2$ is weakly lower semicontinuous, 
we see that
 \ba
\int_{\Omega_N}  \nabla w(x) \,\cdotp \sigma (x) \nabla w(x) dm(x)\leq
\liminf_{n\to \infty}
\int_{\Omega_N}  \nabla w_n(x) \,\cdotp \sigma (x) \nabla w_n(x) dm(x)\leq C_0.
\ea
As this holds for all $N$, we see by applying the monotone convergence theorem
as $N\to \infty$ that (\ref{eq: Matti-Kari 1A}) holds. Thus $w$ is a minimizer
of $A_\sigma$ in $X$.

{\lltext By the above,
$\sigma\nabla w_n\to \sigma\nabla  w$ weakly in $L^2(\Omega_N)$ as $n\to \infty$ for all $N$.
{\ltext As noted above there is a convex function $\Phi:[0,\infty)\to \R$ such that  $\Phi(t)\leq P(t)+c_0\leq 2\Phi(t)$,
$c_0>0$
and that $\Phi(t)$ is increasing for large values of $t$.} 
Thus it follows from the  semicontinuity results for integral operators,
\cite[Thm.\ 13.1.2]{ABM}, Lebesgue's monotone convergence theorem,
and (\ref{Ma 1}) {\lltext that
\ba
& &\hspace{-7mm} \int_\Omega(\Phi(|\nabla w|)+ \Phi(|\sigma\nabla w|)\,dm(x)
\leq 
\lim_{N\to \infty} \liminf_{n\to \infty} \int_{\Omega_N} (\Phi(|\nabla w_n|)+\Phi(|\sigma\nabla w_n|))\,dm
 \\
   & &\nonumber \hspace{-4mm} \leq 
\lim_{N\to \infty} \liminf_{n\to \infty} \int_{\Omega_N} (P(|\nabla w_n|)+P(|\sigma\nabla w_n|))\,dm
+2c_0|\Omega|\leq 
C_2+2c_0|\Omega|.
\ea}
{\ltext It follows from the above
and the inequality
$P(t)\geq c_1t^{8/5}-c_2$ 
that} $\sigma(x) \nabla w (x) \in L^1(\Omega)$.} Consider next $\phi\in C^\infty_0(\Omega)$.
As $w+t\phi\in X$, $t\in\R$
and as $w$ is a minimizer of $A_\sigma$ in $X$ it follows  that
\ba
\left. \frac d{d t} A_\sigma[w+t\phi]\right|_{t=0}=2\int_{\Omega} 
 \nabla \phi(x) \,\cdotp \sigma (x) \nabla w(x) dm(x)=0.
 \ea
This shows that the conductivity equation (\ref{eq: Matti-Kari 1C}) is valid 
in the sense of distributions.

Next, assume that $w$ and $\tilde w$ are both minimizers
of $A_\sigma$ in $X$.
Using the convexity of $A_\sigma$ we see that then the second derivative of $t\mapsto 
A_\sigma[tw+(1-t)\tilde w]$ vanishes at $t=0$.
This implies that $\nabla (w-\tilde w)=0$ for a.e.\ $x\in \Omega$. As $w$ and
$\tilde w $ coincide at the boundary, this yields that $w=\tilde w$ and
thus the minimizer is unique.
\hfill \proofbox 
\medskip

The fact that the minimizer $w$ is continuous will be proven
in the next subsection. 

\subsection{The Beltrami equation}\label{subsec: Beltrami}

It is natural to ask if the minimizer $w$ in (\ref{eq: Matti-Kari 1A}) is the only solution of finite $\sigma$-energy 
$A_\sigma[w]$ to the boundary value problem
\beq\label{uusi3}
\nabla \,\cdotp \sigma \nabla w&=&0\quad\hbox{in }\Omega,\\
\nonumber w|_{\p\Omega}&=&h.
\eeq
It turns out the this is the case {\lltext and to prove this we 
introduce one of the basic tools in this work, the Beltrami differential equation.}

For this end, recall the Hodge-star operator $*$ which in two dimensions is just the  rotation
\ba
* = \left(\begin{array}{cc}  0& -1 \\
 1& 0\\ \end{array}\right).
\ea
Note that $\nabla \,\cdotp (*\nabla u)=w$
for all $w\in W^{1,1}(\Omega)$ and recall that $\Omega\subset \C$ is simply connected.
{\ltext If $\sigma(x)=[\sigma^{jk}(x)]_{j,k=1}^2\in \Sigma_\A(\Omega)$,
where $\A$ satisfies (\ref{IM cond 1}),
and if $u\in W^{1,1}(\Omega)$ is a distributional solution to the
conductivity equation
\beq\label{eq:  BB}
\nabla\cdot\sigma(x)\nabla u(x) = 0, 
\eeq
then  by Lemma \ref{calA11}
we have $P(\nabla u),P(\sigma\nabla u)\in L^1(\Omega)$
and thus in particular $\sigma\nabla u\in L^1(\Omega)$.  By (\ref{eq:  BB}) and 
the} Poincare lemma there is a function $v\in W^{1,1}(\Omega)$ such that 
\beq\label{eq: Matti conjugate}
\nabla v=* \, \sigma(x) \nabla u(x).
\eeq
Then
\beq\label{eq: CC}
\nabla  \,\cdotp \sigma^{*}(x)\nabla v=0\quad\hbox{in }\Omega, \quad \quad   \sigma^{*}(x) = * \, \sigma(x)^{-1} *^t.
\eeq
{\ltext In particular, the above shows that $u,v\in W^{1,P}(\Omega)$.}
Moreover, an explicit calculation, see e.g. \cite[formula (16.20)]{AIM}, reveals that  the function $f = u + iv$ satisfies
\beq\label{eq: anisotropic Beltrami AA}
\partial_{\overline z} f = \mu \partial_z f + 
\nu  \,\overline{\partial_z f},
\eeq
where 
\beq\label{35+}
\mu = \frac{\sigma^{22}- \sigma^{11}-2i \sigma^{12}}{1 + \tr(\sigma) + \det (\sigma)}, \quad
\nu = \frac{1-\det (\sigma)}{1 + \tr(\sigma) + \det (\sigma)}, 
\eeq
 {\lltext and $\partial_{\overline z}=\frac 12(\p_{x_1}+i\p_{x_2})$ with
$\partial_{z}=\frac 12(\p_{x_1}-i\p_{x_2})$.}
Summarizing, for $\sigma\in \Sigma_\A(\Omega)$ any distributional
solution $u\in W^{1,1}(\Omega)$ of (\ref{eq:  BB}) is a real part 
of the solution $f$ of (\ref{eq: anisotropic Beltrami AA}). 
Conversely, the real part of any solution $f\in W^{1,1}(\Omega)$ of (\ref{eq: anisotropic Beltrami AA}) satisfies (\ref{eq: BB}) while the imaginary part is a solution to  (\ref{eq: CC}) and 
 as $\sigma\in \Sigma_\A(\Omega)$, 
{\ltext (\ref{eq: BB})-(\ref{eq: CC}) 
 and
Lemma \ref{calA11} yield that
  $u,v\in W^{1,P}( \Omega)$, and hence  $f\in W^{1,P}( \Omega)$. }

Furthermore,  the ellipticity bound of $\sigma(z)$ is closely related to the distortion of the mapping $f$. 
Indeed, 
in the case when $\sigma(z_0)=\diag(\la_1,\la_2)$, a direct computation
shows that 
%
\beq\label{distor2}
K_\sigma(z_0) = K_{\mu,\nu}(z_0),\quad\hbox{where }K_{\mu,\nu}(z)=  \frac{1+ |\mu(z)| + |\nu(z)|}{1-( |\mu(z)| + |\nu(z)|)}
\eeq
and $K_\sigma(z)$ is the ellipticity of $\sigma(z)$ defined in (\ref{uusi1}). 
Using the chain rule for the complex derivatives, which can be
written as
\beq\label{chain rules}
\p(v\circ F)&=&
(\p v)\circ F\cdotp\p F+(\overline \p v)\circ F\cdotp\overline{\overline \p  F},\\
\overline \p(v\circ F)
&=&
(\p v)\circ F\cdotp\overline \p F+(\overline \p v)\circ F\cdotp\overline{ \p  F},
\eeq
we see that $|\mu(z)|$ and   $|\nu(z)|$ do not change in an orthogonal rotation
of the coordinate axis, $z\mapsto \alpha z$ where $\alpha\in \C$, $|\alpha|=1$.
As for any $z_0\in \Omega$ there exists an orthogonal
rotation of the coordinate axis so that  matrix  
$\sigma(z_0)$ is diagonal in the rotated coordinates, we see that
the identity (\ref{distor2}) holds for all $z_0\in \Omega$.


The equation (\ref{eq: anisotropic Beltrami AA}) is also
 equivalent to the Beltrami equation
\beq\label{eq: pseudo Beltrami2}
\dbar f(z)=\tilde \mu(z)\,\p f(z)\quad\hbox{in }\Omega,
\eeq
{\lltext with the Beltrami coefficient 
\beq\label{eq: modified anisotropic Beltrami A1}
\tilde \mu(z)=\left\{\begin{array}{cl}
\mu(z)+\nu(z)\p_z f(x)\left(\overline{\p_z f(x)}\right)^{-1}&\hbox{if }\p_z f(x)\not =0,\\
\mu(z)&\hbox{if }\p_z f(x)=0,\end{array}\right.
\eeq
 satisfying} $|\tilde \mu(z)| \leq  |\mu(z)| + |\nu(z)|$ pointwise. 
 We define the distortion of $f$ at $z$ be 
\beq\label{distor2 mod}
K(z,f) := K_{\tilde \mu}(z)=\frac{1+ |\tilde \mu(z)|}{1- |\tilde \mu(z)|}\leq K_\sigma(z),\quad z\in \Omega.
\eeq
%
%
{\mmmtext 
Below will also use the notation $K(z,f)=K_f(z)$.} 


 In the sequel we will use frequently these different interpretations of the Beltrami equation. Note that $K(z,f)={(1+|\tilde 
\mu(z)|)}/({1-|\tilde \mu(z)|})$ so that
$K(z,f)=(|\p  f|+|\dbar  f |)/(|\p  f |-|\dbar  f|)$.
As
$\|D f \|^2=(|\p  f |+|\dbar  f |)^2$ and $J(z, f)=|\p  f |^2-|\dbar  f |^2$,
this yields the distortion equality, see e.g.\ \cite[formula (20.3)]{AIM},
\beq\label{eq: distortion ineq.}
\|D f(z)\|^2= K(z, f)J(z, f),\quad\hbox{for a.e.\ }z\in \Omega.
\eeq

We will use extensively the fact that 
if {\ltext a homeomorphism} $F:\Omega\to \Omega'$, {\lltext $F\in W^{1,1}(\Omega)$
is a finite distortion mapping}   with the distortion
 $K_F\in L^1(\Omega)$ then
by \cite{HKO} or \cite[Thm.\ 21.1.4]{AIM} 
  the inverse function $H=F^{-1}:\Omega'\to \Omega$ 
 is  in $W^{1,2}(\Omega')$  and its derivative $DH$ satisfies
\beq\label{eq: W12}
\|DH\|_{L^2(\Omega')}\leq 2\|K_F\|_{L^1(\Omega)}.
\eeq

We will also need few basic notions, see \cite{AIM}, from the theory of Beltrami equations. As the coefficients $\mu, \nu$ are defined only in the bounded domain $\Omega$, outside   $\Omega$  we set  $\mu(z) =  \nu(z) = 0$ and $\sigma(z) =1$,  and consider global solutions to (\ref{eq: anisotropic Beltrami AA}) in $\C$. 
{\ltext In particular, we consider the case when $\Omega$ is the unit disc $\D=B(1)$.} We say that 
a solution  $f\in W^{1,1}_{loc}(\C)$ of the equation  (\ref{eq: anisotropic Beltrami AA}) in $z\in \C$ is a 
{\it principal solution} if 
\ba
& &1. \quad  f:\C\to \C \; \mbox{ is a homeomorphism of } \C \mbox{ \;and }\\
& &2. \quad \; 
f(z)=z+ {\mathcal O}(1/z) \quad \quad  \mbox{as \;} z\to \infty.
\ea
The existence principal solutions is a fundamental fact that holds true  in quite wide generality.
Further, with the principal solution one can classify all solutions, of sufficient regularity, to the Beltrami equation.
{\mmmtext These facts are summarized in the following version of Stoilow's factorization theorem, 
for which proof we cite to  \cite[Thm.\ 20.5.2]{AIM}.}

\begin{theorem} \label{Stoilow}
Suppose $\mu(z)$ is  supported in the unit disk $\D$, {\lltext $|\mu(z)|<1$ a.e.\ and  
$$  \int_\D \, \exp({\mathcal A}(K_{\mu}(z)))\,dm(z)  < \infty, \quad K_{\mu}(z)=\frac {1+|\mu(z)|}{1-|\mu(z)|}
$$ 
where $\A$ satisfies the almost linear growth condition}  (\ref{cond for A function}). 
Then the equation
\beq\label{Beltrami.principal}
& &\dbar \principalPhi(z)= \mu(z)\,\p \principalPhi(z),\quad z\in \C,\\
\label{Beltrami.principal2}
 & &\principalPhi(z)=z+ {\mathcal O}(1/z) \quad \quad  \mbox{as \;} z\to \infty
\eeq 
has a unique solution in $\principalPhi\in W^{1,1}_{loc}(\C)$. The solution
$\principalPhi:\C\to \C$ is a homeomorphism and satisfies $\principalPhi\in W^{1,P}_{loc}(\C)$. 
Moreover, when $\Omega_1\subset \C$ is open, every solution of the equation
\beq\label{Beltrami24}
\dbar f(z)= \mu(z)\,\p f(z),\quad z\in \Omega_1,
\eeq 
with  the regularity  $f \in W^{1,P}_{loc}(\Omega_1)$,
can be written as $f = H \, \circ \, \principalPhi$, where $\principalPhi$ is the solution to 
(\ref{Beltrami.principal})-(\ref{Beltrami.principal2}) and $H$ is holomorphic function in $\Omega_1' = \principalPhi(\Omega_1)$.

\end{theorem}

{\lltext Below we apply  this results with the Poincare lemma to analyze the solutions of conductivity equation in the simply connected domain
 $\Omega$.}


\begin{corollary} \label{laitos}
{\mmmtext  Let $\sigma\in \Sigma_\A(\Omega)$ where
$\A$ satisfies
 (\ref{cond for A function}) and} $u \in W^{1,1}_{loc}(\Omega)$ satisfy
\beq\label{laitos3} \nabla \,\cdotp \sigma \nabla w=0\quad\hbox{in }\Omega \quad 
\hbox{and}\quad  
\int_\Omega  \nabla u(x)\cdot  \sigma(x)\nabla u(x)\,dm(x) < \infty.
\eeq
Then $u = w \circ \principalPhi$, where $\principalPhi:\C \to \C$ is a homeomorphism, $\principalPhi\in  W^{1,P}_{loc}(\C)$, and  $w$ is harmonic in the domain $\Omega'=\principalPhi(\Omega)$. In particular, $u:\Omega\to \R$ is continuous.
\end{corollary}
{\bf Proof.} Let  $v \in W^{1,1}_{loc}(\Omega)$ be the conjugate function described in (\ref{eq: Matti conjugate}), and set $f=u+iv$. Then by Lemma \ref{calA11} we have $f\in W^{1,P}(\Omega)$
and Theorem \ref{Stoilow} yields that
$f =  H \, \circ \, \principalPhi$, where  $\principalPhi:\C\to \C$  a homeomorphism
with $\principalPhi\in  W^{1,P}_{loc}(\C)$
and $H$ is holomorphic in $\principalPhi(\Omega)$. Thus the real part $u = (\Re H) \, \circ \, \principalPhi$  has the required factorization.
\hfill \proofbox 

%



%


%

Theorem \ref{thm: Kari 3BB} and Corollary \ref{laitos} yield  Theorem \ref{thm: Kari 3}.

\subsection{Invariance of  Dirichlet-to-Neumann form in coordinate transformations.}
\label{sec: Invariance}

In this section, we assume that $\sigma\in \Sigma_{\mathcal A}(\Omega)$, where $\mathcal A$ satisfies
(\ref{cond for A function}).
{\mmmtext We say that $F:\Omega\to \Omega'$ satisfies
the condition ${\mathcal N}$ if for any measurable 
set $E \subset \Omega$ we have  $| E|=0 \Rightarrow |F(E)| = 0$.
Also, we say that
$F$ satisfies
the condition ${\mathcal N}^{-1}$ if for any measurable 
set $E \subset \Omega$ we have  $| F(E)|=0 \Rightarrow |E| = 0$.}

{\mmmtext 
Let  $\sigma\in \Sigma_{\mathcal A}(\C)$ be such that 
 {\ltext  $\sigma$
is constant 1 in $\C\setminus\Omega$}. Let
\beq\label{eq: hat mu}
\hat \mu(z)  = \frac{\sigma^{11}(z)- \sigma^{22}(z) + 2i\sigma^{12}(z)}{\sigma^{11}(z) + \sigma^{22}(z) +
2\sqrt {\det \sigma(z)}}
\eeq
be the Beltrami coefficient associated to the isothermal coordinates corresponding
to $\sigma$, see e.g.\ \cite{S}, \cite[Thm.\ 10.1.1]{AIM}. A direct computation shows
that
$K_{\hat\mu}(z)=K_{\tilde\sigma}(z)$ 
and thus $\exp(\A(K_{\hat\mu}))\in L^1_{loc}(\C)$
and
by Theorem \ref{Stoilow}, there exists a homeomorphism $F:\C\to \C$
satisfying the equation (\ref{Beltrami.principal})-(\ref{Beltrami.principal2}) with the Beltrami coefficient $\hat \mu$
such that $F\in W^{1,P}_{loc}(\C)$.
Due to the choice of $\hat \mu$, the conductivity $F_*\sigma$ is isotropic, see e.g.\  
\cite{S}, \cite[Thm.\ 10.1.1]{AIM}. 
Let us next consider the properties of the map $F$. First, as $\exp({\mathcal A}(K_{\hat \mu}))\in L^1_{loc}(\C)$, it follows from
\cite{KKMOZ} that the function  $F$ satisfies the condition $\mathcal N$.
%
Moreover, the fact that  $K_F=K_{\hat \mu}\in L^1_{loc}(\C)$ implies by (\ref {eq: W12}) that it inverse $H=F^{-1}$ 
 is  in $W^{1,2}_{loc}(\C)$. This yields by  \cite[Thm.\ 3.3.7]{AIM} that $F^{-1}$ satisfies the condition $\mathcal N$.
In particular, the above yields that both $F$ and $F^{-1}$ are in $W^{1,P}_{loc}(\C)$.

The following lemma formulates the invariance of the Dirichlet-to-Neumann forms in
the diffeomorphisms satisfying the above properties.

\begin{lemma}\label{lem: change of coord.}
{\ltext Assume that $\Omega,\tilde \Omega\subset \C$ are bounded, simply
connected domains with smooth boundaries and that  $\sigma\in \Sigma_{\mathcal A}(\Omega)$ and $\tilde\sigma\in \Sigma_{\mathcal A}(\tilde \Omega)$
where $\mathcal A$ satisfies
(\ref{cond for A function}).
Let $F:\Omega\to \tilde \Omega$ be a homeomorphism so that
$F$ and $F^{-1}$ are $W^{1,P}$-smooth and $F$ satisfies conditions $\mathcal N$ and $\mathcal N^{-1}$. Assume
that $F$ and $F^{-1}$ are  $C^1$ smooth near the boundary
and assume that $\rho=F|_{\p \Omega}$ is $C^2$-smooth.
Also, suppose
  $\tilde \sigma=F_*\sigma$. 
Then $Q_{\tilde \sigma}[\tilde h]=Q_\sigma[\tilde h\circ \rho]$ for all $
\tilde h\in H^{1/2}(\p \tilde\Omega)$.}
\end{lemma}

{\bf Proof.} 
As 
$F$ has the properties ${\mathcal N}$  and ${\mathcal N}^{-1}$ 
we have {\mmmtext the area formula}
\beq\label{eq: apu 1}
\int_{\tilde \Omega} H(y)\,dm(y)= \int_\Omega H(F(x))J(x,F)\,dm(x)
\eeq
for all simple functions $H:\tilde\Omega\to \C$, {\lltext where
$J(x,F)$ is the Jacobian determinant of $F$ at $x$}. Thus (\ref{eq: apu 1}) holds
for all $H\in L^1 (\tilde \Omega)$. 

{\ltext Let $\tilde h\in H^{1/2}(\p \tilde \Omega)$ and assume that   $Q_{\tilde \sigma}[\tilde h]<\infty.$}
Let $\tilde u:\tilde \Omega\to \R$ be
the unique  minimizer of $A_{\tilde \sigma}[v]$ in
$\tilde X=\{ \tilde v\in W^{1,1}(\tilde \Omega);\
 \tilde  v|_{\p\tilde \Omega}=\tilde h\}$. 
Then $\tilde u$  is 
 the solution of the conductivity equation
\begin{equation}\label{johty 2B}
\nabla\cdot \tilde \sigma\nabla \tilde u = 0,\quad \tilde u|_{\p\tilde\Omega}
=\tilde h.
\end{equation}
We define   $h=\tilde h\circ F|_{\p \Omega}$ and
 $u=\tilde u\circ F:\Omega\to \C$.

{\ltext 
By Corollary \ref{laitos},  $\tilde u$ can be written   
in the form $\tilde u=\tilde w\circ \tilde G$ where $\tilde w$ is harmonic and
 $\tilde G\in W^{1,1}_{loc}(\C)$ is a homeomorphism $\tilde G:\C\to \C$.  
 
 By Gehring-Lehto theorem, see   \cite[Cor.\ 3.3.3]{AIM},
a homeomorphism $F\in W^{1,1}_{loc}(\Omega) $ is differentiable almost everywhere in $ \Omega$,
say in the set $\Omega\setminus A$, where $A$
has Lebesgue measure zero. 
Similar arguments for $\tilde G$ show that {\lltext  $\tilde G$ and the solution $\tilde u$ 
are} differentiable 
almost everywhere, say in the set $\tilde\Omega\setminus A'$, where $A'$
has Lebesgue measure zero.} 

Since $F$ has 
the property ${\mathcal N}^{-1}$, we see that $A''=A'\cup F^{-1}(A')\subset \Omega$ has
measure zero, and for $x\in \Omega\setminus A''$ the chain rule gives
\beq\label{eq: chain rule}
Du(x)=(D\tilde u)(F(x))\,  \,\cdotp DF(x).
\eeq
Note that the facts that $F$ is a map with an exponentially integrable distortion and 
that $\tilde u$ is 
a real part of a map with an exponentially integrable distortion, do
not generally imply, at least  according to the knowledge of the authors, that
their
composition $u$ is in $W^{1,1}_{loc}(\Omega)$.
To overcome this problem, we define for $m>1$
\ba
& &\tilde \Omega_m=\{y\in \tilde  \Omega;\ \|DF^{-1}(y)\|+
\|DF(F^{-1}(y))\|+ 
\|\tilde \sigma (y)\|+ |\nabla \tilde u(y)|< m \}
\ea
and $ \Omega_m=F^{-1}(\tilde \Omega_m)$.
Then $\nabla u|_{\Omega_m}\in L^2(\Omega_m)$ and $\|\sigma\|<m^3$  in $\Omega_m$.

Now for any $m>0$
\beq\label{eq: apu 1A b}
\int_{\tilde \Omega_m} \nabla \tilde  u(y) \,\cdotp \tilde \sigma (y)\nabla \tilde u(y)\,dm(y)
\leq A_{\tilde \sigma}[\tilde u]<\infty.
\eeq
Due to the definition of $\tilde \sigma =F_*\sigma$, we see by using formulae
(\ref {eq: apu 1}) and (\ref{eq: chain rule}) that
\beq\label{eq: apu 2b}
\int_{ \Omega_m} \nabla u(x) \,\cdotp \sigma (x)\nabla u(x)\,dm(x)= \int_{\tilde \Omega_m} \nabla \tilde  u(y) \,\cdotp \tilde \sigma (y)\nabla \tilde u(y)\,dm(y).
 \eeq
 Letting $m\to \infty$ and using monotone convergence theorem, we see that
 \beq\label{eq: apu 2c}
\int_{ \Omega} \nabla u(x) \,\cdotp \sigma (x)\nabla u(x)dm(x)=
\int_{\tilde \Omega} \nabla \tilde  u(y) \,\cdotp \tilde \sigma (y)\nabla \tilde u(y)dm(y)
=A_{\tilde \sigma}[\tilde u]<\infty.\hspace{3mm}
 \eeq
 By Lemma \ref{calA11} this implies that $u\in W^{1,P}(\Omega)\subset W^{1,1}(\Omega)$.

 Clearly, {\ltext as $\rho=F|_{\p \Omega}$ is $C^2$-smooth 
 $h:=\tilde h\circ F\in H^{1/2}(\p \Omega)$ and
 $u|_{\p\Omega}=h$.}  
   Thus $u\in X=\{ w\in W^{1,1}(\Omega);\
 w|_{\p\Omega}=h\}$. 
Since  $\tilde u$ is a minimizer of $A_{\tilde \sigma}$ in $\tilde X$, and $u$ satisfies
$
 A_{\sigma}[u]\leq A_{\tilde \sigma}[\tilde u]=Q_{\tilde \sigma}(\tilde h)
 $ 
 and see that
\ba
 Q_{\sigma}[h]\leq Q_{\tilde \sigma}[\tilde h].
\ea 
Changing roles of $\tilde \sigma$ and $\sigma$ we obtain an opposite inequality,
and prove the claim.

 \hfill \proofbox 
\medskip

{\ltext
In particular, if   $\sigma\in \Sigma_{\mathcal A}(\Omega)$, $\tilde\sigma\in \Sigma_{\mathcal A}(\tilde \Omega)$ 
and $F$ are as in Lemma \ref{lem: change of coord.} and in addition to that,
$\sigma$ and $\tilde \sigma$ are bounded near $\p \Omega$ and $\p\tilde\Omega$,
respectively and $\rho=F|_{\p \Omega}:\p \Omega\to\p \tilde \Omega$ is $C^2$-smooth, 
then the quadratic forms $Q_\sigma$ and $Q_{\tilde \sigma}$
can be written in terms of the Dirichlet-to-Neumann maps 
$\Lambda_\sigma:H^{1/2}(\p \Omega)\to H^{-1/2}(\p \Omega)$ and 
$\Lambda_{\tilde \sigma}:H^{1/2}(\p \tilde \Omega)\to H^{-1/2}(\p \tilde \Omega)$
as in formula (\ref{eq: Q forms}). Then, Lemma \ref{lem: change of coord.} implies that 
\beq\label{eq: push forward of Lambda}
\Lambda_{\tilde \sigma }=\rho _*  \Lambda_\sigma,
\eeq
where  $\rho _*\Lambda_\sigma$ is the push forward of $\Lambda_\sigma$ in $\rho $ 
defined by $ (\rho _*\Lambda_{\sigma })(\tilde h)= 
j\,\cdotp[(\Lambda_\sigma( \tilde h\circ \rho ))\circ \rho^{-1}] $
for $\tilde h\in H^{1/2}(\p \tilde \Om)$, where $j(z)$ is the Jacobian of the map
$\rho ^{-1}:\p \tilde \Omega\to\p  \Omega$.}

\subsection{Counterexamples revisited}
\label{subsec: proof for counterexamples}

In this section we give the proofs of the claims stated in Subsection \ref{sec: Counterexamples}.
We start by proving Theorem \ref{prop: invisibility}. Since the used singular
change of variables in integration is a tricky business  we present the 
arguments in detail.

\noindent {\bf Proof} (of Thm.\ \ref{prop: invisibility}). {\mmmtext (i) 
Our aim is first  to show that we have $ Q_\sigma[h]\leq Q_{\tilde \sigma}[h]$
and then to prove the opposite inequality. The  proofs of these
inequalities are based on different techniques due to the fact that 
$\tilde \sigma$ is not even in $L^1(B(2))$.


 Let $0< r<2$  and $\K(r)=\K\cup F(\overline B(r))$.
Moreover, let $\tilde \sigma_r$ be 
a conductivity that coincide
with ${\tilde \sigma}$ in $B(2)\setminus \K(r)$ and is 
$0$ in  $\K(r)$. Similarly, let 
$\sigma_r$ be 
a conductivity that coincide
with $\sigma$ in $B(2)\setminus \overline B(r)$ and is 
$0$ in  $\overline B(r)$. 
For these conductivities we define the quadratic forms $A^r:W^{1,1}(B(2))\to \R_+\cup\{0,\infty\}$ and $\tilde A^r:W^{1,1}(B(2))\to \R_+\cup\{0,\infty\}$,
\ba
A^r [v]=\int_{B(2)\setminus \overline B(r)} \nabla v\, \,\cdotp\sigma\nabla v\,dm(x),
\quad \tilde A^r [v]=\int_{B(2)\setminus \K(r)} \nabla v\, \,\cdotp\tilde \sigma\nabla v\,dm(x).
\ea
If we minimize  $\tilde A^r [v]$ over $v\in W^{1,1}(B(2))$
with $v|_{\p B(2)}=h$, we see that minimizers exist and that the restriction of
any minimizer to $B(2)\setminus \overline \K(r)$ is the function $\tilde u_r\in W^{1,2}(B(2)\setminus \K(r))$
satisfying
\ba
 \nabla\,\cdotp\tilde \sigma\nabla \tilde u_r =0\quad\hbox{in }B(2)\setminus \K(r),\quad 
\tilde u_r |_{\p B(2)}=h,\quad \nu
\,\cdotp  \tilde\sigma\nabla \tilde u_r|_{\p \K(r)}=0.
\ea
Analogous equations hold for the minimizer $u^r$ of $  A^r $.
As 
$\sigma$ in
$\overline B(2)\setminus B(r)$ 
and $\tilde \sigma$ in $\overline {B(2)\setminus \K(r)}$
are bounded from above and below by positive constants,
we using the change of variables and the chain rule that
\beq\label{Q are the same} 
Q_{\sigma_r}[h]=Q_{\tilde \sigma_r}[h],\quad\hbox{for }h\in H^{1/2}(\p B(2)).
\eeq
As $\sigma(x)\geq { \sigma}_r(x)$ and $\tilde \sigma(x)\geq {\tilde \sigma}_r(x)$ for all $x\in B(2)$,
\beq\label{ineq: a}Q_\sigma[h]\geq Q_{\sigma_r}[h],\quad
Q_{\tilde \sigma}[h]\geq Q_{\tilde \sigma_r}[h].
\eeq
%

Let us consider the  minimization problem
(\ref{def: DN map}) for $\sigma$. It is solved by the unique minimizer $u\in W^{1,1}(B(2))$ satisfying
\ba
\nabla\,\cdotp \sigma\nabla u=0\quad\hbox{in }B(2),\quad u|_{\p B(2)}=h.
\ea
As $\sigma,\sigma^{-1}\in L^\infty(B(2))$ we have  $u\in W^{1,2}(B(2))$ and  
Morrey's theorem \cite{Morrey}
yields
that the solution $u$ is $C^{0,\alpha}$-smooth in the open ball $B(2)$ for some $\alpha>0$.
{\ltext Thus $u|_{B(R)}$ is in the Royden algebra  ${\mathcal R}(B(R))=C(B(R))\cap L^\infty(B(R))
\cap W^{1,2}(B(R))$ for all $R<2$.}

By e.g.\ \cite[p. 443]{IM}, for any $0<R<2$ the p-capacity of the disc $B(r)$
in $B(R)$ goes to zero as $r\to 0$ for all $p>1$.
Using this, and that $u\in W^{1,2}(B(2))\subset L^q(B(2))$ for $q<\infty$,
we see that (cf.\ \cite{KSVW} for  explicit estimates in the case when $\sigma=1$) 
\ba
\lim_{r\to 0} Q_{\sigma_r}[h]=Q_\sigma[h],
\ea
that is, the effect of an insulating disc of radius $r$ in the boundary measurements
vanishes as $r\to 0$. This and
the inequalities (\ref{Q are the same}) and (\ref{ineq: a})  yield $Q_{{\tilde \sigma}}[h]\geq Q_{\sigma}[h]$. Next we consider the opposite
inequality.


{\ltext 
Let  $\tilde u=u\circ F^{-1}$ in $B(2)\setminus \K$.
 As $F$ is a homeomorphism,
we see that if $x\to 0$ then $d(F(x),\K)\to 0$ and vice versa. Thus,
as $u$ is continuous at zero,  
we see that $\tilde u\in C(B(2)\setminus \K^{int})$ and
$\tilde u$ has the constant value $u(0)$ on $\p \K$.
Moreover, as $F^{-1}\in C^1(B(2)\setminus \K)$, $\|DF^{-1}\|\leq C_0$ in ${B(2)\setminus \K}$ 
 and $u$ is in the Royden algebra  ${\mathcal R}(B(R))$ for all $R<2$,  we have
 by \cite[Thm.\ 3.8.2]{AIM} that the chain rule holds  implying} that
$D\tilde u=((Du)\circ F^{-1})\,\cdotp DF^{-1}$   a.e.\ in $B(2)\setminus \K$.
Let $0<R'<R''<2$.
Then
\ba |D\tilde u(z)|\leq C_0\|Du\|_{C(\overline B(R''))},\quad
\hbox{for }z\in F(B(R''))\setminus \K.
\ea
As $F$ and $F^{-1}$ are $C^1$ smooth up to $\p B(2)$, $\tilde u\in W^{1,1}(B(2)\setminus B(R'))$. These give
 $\tilde u\in W^{1,1}(B(2)\setminus \K)$.
Let  $\tilde v\in W^{1,1}(B(2))$ be a function that coincides with} $\tilde u$ in $B(2)\setminus \K$
and with $u(0)$ in $\K$. 

Again,
using the chain rule and the area formula  as in the proof of Lemma \ref{lem: change of coord.} we see 
that 
$
\tilde A^r [\tilde v]=A^r [u]$ for 
$r>1.$
Applying monotone convergence theorem twice, we obtain
\beq\label{ineq A}
Q_{\tilde \sigma}[h]\leq A_{\tilde \sigma} [\tilde v]=
\lim_{r\to 0}\tilde A^r[ \tilde v]=\lim_{r\to 0}A^r[u]=
Q_\sigma[h].
\eeq


{\lltext  As we have already proven the opposite inequality, this proves 
the claim (i).}

(ii) Assume that $\tilde \sigma$
is a  cloaking conductivity obtained by
the transformation map $F$ and background conductivity $\sigma\in L^\infty(B(2))$,
$\sigma\geq c_1>0$ but  {\lltext that opposite to the claim, we have} $\tr(\tilde \sigma)\in L^1(B(2)\setminus \K)$.
{\ltext
Using formula (\ref{uusi2}) and the facts that $\det(\tilde \sigma)=\det(\sigma\circ F^{-1})$
is bounded from above and below by strictly positive constants
and $\tr(\tilde \sigma)\in L^1(B(2)\setminus \K)$, we see that
$\tr(\tilde \sigma^{-1}) \in L^1(B(2)\setminus \K)$.  Hence
by  Lemma \ref{lem: Inequality}, $K_{\tilde \sigma} \in L^1(B(2)\setminus \K)$.
Let $ G:B(2)\setminus \K\to B(2)\setminus \{0\}$ be the inverse map of $F$.
Using the formulas (\ref{uusi1}), (\ref {transf}), and (\ref{eq: distortion ineq.})
we see that 
\ba
\|\tilde\sigma(y)\|= \frac {\|DF(x)\,\cdotp \sigma(x)\,\cdotp DF(x)^t\|}{J(x,F)}
\geq \frac {\|DF(x)\|^2}{J(x,F)K_\sigma(x)}=
\frac {K_F(x) }{K_\sigma(x)},\quad x=F^{-1}(y).
\ea
As $K_ G=K_F\circ F^{-1}$, cf.\ \cite[formula (2.15)]{AIM} and 
$\|\tilde\sigma(y)\|\leq K_{\tilde \sigma}(y)$, the above yields  
 $K_ G  \in L^1(B(2)\setminus \K)$. Hence, we see using
 (\ref{eq: W12}) that}
$F= G^{-1}$ is  in $W^{1,2}(B(2)\setminus \{0\})$ and
$
\|DF\|_{L^2(B(2)\setminus \{0\})}\leq 2\|K_G\|_{L^1(B(2)\setminus\K)}.
$
By the removability of singularities in Sob\-ol\-ev spaces,
see \cite{KKM}, this implies that 
$F: B(2)\setminus \{0\}
\to  B(2)\setminus \K$ can
be extended to a function $F^{ext}:B(2)\to \C$, 
$F^{ext}\in W^{1,2}(B(2))$. As the distortion $K_F$ of 
the map $F$ is finite a.e., also the map $F^{ext}$ is a finite
distortion map, see  \cite[Def.\ 20.0.3]{AIM}.
Thus, as $F^{ext}\in W^{1,2}_{loc}(B(2))$, it follows  
from the continuity theorem of finite distortion maps \cite[Thm.\ 20.1.1]{AIM}
or \cite{Manfredi} that 
$F^{ext}:B(2)\to \C$ is continuous. Let $y_0=F(0)$.
%
%
Then the set $F^{ext}(\overline B(2))=(\overline B(2)\setminus \K)\cup\{y_0\}$
is not closed as $\p K$ contains more that one point and thus it is not compact.
This is a contradiction with the fact that $F^{ext}$ is continuous. This proves
the claim (ii).
}
\hfill \proofbox  \medskip

Next we prove the claim concerning the last counterexample. 
\smallskip

\noindent {\bf Proof of Theorem \ref{thm: hologram}}. Let us start by reviewing the  properties of the  Iwaniec-Martin maps.
Let  ${\mathcal A}_1:[1,\infty]\to [0,\infty]$ 
be a strictly increasing positive smooth function with $\A_1(1)=0$ which is
satisfies the condition (\ref{IM cond 1}).
Then by   \cite[Thm.\ 11.2.1]{IM}   there 
exists  a $W^{1,1}$-homeomorphism $F:B(2)\setminus\{0\}\to B(2) \setminus \overline B(1)$ which 
Beltrami coefficient $\mu$  satisfies 
\beq
\label{IM3 pre}
\int_{B(2)\setminus\{0\}} \, \exp\big(
{\mathcal A}_1(K_{\mu}(z) )\big) \,dm(z)  < \infty,\quad\hbox{where}\quad  
  K_{\mu}(z) :=\frac {1+|\mu(z)|}{1-|\mu(z)|}.
\eeq
The function $F$ 
 can be obtained using  the construction procedure of \cite[Sec.\ 20.3]{AIM}
(see \cite[Thm.\ 11.2.1]{IM}  the original construction)
as follows:  Let $S_{}(t)$ be solution of the equation 
\beq\label{IM example A}
\A_{1}(S_{}(t))=1+\log (t^{-1}),\quad 0<t\leq 1. 
\eeq
Then $S:(0,1]\to [1,\infty)$ is well defined decreasing function, $S(1)=1$ and
with suitably chosen $c_1>0$, the function
\beq\label{IM example}
 F_{}(z)=\frac z{|z|}\rho_{}(|z|),\quad \rho_{}(s)=1+c_1\bigg(\exp\big(\int_0^{s}\frac {dt}{t\,S_{}(t)}\big)-1\bigg)
\eeq
 is a homeomorpism $F:B(2)\setminus \{0\}\to B(2)\setminus \overline B(1)$. 
 We say that $F$ is the  Iwaniec-Martin map corresponding to the weight function $\A_1(t)$.


Next let  ${\mathcal A}:[1,\infty]\to [0,\infty]$ 
be a strictly increasing positive smooth function with $\A(1)=0$ which is
satisfies the condition (\ref{IM cond 1}) and let
$F_1$ be the  Iwaniec-Martin map corresponding to the weight function $\A_1(t)=\A(4t)$.

Using the inverse of the map $F_1$ we define  $\sigma_1=(F_1^{-1})_* 1$ on $B(2)\setminus \{0\}$
and consider this function as an a.e.\ defined measurable function on $B(2)$.
Using the definition of push forward,
(\ref{IM3 pre}), we see that $\det(\sigma_1)=1$
 and $K_{\sigma_1}(z)=K(F_1^{-1}(z),F_1^{-1})=K_\mu(z)$. Thus 
 Lemma \ref{lem: Inequality} and the fact that $F_1$ satisfies (\ref{IM3 pre}) 
 with the weight function  $\A_1(t)=\A(4t)$
 yield that
 $\sigma_1$ satisfies (\ref{IM3}) with the weight function $\A(t)$.

{Recall that the conductivity 
$\gamma_1$ that is identically 1 in $B(2) \setminus \overline B(1)$
and zero in $\overline B(1)$. 
Next, we consider the minimization problem
(\ref{def: DN map}) with the conductivities
$\gamma_1$ and $\sigma_1$. 
To this end, we  make analogous definitions to the proof of Thm.\ \ref{prop: invisibility}.
For $1<r<2$ let
 $\gamma_r$ be 
a conductivity that is 1
 in $B(2)\setminus  B(r)$ and is 
$0$ in  $ B(r)$. Similarly,
let
$\sigma_r$ be 
a conductivity that coincide
with $\sigma_1$ in $B(2)\setminus B(r-1)$ and is 
$0$ in  $B(r-1)$. 
 
As in (\ref{Q are the same}) and (\ref{ineq: a}), we see  for $h\in H^{1/2}(\p B(2))$ and $r>1$
\beq
Q_{\sigma_r}[h]=Q_{\gamma_r}[h],\quad
\label{Q estimated} 
Q_{\sigma_r}[h]\leq Q_{\sigma_1}[h],\quad Q_{\gamma_r}[h]\leq Q_{\gamma_1}[h].\eeq

Let $h\in H^{1/2}(\p B(2)).$ {\lltext For $1\leq r<2$
the} solution of the boundary value problem
\ba
\Delta w_r=0\quad\hbox{in }B(2)\setminus \overline B(r),\quad w_r|_{\p B(2)}=h,\quad \p_\nu w_r|_{\p B(r)}=0
\ea
satisfies $Q_{\gamma_r}[h]=\|\nabla w_r\|^2_{L^2(B(2)\setminus \overline B(r))}$
and it is easy to see (c.f.\ \cite{Kato}) that 
\beq\label{Q has limit} 
\lim_{r\to 0} Q_{\gamma_r}[h]=Q_{\gamma_1}[h],\quad\hbox{for }h\in H^{1/2}(\p B(2)).
\eeq
Let $w=w_1$. {\ltext Note that $w\in W^{1,2}(B(2)\setminus \overline B(1)).$}


Let us consider the function $v=w\circ {F_1}$. {\ltext As ${F_1}$ is $C^1$-smooth in 
$\overline B(2)\setminus \{0\}$  and the function
$w$ is
$C^1$-smooth in $ \overline B(R)\setminus \overline B(1)$ for all $1<R<2$ we have 
by the chain rule 
$Dv(x)=(Dw)({F_1}(x))\,\cdotp D{F_1}(x)$ for all $x\in B(2)\setminus \{0\}$}. {\ltext As 
$Dw\in L^2(B(2)\setminus B(R))$ and 
$Dw\in L^\infty(B(R)\setminus \overline B(1))$ for all $1<R<2$, and
\ba
D{F_1}(x)=\frac {\rho(|x|)}{|x|}(I-P(x))+ \rho'(|x|)P(x)   
\ea
where $P(x):y\mapsto |x|^{-2}(x\, \,\cdotp y)x$ is the projector to the radial direction at the point $x$, we 
using (\ref{IM example}) that
see that
$\|D{F_1}(x)\|\leq  C|x|^{-1}$ with some $C>0$ and 
\beq\label{eq: example}
Dv\in L^p(B(2)\setminus \{0\}),\quad\hbox{for any $p\in (1,2)$.}
\eeq
Thus  $v\in W^{1,p}(B(2)\setminus \{0\})$ with any $p\in (1,2)$ and
by the removability of singularities 
in Sobolev spaces, see e.g.\ \cite{KKM}, function $v$ can be
considered as a measurable function in $B(2)$ for which $v \in W^{1,p}(B(2))$.  Thus
$v$ is in the domain of definition of the quadratic form $A_{\sigma_1}$. 

As $w\in C^1(\overline B(R)\setminus \overline B(1))$ for all $1<R<2$ 
and ${F_1}$ is $C^1$-smooth in $\overline B(2)\setminus \overline B(1))$, we can use again the chain rule, 
the area formula, and the monotone convergence theorem to obtain
\beq\label{eq: leq}
& &Q_{\sigma_1}[h]\leq A_{\sigma_1}[v]=
\lim_{R \to 2}\lim_{\rho \to 0}\int_{B(R)\setminus \overline B(\rho )}\nabla v\, \,\cdotp \sigma_1\nabla v\,dm(x)\\ \nonumber
& &=\lim_{R \to 2}\lim_{\rho \to 0} \int_{{F_1}(B(R)\setminus \overline B(\rho ))}\nabla w\, \,\cdotp \gamma_1\nabla w\,dm(x) 
=Q_{\gamma_1}[h].
\eeq

Next, consider an inequality opposite to (\ref{eq: leq}). We have
by (\ref{Q estimated}) and (\ref{Q has limit}) 
\beq\label{eq: leq222}
Q_{\sigma_1}[h]\geq \lim_{r\to 1}Q_{\sigma_r}[h] 
=\lim_{r\to 1}Q_{\gamma_r}[h]=Q_{\gamma_1}[h].
\eeq
The above inequalities  prove the claim.}
\hfill \proofbox 

\section{Complex geometric optics solutions}\label{sec: proofs}

{\lltext  In this section we assume that $\A$ satisfies the  almost linear growth condition (\ref{cond for A function}).}
 
\subsection{Existence and properties  of the complex geometric optics solutions}
\label{subsec: CGO solutions}

{\ltext
Let us start with the observation that if $\sigma_0\in \Sigma(\Omega_0)$ is  a conductivity  in  a smooth 
simply connected domain $\Omega_0\subset\C$, 
and $\sigma_1$ is a conductivity in a larger smooth domain $\Omega_1$ which coincides
with $\sigma_0$ in $\Omega_0$ and is one  in $\Omega_1\setminus \Omega_0$,
then $Q_{\sigma_0}$ determines $Q_{\sigma_1}$ by formula
\ba
Q_{\sigma_1}[h]=\inf\{ \int_{\Omega_1\setminus \Omega_0} |\nabla v|^2\,dm(z)+
Q_{\sigma_0}[v|_{\p \Omega_0}]\ ;\ v\in W^{1,2}(\Omega_1\setminus \overline \Omega_0),\ 
v|_{\p \Omega_1}=h\}.
\ea
%
 This observation implies that we may consider inverse problems by assuming that the conductivity
 $\sigma$ is the identity near $\p \Omega$
without loss of generality. Also, we may assume that $\Omega =\D$, which
we do below.}

 The main result of this section is the following
 uniqueness and existence theorem for the 
complex geometrical optics
solutions.

\begin{theorem}\label{thm:CGO} Let $\sigma\in \Sigma_\A(\C)$ be a conductivity 
such that 
 $\sigma(x)=1$ for $x\in \C\setminus \Omega$. Then for every $k\in \C$ there is a unique solution $u(\,\cdotp,k)\in W^{1,P}_{loc}(\C)$
 for 
\ba
& & \nabla_z\cdot\sigma(z)\nabla_z u(z,k)=0, \quad \hbox{for a.e. } z \in \C, \\ 
& &u(z,k)  = e^{ikz}(1+\mathcal O(\frac{1}{z}))\quad \hbox{ as }|z|\to \infty.
\ea 
\end{theorem}

We point out that the regularity $u\in W^{1,P}_{loc}(\C)$ is optimal
in the sense that 
 a slightly stronger assumption $u\in W^{1,2}_{loc}(\C)$
would not be valid for the solutions and a slightly weaker 
assumption
$ 
u \in \bigcap_{1 < q < 2} W^{1,q}_{loc}(\C)
$
would not be strong enough for obtaining the uniqueness, see \cite{AIM} and 
the equivalence of the conductivity equation and the Beltrami equation discussed in
(\ref{eq:  BB})-(\ref{eq: anisotropic Beltrami AA}).

We prove Thm.\ \ref{thm:CGO}  in several steps. Recalling the reduction to the Beltrami equation (\ref{eq: anisotropic Beltrami AA}), we 
 start with the following lemma, {\mmmtext where we denote} 
 \ba
B_\A(\D)=\{\mu \in  L^\infty(\C)&;& \supp(\mu)\subset \overline \D,\
0\leq \mu(x)< 1\quad \hbox{a.e.},\ \hbox{and }\\
& &\int_\D \exp(\A(K_\mu(z)))\,dm(z)<\infty\}.\ea

\begin{lemma}\label{lemma 1} Assume that 
$\mu \in B_\A(\D)$ and $f\in W^{1,P}_{loc}(\C)$
satisfies
\beq\label{eq: repeated}
& &\dbar f(z)=\mu(z) {\p f(z)},\quad\hbox{for a.e. }z\in \C,
\\
& &\label{23b}
f(z)=\beta e^{ikz}(1+\O(\frac 1 z ))\quad\hbox{for }|z|\to \infty,
\eeq
where $\beta\in \C\setminus \{0\}$ and 
$k\in \C$.
Then
\beq\label{eq: representation}
f(z)=\beta e^{ik\principalPhi(z)},
\eeq
where $\principalPhi\in W^{1,P}_{loc}(\C)$ is a homeomorphism $\principalPhi:\C\to \C$,
$\dbar \principalPhi(z)=0$ for $|z|>1$,
$K(z,\principalPhi)=K(z,f)$ for a.e.\ $z\in \C$, and
\beq\label{eq: asymptotics}
\principalPhi(z)=z+\O(\frac 1z)\quad\hbox{for }|z|\to \infty.
\eeq
\end{lemma} 

\medskip

\noindent {\bf Proof.} 
{\mmmtext By Theorem \ref{Stoilow} we have} for $f$ the
Stoilow factorization $f=h\circ \principalPhi$ where $h:\C\to \C$ is holomorphic function
and $ \principalPhi$ is the principal solution of (\ref{eq: repeated}).
This and 
the formulae (\ref{23b}) and (\ref{eq: asymptotics}) imply
\ba
\frac {h(\principalPhi(z))}{\beta e^{ik\principalPhi(z)}}=
\frac {f(z)}{\beta e^{ik\principalPhi(z)}}\to 1\quad \hbox{when }|z|\to \infty.
\ea
Thus, $h(\zeta)=\beta e^{ik\zeta}$ for all $\zeta\in \C$, and
$f$ has the representation (\ref{eq: representation}). {\lltext
The claimed properties of $\principalPhi$ follows from the formula  (\ref{eq: representation})
and the similar properties of $f$.}
 \proofbox 
\medskip

Next we consider case where $\beta=1$.
Below we will use the fact that 
if $\principalPhi:\C\to \C$ is a homeomorphism such that $\principalPhi\in W^{1,1}_{loc}(\C)$,
$\principalPhi(z)-z=o(1)$ as $z\to \infty$ and
{\lltext that $\principalPhi$ is analytic outside the disc $\overline B(r)$, $r>0$ 
 then by \cite[Thm.\ 2.10.1 and (2.61)]{AIM} 
\beq\label{eq: Kari 1}
\hbox{$|\principalPhi(z)|\leq |z|+3r$ for $z\in \C$ and $|\principalPhi(z)-z|\leq r$ for $|z|>2r$}.
\eeq
In particular,} the map $\principalPhi$ defined in Lemma \ref{lemma 1}  satisfies this with $r=1$.

\begin{lemma}\label{lemma 2 A}  Assume that $\nu,\mu:\C\to \C$
are measurable functions satisfying 
\beq\label{Ap1 B}
& &\mu(z)=\nu(z)=0,\quad \hbox{ for }z\in \C\setminus \D,\\
& &|\mu(z)|+|\nu(z)|<1,\quad \hbox{ for a.e. }z\in \D, \label{Ap1 Bb}
\eeq
and that
$
K_{\mu,\nu}(z)
$
defined in (\ref{distor2})
satisfies
\beq
& &\int_\D \exp(\A(K_{\mu,\nu}(z)))\,dm(z)<\infty.\label{Ap2 B}
\eeq
Then for $k\in \C$ the equation
\beq\label{eq: eq 1}
\partial_{\overline z} f = \mu \partial_z f + 
\nu  \,\overline{\partial_z f},\quad z\in\C
\eeq
has at most one solution 
$f\in W^{1,P}_{loc}(\C)$ satisfying  
\beq\label{eq: eq 2}
f(z)= e^{ikz}(1+\O(\frac 1 z ))\quad\hbox{for }|z|\to \infty.
\eeq
\end{lemma} 

\medskip

\noindent {\bf Proof.} Observe that we can write equation (\ref{eq: eq 1})
in the form
\beq\label{eq: modified anisotropic Beltrami}
\partial_{\overline z} f = \tilde \mu \,{\partial_z f},\quad z\in\C,
\eeq
where the coefficient  $\tilde \mu$ is given by (\ref{eq: modified anisotropic Beltrami A1}).
Since $|\tilde \mu(z)|\leq |\mu(z)|+|\nu(z)|$, we see that
$\tilde \mu\in B_\A(\D)$.

Next, assume equation (\ref{eq: modified anisotropic Beltrami})
has two solutions $f_1$ and $f_2$  having the asymptotics (\ref{eq: eq 2}). 
Let $\e>0$ and consider function $
f_\e(z)=f_1(z)-(1+\e)f_2(z).
$
Then, $f_\e\in W^{1,P}_{loc}(\C)$, function $f_\e$ satisfies (\ref{eq: eq 1}), and 
\ba
f_\e(z)=- \e e^{ikz}(1+\O(\frac 1 z ))\quad\hbox{for }|z|\to \infty.
\ea
By Lemma \ref{lemma 1} and (\ref{eq: Kari 1}), there is $\principalPhi_\e(z)$ such that
\ba
f_\e(z)=f_1(z)-(1+\e)f_2(z)=-\e e^{ik\principalPhi_\e(z)}
\ea 
and $|\principalPhi_\e(z)|\leq\, |z|+3$. Then for any $z\in \C$
we have that
\ba
f_1(z)-f_2(z)=\lim_{\e\to 0}f_\e(z)=0.
\ea
Thus $f_1=f_2$. \proofbox 
\medskip

\subsection{Proof of Theorem \ref{thm:CGO}}

In following, we use general facts for weakly monotone mappings,
and for this end, we recall some basic facts. Let $\Omega\subset \C$
be open and $u\in W^{1,1}(\Omega)$ be real valued.   
We say that $u$ is weakly monotone,  if the both functions $u(x)$ and $-u(x)$ satisfy 
the maximum principle in the following weak sense: For any $a\in \R$
and relatively compact open sets $\Omega'\subset \Omega$, 
\ba
\max(u(z)-a,0)\in W^{1,1}_0(\Omega')\hbox{ implies that $u(z)\leq a$ for a.e.\ $z\in \Omega'$},
\ea
see \cite[Sec.\ 7.3]{IM}. 
{\lltext We remark that if $f\in W^{1,1}_{loc}(\Omega_1)$ 
and $f:\Omega_1\to \Omega_2$ is homeomorphism, where $\Omega_1,\Omega_2\subset \C$
are open, the real part of $f$ is 
are weakly monotone.  By \cite[Lem.\ 20.5.8]{AIM}, 
 if $f\in W^{1,1}(\Omega)$ is the solution
of the Beltrami equation $\dbar f=\mu \p f$ with a 
Beltrami coefficient $\mu$ satisfying $|\mu(z)|<1$ for a.e.\ $z\in \C$, then
the real and the imaginary parts of $f$ are weakly monotone functions. 
An important property of weakly monotone functions is
that their modulus of continuity can be estimated in an explicit way.
{\mmmtext Let $M_P(t)$ be the $P$-modulus, that is, the function
 determined by the condition 
\ba
\hbox{For $M=M(t)$ we have }\int_1^{1/t} P( s M)\frac {ds}{s^3}=P(1),
\ea
cf.\ (\ref{Ap171}).
The function  $M_P:[0,\infty)\to [0,\infty)$ is continuous at zero and $M_P(0)=0$.
Then by \cite[Thm.\ 7.5.1]{IM} 
it holds that
 if $z',z\in \Omega$  satisfy
$B(z,r)\subset \Omega$, $r<1$ and $|z'-z|<r/2$, 
and $f\in W^{1,P}(\Omega)$ is a weakly monotone function, then
for almost every $z,z'\in B(z,r)$ we have
\beq\label{modulus of cont.}
|f(z')-f(z)|\leq 32\pi r\, \|\nabla f\|_{(P,r)} M_P(\frac {|z-z'|}{2r}),
\eeq
where
\ba
\|\nabla f\|_{(P,r)}=\inf\left\{\frac 1\lambda\,;\ \lambda>0,\ \frac 1{\pi r^2}\int_{B(z,r)}P(\lambda |\nabla f(x)|)
\,dm(z)\leq P(1)\right\}.
\ea
As we will see, this can be used to estimate  the modulus of continuity
of principal solutions of Beltrami equations corresponding to $\mu\in B_\A(\D)$.}

{\ltext Below, we use the unimodular function $e_{k}$  given by
$e_{k}(z)=e^{i(kz+\overline{k}\overline{z})}$.}
The following lemma  shows the existence of
the complex geometric solutions for degenerated conductivities.

\begin{lemma}\label{lemma 2}  Assume that $\mu$ and  $\nu$
 satisfy (\ref{Ap1 B}), (\ref{Ap1 Bb}), (\ref{Ap2 B})
{\ltext and let $k\in \C\setminus\{0\}$.} Then  the equation (\ref{eq: eq 1}) has a solution
$f\in W^{1,P}_{loc}(\C)$ 
satisfying
asymptotics (\ref{eq: eq 2}). Moreover, this solution can be
written in the form
\beq\label{representation 3}
f(z) = e^{ik\uusivarphi(z)}
\eeq
where $\uusivarphi:\C \to \C$ is a homeomorphism satisfying
the asymptotics $ \uusivarphi(z) = z +\mathcal O({z}^{-1}) $. Moreover,
 for $R>1$ 
\beq\label{Kari 1}
\hspace{-7mm}\int_{B(R)} P(|D\uusivarphi(x)|)\,dm(x)\leq C_\A(R) \int_{B(R)} \exp(\A(K_{\mu,\nu}(z)))\,dm(z),
\eeq
where $C_\A(R)$ depends on $R$ and the weight function $\A$.
{\mmmtext In addition,
\beq\label{eq: for log for phi}
\dbar \uusivarphi(z)=
 \mu(z)  {\p  \uusivarphi(z)}
-\frac {\lambdapoistettu \overline k}{k}
 \nu(z) e_{-k}( \uusivarphi(z)) \overline {\p  \uusivarphi(z)},
\quad\hbox{for a.e. }z\in \C.
\eeq
} 
%
\end{lemma}

\noindent {\bf Proof.} 
Let us approximate the functions $\mu$ and $\nu$ with
functions
\beq\label{truncted mu}
& &\mu_n(z)=\left\{\begin{array}{cl} 
\mu(z)& \hbox{if }| \mu(z)|+|  \nu(z)|\leq 1-\frac 1n,\\
\frac {\mu(z)} {|\mu(z)|}(1-\frac 1n)& \hbox{if }|\mu(z)|+|  \nu(z)|>1-\frac 1n,\\
\end{array}\right.\\
\label{truncted nu}
& &  \nu_n(z)=\left\{\begin{array}{cl} 
  \nu(z)& \hbox{if }|  \mu(z)|+|  \nu(z)|\leq 1-\frac 1n,\\
\frac {  \nu(z)} {|  \nu(z)|}(1-\frac 1n)& \hbox{if }|  \mu(z)|+|  \nu(z)|>1-\frac 1n,\\
\end{array}\right.
\eeq
where $n\in\Z_+$.
Consider the equations
\beq\label{eq: eq 1 n}
\dbar f_n(z)= \mu_n(z){\p f_n(z)}+ \nu_n(z)\overline {\p f_n(z)},\quad\hbox{for a.e. }z\in \C,\\
\label{e: asymptoti}
f_n(z)=\lambdapoistettu e^{ikz}(1+\O(\frac 1 z ))\quad\hbox{for }|z|\to \infty.
\eeq
By Lemma \ref{lemma 2 A}  equations  (\ref{eq: eq 1 n})-(\ref{e: asymptoti}) have
at most one solution $f_n\in W^{1,P}_{loc}(\C)$.
The existence of the solutions can be seen as in the proof 
of  \cite[Lem.\ 3.5]{ALP}: By \cite[Lem.\ 3.2]{ALP},  solutions $f_n$ for (\ref{eq: eq 1 n})-(\ref{e: asymptoti})
can be constructed  via the formula $f_n=h\circ g$, where 
and $g$ is the principal solution of 
$\dbar g=\hat \mu{\p g}$, 
 constructed in Thm.\ \ref{Stoilow} and
$h$ is
the solution of $\dbar h=(\hat \nu\circ g^{-1})\overline {\p h}$, $h(z)=e^{ikz}(1+\O(\frac 1 z ))$
constructed in \cite[Thm.\ 4.2]{AP}  where
$\hat \nu  =(1+\nu_n)\overline\mu_n $ and  $\hat \mu=\mu_n(1+\hat \nu)$ and moreover, it holds that $f_n\in W^{1,2}_{loc}(\C)$.
%

Let us now define the coefficient $\tilde \mu$ according to  formula (\ref{eq: modified anisotropic Beltrami A1}),
and define an approximative coefficient
$\tilde \mu_n$ using formula (\ref{eq: modified anisotropic Beltrami A1})
where $\mu$ and $\nu$ are replaced by  $\mu_n$ and $\nu_n$ and
$f$ by $f_n$. We can write the equation
 (\ref{eq: eq 1 n})  in the form
\beq\label{eq: eq 1 modified}
\dbar f_n(z)= \tilde \mu_n(z){\p f_n(z)},\quad\hbox{for a.e. }z\in \C
\eeq
where $|\tilde \mu_n|\leq 1-n^{-1}.$}

By  (\ref{e: asymptoti}), (\ref{eq: eq 1 modified}), and Lemma \ref{lemma 1}, function $f_n$ can be written in the form
\beq\label{eq: representation n}
f_n(z)=\lambdapoistettu e^{ik \uusivarphi_n(z)},
\eeq
where $ \uusivarphi_n$ is a homeomorphism, 
$\dbar  \uusivarphi_n(z)=0$ for $|z|>1$,
$K(z, \uusivarphi_n)=K(z,f_n)$ for a.e.\ $z\in \C$, 
and
\beq\label{eq: asymptotics n}
 \uusivarphi_n(z)=z+\O(\frac 1z)\quad\hbox{for }|z|\to \infty.
\eeq
Then
\ba
|\dbar f_n(z)|=|\tilde \mu_n(z)|\, |\p f_n(z)|\leq |\tilde \mu(z)|\, |\p f_n(z)|.
\ea
{\ltext

Let us consider next 
$a,b>0$ and  $0\leq t\leq (ab)^{1/2}$.  Using the definition (\ref{hotelli57}) of $P(t)$ 
we see
\ba
& P(t)\leq \exp(\A(a)),\hspace{20mm} &\hbox{for }t^2\leq e^{\A(a)},\\
& P(t)\leq  \frac {ab}
 {\A^{-1}(\log \exp(\A (a)))}=b,\quad  &\hbox{for }t^2> e^{\A(a)},
\ea   
which imply the inequality $ P (t)\leq b+\exp(\A(a))$. Due to 
the distortion equality 
(\ref{eq: distortion ineq.}) we can use this for 
$a=K(z, \uusivarphi_n)$, $b=J(z, \uusivarphi_n),$ and $t=|D \uusivarphi_n(z)|$
and obtain 
\beq\label{eq: Lassi star1} 
 P (|D \uusivarphi_n(z)|)
\leq J(z, \uusivarphi_n)+\exp(\A(K(z, \uusivarphi_n))). 
\eeq

%
Then, we see using (\ref{eq: Kari 1}) and the fact that $ \uusivarphi_n$ is a homeomorphism that
\beq
& &\int_{B(R)} P(|D \uusivarphi_n(z)|)\,dm(z)\leq \nonumber
 \int_{B(R)} J(z, \uusivarphi_n)\,dm(z)+
 \int_{B(R)} e^{\A(K(z, \uusivarphi_n)}dm(z) \\  \label{eq: essential2}
 &&\leq
 m( \uusivarphi_n(B(R)))+
 \int_{B(R)} \exp(\A(K_{\tilde \mu}(z))) dm(z)\\  \nonumber
 &&\leq
 \pi  (R+3)^2+
 \int_{B(R)} \exp(\A(K_{\tilde \mu}(z))) dm(z)
\eeq
is finite by the assumption (\ref{Ap2 B}).} We emphasize that
the fact that $ \uusivarphi_n$ is a homeomorphism is the essential
fact which together with the inequality (\ref{eq: Lassi star1})
yields the Orlicz estimate (\ref{eq: essential2}).

The estimate (\ref{eq: essential2})  together with the inequality
(\ref{modulus of cont.}) implies that 
the functions $ \uusivarphi_n$ have
uniformly bounded modulus of continuity {\mmmtext in all compact sets of $\C$}. Moreover, by  (\ref{eq: Kari 1}),
$| \uusivarphi_n(z)|\leq |z|+3.$

{\mmmtext Next we consider the Beltrami equation for $ \uusivarphi$.
To this end, let $\psi\in C^\infty_0(\C)$ and $R>1$ 
be so large that $\supp(\psi)\subset B(R)$.
Since the family $\{ \uusivarphi_n\}_{n=1}^\infty$
is uniformly bounded in the space $W^{1,P}(B(R))$ and
$W^{1,P}(B(R))\subset W^{1,q}(B(R))$ for some $q>1$,  
we see} that there is a subsequence
 $ \uusivarphi_{n_j}$ that converges 
weakly in  $W^{1,q}(B(R))$ {\ltext to some limit $  \uusivarphi$ when $j\to \infty$.}
Let us denote 
\ba
\kappa_n(z)=-\frac {\lambdapoistettu \overline k}{\,k}
 \nu_n(z) e_{-k}( \uusivarphi_n(z)),\quad
\kappa(z)=-\frac {\lambdapoistettu \overline k}{\,k}
 \nu(z) e_{-k}( \uusivarphi(z)).
\ea
Moreover, functions $ \uusivarphi_n$ are uniformly bounded 
and have a uniformly bounded modulus of continuity {\ltext in compact sets by (\ref{modulus of cont.})} and thus
by Arzela-Ascoli theorem there is a subsequence $ \uusivarphi_{n_j}$, denoted also
by  $ \uusivarphi_{n_j}$,
that converges uniformly to some function $\uusivarphi'$ in $B(R)$ for all $R>1$.
As $ \uusivarphi_{n_j}$ converge in $C(\overline B(R))$
 uniformly to $\uusivarphi'$ and weakly in $W^{1,q}(B(R))$ to $  \uusivarphi$
we see using convergence in distributions that $\uusivarphi'=  \uusivarphi$. 
 Thus, we see that 
\ba
\lim_{j\to \infty }e_{-k}( \uusivarphi_{n_j}(z))=e_{-k}(  \uusivarphi(z))
\quad\hbox{uniformly for }z\in B(R),
\ea
and by dominated convergence theorem 
$\kappa_n\to \kappa$ in $L^p(B(R))$ where $1/p+1/q=1$.

{\ltext As $\uusivarphi_n:\C\to\C$ is a homeomorphism and $\uusivarphi_n\in W^{1,1}_{loc}(\C)$,
we can use chain rules (\ref{chain rules}) a.e.\
 by the Gehring-Lehto theorem, see   \cite[Cor.\ 3.3.3]{AIM},
and see using 
 the equations
(\ref{eq: eq 1 n}) and (\ref{eq: representation n}) that 
\beq\label{eq: for log for n}
\dbar \uusivarphi_n(z)=
 \mu_n(z)  {\p  \uusivarphi_n(z)}
-\frac {\lambdapoistettu \overline k}{k}
 \nu_n(z) e_{-k}( \uusivarphi_n(z)) \overline {\p  \uusivarphi_n(z)},\quad\hbox{for a.e. }z\in \C.
\eeq}

{\lltext  Recall that  there is convex function $\Phi:[0,\infty)\to[0,\infty)$ such that $\Phi(t)\leq P(t)+c_0\leq 2\Phi(t)$.
By \cite[Thm.\ 13.1.2]{ABM} the map 
$\phi\mapsto \int_{B(R)} \Phi(|D \phi(x)|)\,dm(x)$
is weakly lower
semicontinuous in $W^{1,1}(B(R))$. 
By (\ref{eq: essential2}) the integral of $\Phi(|D \uusivarphi_n|)$ in uniformly bounded in
$n\in \Z_+$ over any disc $B(R)$. In particular, this yields that
  $ \uusivarphi\in W^{1,P}(B(R))$ for $R>1$ and that
  (\ref{Kari 1}) holds.}
  
Furthermore, as $| \uusivarphi(z)|\leq\, |z|+3$, this yields 
 that
\beq\label{eq: representation 2}
  f(z):=\lambdapoistettu e^{ik \uusivarphi(z)}\in W^{1,P}_{loc}(\C).
\eeq
Define next $ \uusivarphi_n(\infty)=  \uusivarphi(\infty)=\infty$. As $ \uusivarphi_n$ and $ \uusivarphi$ are 
conformal at infinity, we see  using Cauchy formula
for $( \uusivarphi_n(1/z)- \uusivarphi(0))^{-1}$ that
that
\beq\label{eq: asymptotics n2}
 \uusivarphi(z)=z+\O(\frac 1z)\quad\hbox{for }|z|\to \infty.
\eeq

As    $D  \uusivarphi_{n_j}$ converges 
weakly in  $L^{q}(B(R))$ to $D\uusivarphi$ and  their norms are uniformly bounded, 
we have
\ba
& &\hspace{-8mm}|\int_\C  (\dbar   \uusivarphi-\mu { \p  \uusivarphi}- \kappa \overline { \p  \uusivarphi})\psi\,dm(z)|=
\lim_{j\to \infty}
|\int_\C  (\dbar  \uusivarphi_{n_j}-\mu { \p \uusivarphi_{n_j}} -\kappa \overline {\p  \uusivarphi_{n_j}})\psi\,dm(z)|
\\
&\leq &
\lim_{j\to \infty}
|\int_\C( (\mu_{n_j}-\mu) {\p  \uusivarphi_{n_j}}+  (\kappa_{n_j}-\kappa) \overline {\p  \uusivarphi_{n_j}})\psi\,dm(z)|
\\
&\leq &
\lim_{j\to \infty}(\|\mu_{n_j}-\mu\|_{L^p(B(1))}+
\|\kappa_{n_j}-\kappa\|_{L^p(B(1))} )\| \p  \uusivarphi_{n_j}\|_{L^q(B(1))}
\|\psi\|_{L^{\infty}(B(1))}=0.
\ea
This implies that $  \uusivarphi(z)$ satisfies (\ref{eq: for log for phi}).

Next we show that $   \uusivarphi$ is homeomorphism.
{\lltext As $K(z)=K_{\nu,\mu}\in L^1_{loc}(\C)$, we have $K(z;\uusivarphi_n)\in L^1_{loc}(\C)$ thus by  (\ref {eq: W12})  the inverse maps
$ \uusivarphi_n^{-1}$  satisfy
$ \uusivarphi_n^{-1}\in W^{1,2}_{loc}(\C;\C)$ and for all $R>1$ the norms}
$\|\uusivarphi_n^{-1}\|_{W^{1,2}(B(R))}$, $n\in \Z_+$ are uniformly bounded. 
Thus by 
the formula (\ref{modulus of cont.}) the family $ (\uusivarphi_n^{-1})_{n=1}^\infty$
 has a uniform modulus on continuity in compact sets. 
Hence, we see that there is a continuous function $\psi:\C\to \C$ such that
$ \uusivarphi_n^{-1}\to \psi$ uniformly on compact sets when $n\to \infty$.
As $ \uusivarphi_n$ are conformal at infinity, we see using again the Cauchy formula
that $ \uusivarphi_{n_j}^{-1}\to \psi$ uniformly on the Riemann sphere $\S^2$ as $j\to \infty$.  
Then, 
\ba
\psi\circ  \uusivarphi(z)=\lim_{j\to \infty}  \uusivarphi_{n_j}^{-1}( \uusivarphi(z))=\lim_{j\to \infty}  \uusivarphi_{n_j}^{-1}(   \uusivarphi_{n_j}(z))=z
\ea
which implies that $   \uusivarphi:\S^2\to \S^2$  is a continuous injective map and hence  a homeomorphism. 

As $  \uusivarphi:\C\to\C$ is a homeomorphism and $  \uusivarphi\in W^{1,1}_{loc}(\C)$,
we can by the Gehring-Lehto theorem use chain rules (\ref{chain rules}) a.e.\
and see using the
  equation
(\ref{eq: for log for phi}) that $f(z)=\lambdapoistettu e^{ik \uusivarphi(z)}$ satisfies  equation (\ref{eq: eq 1}).
By (\ref{eq: asymptotics n2}), $f(z)$ satisfies asymptotics (\ref{eq: eq 2}).
This proves the claim.
%
\hfill\proofbox \bigskip

The above uniqueness and existence results have now
proven Theorem \ref{thm:CGO}. \proofbox

\section{Inverse conductivity problem with degenerate isotropic
conductivity}

In this section we consider  exponentially integrable scalar conductivities $\sigma$.
In particular, we assume that
$\sigma$ is 1 in an open set containing $\C\setminus \D$ and its
 ellipticity function $K(z)=K_\sigma(z)$ 
of the conductivity $\sigma$ satisfies
an Orlicz space estimates 
\beq\label{Las 9.1}
\int_{\B(R_1)} \exp(\exp(qK(x)))\, dm(x)\leq C_0 \quad\hbox{for some }C_0,q>0
\eeq
with $R_1=1$.
Note that by John-Nirenberg lemma, (\ref{Las 9.1})
is satisfied if
\beq\label{Las 9.2}
\exp(qK(x))\in BMO(\D),\quad\hbox{for some }q>0.
\eeq
{\ltext As noted before, we may assume without  loss of generality that $\Omega$
is the unit disc $\D$.}

\subsection{Estimates for principal solutions in Orlicz spaces}

Let us consider next the   principal solution of the Beltrami equation
\beq\label{Las 9.3}
& &\dbar \uusiPhi (z)=\mu(z)\, {\p \uusiPhi (z)},\quad z\in \C \\
& &\uusiPhi (z)=z+O(\frac 1z)\quad \hbox{when }|z|\to \infty.\label{Las 9.4}
\eeq
{\mmmtext For this end, let $R_0\geq 1$ and 
\ba
& &B^p_{exp,N}({B(R_0)})=\{\mu:\C\to \C ;\ |\mu(z)|<1\ \hbox{for a.e.\ }z,\ \ \supp(\mu)\subset {B(R_0)},\\
& &\hspace{4.6cm}\hbox{and }\int_{B(R_0)} \exp(pK_\mu(z))\,dm(z)\leq N\}
\ea
and $B^p_{exp}({B(R_0)})=\bigcup_{N>0}B^p_{exp,N}({B(R_0)})$.
The reason that we use the radius $R_0$ is to be able to apply the obtained results for
the inverse function  of the solution of the Beltrami equation satisfying another Beltrami
equation with modified coefficients, see (\ref{Las 9.15}).

Assume that $p>2$ and $\mu\in B^p_{exp}({B(R_0)})$. Then by \cite[Thm.\ 1.1]{AGRS} we have
the $L^2$-estimate
\beq\label{Kari B1}
\|(\mu S)^m\mu\|_{L^2(\C)}\leq C(p,\beta)m^{-\beta/2}
\int_{B(R_0)} \exp(pK_\mu(z))\,dm(z),\quad 2<\beta<p.
\eeq
In particular, as $ \uusiPhi$ satisfies
\ba
\dbar  \uusiPhi=\dbar ( \uusiPhi-z)=\mu\p ( \uusiPhi -z)+\mu=\mu S\dbar ( \uusiPhi -z)+\mu=\mu S\dbar  \uusiPhi +\mu,
\ea
where
\ba
S\phi (z) =
 \frac 1 {\pi }\int _\C \frac{\p_w \phi (w)}{w-z}\, dm (w),
\ea
is the Beurling operator, (\ref{Kari B1}) yields 
\beq\label{eq: Beurling series}
\dbar  \uusiPhi=\sum_{m=0}^\infty (\mu S)^m\mu,
\eeq
where the series converges in $L^2(\C)$. To analyze
the convergence more precisely, we need a refinement of the $L^p$-scale.
In particular, we will use the Orlicz spaces $X^{j,q}(S)$, $j\in \Z_+$, $q\in \R$, $S\subset \C$ that are defined
{\ltext by
\beq\label{Las 1}
u\in X^{j,q}(S)\quad\hbox{if and only if}\quad \int_S M_{j,q}(u(x))\,dm(x)<\infty
\eeq
where 
\beq\label{Las 2}
 M_{j,q}(t)=|t|^j\log^q(e+|t|).
 \eeq
 We use shorthand notations $X^{q}(S)=X^{2,q}(S)$ and $ M_{q}(t)= M_{2,q}(t)$.}
Although (\ref{Las 1})-(\ref{Las 2}) do not define a norm in $X^{j,q}(S)$,
there is an equivalent norm 
  \beq\label{copied app 2}
 \| u\|_{X^{j,q}(S)}=\sup\{ \int_S |u(x)v(x)|\, dm(x)\ ;\ \int_D G_{j,q}(|u(x)|)\, dm(x)\leq 1\}.
 \eeq
where $G_{j,q}(t)$ is such function that $(M_{j,q},G_{j,q})$ are a Young complementary pair
 (c.f.\ the Appendix) and,
in particular, 
the  following lemma holds.

\begin{lemma}\label{lem: LassiA} We have  for
 $j=1,2,\dots,$ $q\geq 0$

(i) $\int_{B(R_0)} M_{j,q}(u(x) )dm(x)\leq 2\|u\|_{X_{j,q}({B(R_0)})}^j\log^q(e+\|u\|_{X^{j,q}({B(R_0)})})$,

(ii) $\|u\|_{X^{j,q}({B(R_0)})}\leq \phi(\int_{B(R_0)} M_{j,q}(u(x) )dm(x)),$ $\phi(t)=
  t^{1/j}(1+2\log^q(e+t^{-1/j} )).$
\end{lemma}

{\bf Proof.}  (i). Let us denote $M(t)=M_{j,q}(t)$. For this function we use the
equivalent norms
$\|u\|_{M}$ and $\|u\|_{(M)}$ defined in the Appendix. To show the claim we use the inequality 
\beq\label{Eq Lassi 3}
\log(e+st)\leq  2\log(e+s)\,\log(e+t),\quad t,s\geq 0.
\eeq
Let us consider function $w\in X^{j,q}({B(R_0)})$.
By (\ref{Eq Lassi 3}) we have for $k>0$
\beq\label{Lassi-Matti 1}
\int_{B(R_0)} M_{j,q}(kw)\,dm&=&k^j\int_{B(R_0)} |w|^j\log^q(e+k|w| )\,dm\\
& \leq& 2k^j\log^q(e+k )\int_{B(R_0)} M_{j,q}(w)\,dm.\nonumber
\eeq
A function $u\in X^{j,q}({B(R_0)})$ can be written as $u=kw$ where  $k=\|u\|_{(M)}$ and
$\|w\|_{(M)}=1$. 
Then by (\ref{app 5})-(\ref{app 6})  we have $\int_{B(R_0)} M_{j,q}(w )dm=1$ 
hence (\ref{Lassi-Matti 1}) and (\ref{app 4}) yield the claim (i).

(ii) Using  (\ref{Lassi-Matti 1})  and the definition  (\ref{app 2}) of the Orlicz norm,
we see that for all $k>0$
\ba
\|u\|_{X^{j,q}({B(R_0)})}\leq \frac 1k(1+\int_{B(R_0)} \hspace{-3mm}M_{j,q}(ku)\,dm)
\leq \frac 1k(1+2k^j\log^q(e+k ) \int_{B(R_0)} \hspace{-3mm}M_{j,q}(u)\,dm).
\ea
Let $T=\int_{B(R_0)} M_{j,q}(u)\,dm$. Substituting above $k=T^{-1/j}$
we obtain  (ii).
%
%
%
%
%
\hfill \proofbox

\begin{theorem}\label{thm: New 9.1} Assume that  $\mu\in B^p_{exp}({B(R_0)})$,
$2<p<\infty$. Then {\ltext the equations (\ref{Las 9.3})-(\ref{Las 9.4})
have a unique solution  $\uusiPhi\in W^{1,1}_{loc}(\C)$  which}
satisfies for $0\leq q\leq p/4$
\beq\label{Las 9.5 A}
\dbar \uusiPhi \in X^q(\C)
\eeq
and the series (\ref{eq: Beurling series}) converges in $X^q(\C)$.
The convergence of the series (\ref{eq: Beurling series}) in $X^q(\C)$
is uniform for $\mu\in B^p_{exp,N}({B(R_0)})$ with any $N>0$.
 Moreover, for $\mu\in B^p_{exp,N}({B(R_0)})$  {\ltext the Jacobian
$J_\uusiPhi(z)$ of $\uusiPhi$ satisfies
\beq\label{Las 9.6}
\| J_ \uusiPhi\|_{X^{1,q}(B(R))}\leq C
\eeq
where $C$ depends on $p,q,N$, and $R$.}
Moreover, let $s>2$ and assume that $\mu_m,\tilde \mu_m\in B_{exp,N}^p({B(R_0)})$ and $0\leq q\leq p/4$.
Then it holds that
\beq\label{Las 9.6 B}
\lim_{m\to \infty}\|\mu_m- \tilde \mu_m\|_{L^s({B(R_0)})}=0\quad\Rightarrow\quad
\lim_{m\to \infty}\|\dbar  \uusiPhi_m-\dbar \tilde  \uusiPhi_m\|_{X^q(\C)}=0\
\eeq
where $ \uusiPhi_m$ and $\tilde  \uusiPhi_m$ are the solutions of
 (\ref{Las 9.3})-(\ref{Las 9.4}) corresponding
to $\mu_m$, $\tilde \mu_m$, respectively.
\end{theorem}

{\bf Proof.} 
{\ltext Let $ \uusiPhi^\alpha(z)$, $|\lambda|\leq 1$, $z\in \C$ be
the principal solution
corresponding to  
the Beltrami coefficient
$\lambda\mu$, that is, solution with the 
 Beltrami equation (\ref{Beltrami.principal})-(\ref{Beltrami.principal2}) with coefficient $\lambda\mu$.
 These solutions, in particular  $ \uusiPhi^\lambda=\uusiPhi^1$, exist and are unique
 by Thm.\ \ref{Stoilow}.}
 It follows from \cite[Thm.\ 1.1 and 5.1]{AGRS}
that the Jacobian  
determinant $J \uusiPhi^\alpha(z)$ of $ \uusiPhi^\lambda$ satisfies 
\beq\label{Eq Lassi 11}
\int_{B(R_0)} J \uusiPhi^\lambda\, \log^{2q}(e+J \uusiPhi^\lambda)\, dm(z)\leq C<\infty,
\eeq
where $C$ is independent of $\lambda$ and $\mu\in B^p_{exp,N}({B(R_0)})$ and
depends only on $N,p$, and $q$.
Thus (\ref{Las 9.6}) follows from Lemma \ref{lem: LassiA} ii.

We showed already that  when $p>2$, $\dbar \uusiPhi \in L^2(\C)$ and that
the series (\ref{eq: Beurling series}) converges in $L^2(\C)$.
To show the convergence of (\ref{eq: Beurling series}) in $X^q(\C)$ and
to prove (\ref{Las 9.6 B}), we present  few lemmas 
in terms of Orlicz spaces $X^q({B(R_0)})$ and the function $M_q$ defined in (\ref{Las 2}).
Note that as $\mu$ vanishes in $\C\setminus {B(R_0)}$, $\|(\mu S)^n\mu \|_{X^q(\C)}
=\|(\mu S)^n\mu \|_{X^q({B(R_0)})}$.

\begin{lemma}\label{lem: LassiB} 
Let $N\in\Z_+$,  
$2<2q<\beta<p$, and $\mu \in B_{exp,N}^p({B(R_0)})$. Then
\beq \label{eq:Lassi 4a}
\int_{B(R_0)} M_q(\psi_n(x) )dm(x)\leq cn^{-(\beta-q)}<cn^{-q},
\eeq
where $\psi_n=(\mu S)^n\mu$ and $c>0$ depends only on $N,p,\beta$, and $q$.
\end{lemma}

{\bf Proof.}  Let 
$
E_n=\{z\in {B(R_0)};\ |\psi_n(z)|\geq A^n\},
$
where $A>1$ is a constant to be chosen later. By (\ref{Kari B1}),
\beq\label{Eq Lassi 5 modified}
\|\psi_n\|_{L^2({B(R_0)})}\leq C_{N,\beta,p}n^{-\beta/2}.
\eeq
{\ltext Thus
\beq\label{Eq Lassi 8} 
|E_n|\leq C^2_{N,\beta,p}A^{-2n}n^{-\beta}.
\eeq
Using (\ref{Eq Lassi 5 modified}) we obtain
\beq\label{Eq Lassi 9}
\hspace{-5mm}\int_{{B(R_0)}\setminus E_n}|\psi_n|^2\log^q(e+|\psi_n|)\,dm\leq \|\psi_n\|_{L^2({B(R_0)})}^2\log^q(e+A^n) \leq 
C_1 n^{-\beta+q}
\eeq
where $C_1=C^2_{N,\beta,p}\log^q(e+A)$.
}

{\ltext The principal solution
corresponding to 
the Beltrami coefficient
$\lambda\mu$ can be written in the form $ \uusiPhi^\lambda=(I-\lambda \mu S)^{-1}(\lambda \mu)$. 
Expanding  $ \dbar_z\uusiPhi^\alpha(z)$ as a power series in $\lambda$, we see that} by (\ref{eq: Beurling series}) we can write using any $\rho$, $0<\rho\,<1$ 
\ba
\chi_{E_n}(z)\psi_n(z)=\frac 1{2\pi i}\int_{|\lambda|=\rho}
\lambda^{-n-2} \chi_{E_n}(z) \dbar_z  \uusiPhi^\alpha(z)\,d\lambda.
\ea
This gives
\beq\label{Eq Lassi 10}
\|\chi_{E_n}\psi_n\|_{X^q({B(R_0)})}\leq \rho^{-(n+2)}\sup_{|\lambda|=\rho}
\|\chi_{E_n} \dbar_z  \uusiPhi^\lambda \|_{X^q({B(R_0)})}.
\eeq

Using the facts that $|\lambda|=\rho$ and that the Beltrami coefficient  $\uusiPhi^\lambda$ is bounded
by $|\lambda|$, we have, by the distortion
equality (\ref{eq: distortion ineq.}),
$
|\dbar_z  \uusiPhi^\alpha(z)|^2\leq {\rho^2}({1-\rho^2})^{-1}  J \uusiPhi^\alpha(z).
$
Hence,
\ba
I:=\int_{E_n} M_q(\dbar_z  \uusiPhi^\alpha(z))\,dm(z)\leq
\frac {\rho^2}{1-\rho^2} \int_{E_n}
J \uusiPhi^\lambda\log^q(e+\bigg(\frac {\rho^2}{1-\rho^2}J \uusiPhi^\lambda\bigg)^{\frac 12})\,dm.
\ea 
{\ltext Let $\hat C$ denote next a generic constant which is function of 
$N,\beta,p,q$ and $\rho$ but not of $A$. The above} 
implies by   (\ref{Eq Lassi 3}), (\ref{Eq Lassi 11}),  {\lltext and
the inequality $\log(e+t^{1/2})\leq 1+\log (e+t),$ $t\geq 0$,} 
\ba
I
&\leq&
\hat C
\int_{E_n} J \uusiPhi^\alpha(z)
(1+\log(e+J \uusiPhi^\lambda))^q\,dm(z)\\
&\leq&
\hat C
\left(\int_{E_n} J \uusiPhi^\alpha(z)\,dm(z)\right)^{1/2}
\left(\int_{E_n}  J \uusiPhi^\alpha(z)(1+\log(e+ J \uusiPhi^\alpha(z)))^{2q}\,dm(z)\right)^{1/2}\\
&\leq&
\hat C
\left(\int_{E_n} J \uusiPhi^\alpha(z)\,dm(z)\right)^{1/2}.
\ea
{\ltext By the area distortion theorem \cite{As}, 
\ba
\int_{E_n}  J \uusiPhi^\alpha(z)\,dm(z)
\leq | \uusiPhi^\alpha(E_n)|\leq \hat C |E_n|^{1/M}\leq \hat C A^{-2n/M},
\ea
where $M=(1+\rho)/(1-\rho)>1,$
{\lltext  and thus
$
I\leq  \hat C A^{-n/M}.
$
By Lemma \ref{lem: LassiA} (ii) also 
$\| \chi_{E_n}\dbar  \uusiPhi^\lambda\|_{X^q}\leq \hat C  A^{-n/M}$.
 We take $\rho>e^{-1/2}$ and
$A=e^{M}$ we} see using (\ref{Eq Lassi 10}) and 
Lemma \ref{lem: LassiA} again that also $\int_{E_n}M_q(\psi_n)dm(z)\leq \hat C e^{-n/2}$ 
as $n\to \infty$. Thus the assertion follows from (\ref{Eq Lassi 9}). \hfill \proofbox 
\medskip

Lemmas \ref{lem: LassiA} (ii) and  \ref{lem: LassiB} and the fact that $\mu$ vanishes outside ${B(R_0)}$ yield that  for any $N>0$}
the series (\ref{eq: Beurling series}) converges in $X^q(\C)$,
and the convergence of the series (\ref{eq: Beurling series}) in $X^q(\C)$
is uniform for $\mu\in B^p_{exp,N}({B(R_0)})$ where $q>1$ and $p>2q$.
Thus to prove Theorem \ref{thm: New 9.1} it remains to show (\ref{Las 9.6 B}).

\begin{lemma}\label{lem: aux}
Let $2<2q<p$, $N>0$, $2<\beta<p$, $s>2$,
$\mu,\nu\in B_{exp,N}^p({B(R_0)})$,  and $B_n=(\mu S)^n\mu-(\nu S)^n\nu$. Then 
\beq \label{eq:Lassi 4}
\sup_{n\in \Z_+} \int_\C M_q(B_n(x) )dm(x)\leq C
\eeq
where $C>0$ depends only on $N,p$, and $q$. Moreover,
there is $T>1$ such that
\beq\label{90 and half}
\|B_n\|_{L^2(\C)}\leq C_{N,\beta,p,s,T}
\min\left(
nT^{n}\|\mu-\nu\|_{L^s({B(R_0)})},n^{-\beta/2}\right).
\eeq
\end{lemma}

{\bf Proof.}
Lemmas \ref{lem: LassiA}  and \ref{lem: LassiB} yield (\ref{eq:Lassi 4}).
Next, let us observe that for $z\in\C$ 
\ba
B_n(z)=(\mu S)^n\mu-(\nu S)^n\nu=
\sum_{j=0}^{n} A_j(z),\quad A_j(z)=(\mu S)^j (\mu-\nu)(S\nu)^{n-j} \chi_{B(R_0)}.
\ea
As $\|\nu\|_{L^\infty}\leq 1$ and 
$\|S\|_q:=\|S\|_{L^q(\C)\to L^q(\C)}<\infty$
for $1<q<\infty$, we have that
\ba
&&\hspace{-5mm}\int_\C |A_j(z)|^q\,dm(z)
\leq (\|S\|_q^q)^j \int_{B(R_0)} 
|\mu(z)-\nu(z)|^q|((S\nu)^{n-j}\chi_{B(R_0)})(z)|^q\,dm(z)\\
&& \hspace{-5mm}\leq \|S\|_q^{jq}
\left(\int_{B(R_0)}  \hspace{-3mm}|\mu(z)-\nu(z)|^{q\rho}dm(z)\right)^{\frac 1\rho}
\left(\int_{B(R_0)}  \hspace{-3mm} |((S\nu)^{n-j}\chi_{B(R_0)})(z)|^{q\rho'}dm(z)\right)^{\frac1{\rho'}}
\ea
where $\rho^{-1}+(\rho')^{-1}=1$, $1<\rho<\infty$.
Thus 
\ba
\|A_j(z)\|_{L^q(\C)}
\leq (\|S\|_q)^j \|\mu-\nu\|_{L^{\rho q}({B(R_0)})}
(\|S\|_{q\rho'}^q)^{n-j}\|\nu \|_{L^{q\rho'}({B(R_0)})}^q,
\ea
where $\|\nu \|_{L^{q\rho'}({B(R_0)})}\leq \pi R_0^2$.
Thus by choosing $q=2$ and $\rho$ so that $s=q\rho>2$ yielding
$q\rho'=2s/(s-2)$, we obtain
\ba
\|(\mu S)^n\mu-(\nu S)^n\nu\|_{L^2(\C)}\leq 
(n+1)\pi^2R_0^4
 (1+\|S\|_{(2s/(s-2))}^2)^n \|\mu-\nu\|_{L^s({B(R_0)})}.
\ea
This and (\ref{Kari B1})  show that (\ref{90 and half}) is valid. 
 \hfill \proofbox \medskip

Now we are ready to prove (\ref{Las 9.6 B}) which
finishes the proof of Theorem \ref{thm: New 9.1}.
Let $B_{n,m}=(\mu_m S)^n\mu_m-(\tilde \mu_m S)^n\tilde \mu_m$.
By Schwartz inequality, (\ref{eq:Lassi 4}), (\ref{90 and half}) and Lemma \ref{lem: LassiA}  yield
\beq \label{eq:Lassi 6}
\int_{B(R_0)} M_q(B_{n,m}(z) )dm(z)&\leq&
\int_{{B(R_0)}} |B_{n,m}|^2 \log^{q}(e+ |B_{n,m}|)\,dm(z)\\
&\leq& \nonumber
\left(\int_{{B(R_0)}}  M_{2q}(B_{n,m}(z))dm(z)\right)^{1/2}\|B_{n,m}\|_{L^2({B(R_0)})}
\\
&\leq& \nonumber
C \min\left(
nT^{n}\|\mu_m-\tilde \mu_m\|_{L^s({B(R_0)})},n^{-\beta/2}\right)
\eeq
where $C$ depends only on $q,p,\beta,s,T,$ and $N$.

Let $\e>0$. As $\mu_m$ and $\tilde \mu_m$ vanish outside ${B(R_0)}$, 
\ba
\|\dbar  \uusiPhi_m-\dbar \tilde  \uusiPhi_m\|_{X^q(\C)}=
\|\dbar  \uusiPhi_m-\dbar \tilde  \uusiPhi_m\|_{X^q({B(R_0)})}\leq \sum_{n=0}^\infty \|B_{n,m}\|_{X^q({B(R_0)})}.
\ea
Thus by (\ref{eq:Lassi 6}) and
Lemma \ref{lem: LassiA} (ii) we can take $n_0\in \N$ so large that
for all $m$
\ba
\sum_{n=n_0}^\infty \|B_{n,m}\|_{X^q({B(R_0)})}\leq \frac \e 2.
\ea
Applying again (\ref{eq:Lassi 6}) and
Lemma \ref{lem: LassiA} (ii) we can choose $\delta>0$ so that
\ba
\sum_{n=0}^{n_0-1} \|B_{n,m}\|_{X^q({B(R_0)})}\leq \frac \e2\quad
\hbox{when }\|\mu_m-\tilde\mu_m\|_{L^s({B(R_0)})}\leq  \delta.
\ea
This proves Theorem \ref{thm: New 9.1}.
\hfill \proofbox

\begin{lemma}\label{Lassi lemma 1} 
Assume that $K_\mu$ corresponding to $\mu$ supported in $\D$
satisfies  (\ref{Las 9.1})  with $q,C_0>0$ and $R_1=1$. Let $\uusiPhi$ is the principal solution of
Beltrami equation corresponding to $\mu$.
Then
for all $\beta,R>0$
the inverse function $\uusiPsi=\uusiPhi^{-1}:\C\to \C$ of $\uusiPhi$ satisfies
\ba
\int_{B(R)} \exp(\beta K_{\mu}(\uusiPsi(z)))\,dm(z)<C
\ea
where $C$ depends on $q,C_0,\beta,$ and $R$.
\end{lemma} 
 
\medskip

\noindent {\bf Proof.}
Since $\uusiPhi$ satisfies the  condition ${\mathcal N}$
by \cite[Cor.\ 4.3]{AGRS}
we may change variable in integration to see that
\beq\label{Lassi star}
\int_{B(R)} \exp(\beta K_{\mu}(\uusiPsi(z)))\,dm(z)
=\int_{\uusiPsi(B(R))} \exp(\beta K_{\mu}(w)) J_\uusiPhi(w)\,dm(w).
\eeq
Using (\ref{eq: Kari 1}) for function $\principalPhi$ and $R>3$
we see that $\uusiPsi(B(R))\subset \tilde B=B(\tilde R)$, $\tilde R=R+1$.
By (\ref{Las 9.1}), $\exp(K_\mu(z))\in L^q(\tilde B)$ for all $q>1$
and thus {\ltext by (\ref{Las 9.6})}
$
J_\uusiPhi \in X^{1,q}(B(R))$ for $R>0.
$

Let us next use properties of Orlicz spaces and the notations discussed in the  Appendix {\ltext using
 a Young complementary pair $(F,G)$ where $F(t)=
\exp(t^{1/p})-1$ and $G(t)$ satisfies 
$G(t)= C_pt(\log(1+C_pt))^p$ for $t>T_p$ with suitable $C_p,T_p>0$,
see \cite[Thm.\ I.6.1]{KR}.}

By using $u(z)=\exp(\beta K_\mu(z))$ and $v=J_\uusiPhi(z)$ we
obtain from Young's inequality (\ref{app 9}) the inequality  
\beq\label{Lassi 2}
& &\hspace{-.7cm}\int_{B(R)} \exp(\beta K_{\mu}(\uusiPsi(z)))\,dm(z)\\
\nonumber
& &\hspace{-.7cm}\leq \int_{\tilde B} F(\exp(\beta K_{\mu}(w)))\,dm(w)+
 \int_{\tilde B} G(J_\uusiPhi(w))\,dm(w)\\
\nonumber
 & &\hspace{-.7cm}\leq \int_{\tilde B} \exp((\exp(\beta K_{\mu}(w)))^{1/p})\,dm(w)+
 \int_{\tilde B} C_pJ_\uusiPhi(w)(\log(1+C_pJ_\uusiPhi(w)))^p\,dm(w).
\eeq
We apply this by using  $p>\beta/q$, so that 
$
(\exp(\beta K_{\mu}(w)))^{1/p}
\leq \exp(qK_{\mu}(w)).
$ 
Thus 
\ba
\int_{\tilde B} \exp((\exp(\beta K_{\mu}(w)))^{1/p})\,dm(w)\leq
\int_{\tilde B} \exp(\exp(q K_{\mu}(w)))\,dm(w)<\infty.
\ea
The last term in (\ref{Lassi 2}) is finite by (\ref{Las 9.6}), and thus the 
claim follows.\hfill \proofbox 
\medskip

\subsection{Asymptotics of the phase function of the exponentially growing solution}

{\ltext Let $\mu\in B^p_{exp}({B(R_0)})$,   $k\in \C\setminus\{0\}$
and $\lambda\in \C$ satisfy $|\lambda|\leq 1$. Then using Lemmas  
\ref{lemma 2 A} and \ref{lemma 2},
with the affine weight $\A(t)=pt-p$ corresponding to the gauge function $Q$, we see
that the equation 
\beq\label{Las 9.3 B}
& &\dbar_z \uusif_k (z)=\lambda \frac {\overline k}k e_{-k}(z)\mu(z)\, \overline{\p_z \uusif_k (z)},\quad
\hbox{for a.e. } z\in \C,\\
& &\uusif_k (z)=e^{ikz}(1+O(\frac 1z)),\quad \hbox{as }|z|\to \infty\label{Las 9.4 B}
\eeq
has the unique solution $\uusif_k\in W^{1,Q}_{loc}(\C)$.
Moreover, this solution can be
written in the form
\beq\label{representation 3 BB}
 \uusif_k (z) = e^{ik\uusivarphi_k(z)}
\eeq
where $\uusivarphi_k:\C \to \C$ is a homeomorphism satisfying
\beq\label{eq: for log for ph BBi}
& & \dbar \uusivarphi_k(z)=-\frac {\lambda\, \overline k}{k}
\mu(z) e_{-k}( \uusivarphi_k(z)) \overline {\p  \uusivarphi_k(z)},\quad\hbox{for a.e. }z\in \C,\\
 & &\uusivarphi_k(z) = z +\mathcal O(\frac{1}{z}),\quad \hbox{as }|z|\to \infty.\label{eq: asymptotics BBi}
\eeq
We denote  below  $\uusif_k(z)=\uusif(z,k)$ and  $\uusivarphi_k(z)=\uusivarphi(z,k)$ 
and estimate next functions $\uusivarphi_k$
in  the Orlicz space $X^q(\C)$.} The following lemma is a generalization
of results of \cite{AP} to  Orlicz space setting.

\begin{lemma}\label{Lemma new 9.3} 
 Assume that $\nu\in B^p_{exp}({B(R_0)})$ 
for all $0<p<\infty$. 
{\lltext For $k\in \C\setminus\{0\}$ let $\uusiPhi_k\in W^{1,1}(\C)$ be the solution of
\beq\label{eq: uusi Psi}
& & \dbar \uusiPhi_k(z)=-\frac {\, \overline k\,}{k}
\nu(z)e_{-k}(z)  \p  \uusiPhi_k(z),\quad\hbox{for a.e. }z\in \C,\\
 & &\uusiPhi_k(z) = z +\mathcal O(\frac{1}{z}).\label{eq: asymptotics uusi Psi}
\eeq
  Then} for all $\e>0$ there exist $C_0>0$ 
such that
  $\dbar_z \uusiPhi_k (z)=g_k(z)+h_k(z)$ where 
   $g_k,h_k\in X^q(\C)$ are supported in ${B(R_0)}$ and
\beq\
& &\label{Las 9.7}
\sup_{k\in \C\setminus\{0\}}  \| h_k\|_{X^q} <\e,\\
& &\label{Las 9.8}
\sup_{k\in \C\setminus\{0\}}  \| g_k\|_{X^q} <C_0,\\
& &\label{Las 9.9}
\lim_{k \to  \infty}  \hat g_k(\xi)=0,
\eeq
where for all compact sets $S\subset \C$ the convergence in (\ref{Las 9.9}) is uniform for $ \xi \in S$.
\end{lemma}

\noindent {\bf Proof.} 
{\ltext  Let us denote $\tilde \nu_k(z)=\lambdapoistettu{\overline k}k^{-1} \nu(z)$ for $k\in \C\setminus\{0\}$.
Note that then for any $p>0$ there is $N>0$ such that $\tilde \nu_k( \,\cdotp,k)e_{-k}(\,\cdotp)\in B^p_{exp,N}({B(R_0)})$ 
for all $k\in \C\setminus\{0\}$.
By Theorem \ref{thm: New 9.1},
\ba
\lim_{n\to \infty}\| \dbar \uusiPhi_k-\sum_{n=0}^\infty (\tilde \nu_k e_{-k}S)^n(\tilde \nu_k e_{-k})\|_{X^q(\C)}=0
\ea
uniformly in $k\in \C\setminus\{0\}$.}
We define
\ba
g_k(z)=\sum_{n=0}^m (\tilde \nu_k e_{-k}S)^n(\tilde \nu_k e_{-k}),\quad
h_k(z)=\sum_{n=m+1}^\infty (\tilde \nu_k e_{-k}S)^n(\tilde \nu_k e_{-k})
\ea
For given $\e>0$ we can choose $m$ so large that (\ref{Las 9.7}) holds for all $k\in \C\setminus\{0\}$ and then using Lemma \ref{lem: LassiB}  choose
{\ltext  $C_0$ so that (\ref{Las 9.8}) holds for all $k\in \C\setminus\{0\}$.}

Next, we show  (\ref{Las 9.9}) when $\e$ and $m$ fixed so that
 (\ref{Las 9.7}) and  (\ref{Las 9.8}) hold.  
{\ltext We can write
\ba
g_k(z)=\sum_{n=0}^m e_{-nk}G_n,\quad G_n=\left (\frac {\overline k}k\lambdapoistettu\right)^{n+1} 
\nu S_{n}(k)\nu\dots \nu S_1(k)\nu
\ea
where $S_j(k)$ is} the Fourier-multiplier 
\ba
(S_j(k)\phi)^{\wedge}(\xi)=m(\xi+jk)\widehat \phi(\xi), \quad
m(\xi)=\frac {\overline \xi}\xi.
\ea
The proof of  \cite[Lemma 7.3]{AP}  for $n\ge 1$ and the Riemann-Lebesgue lemma for $n=0$ yields that for any $\tilde \e>0$
there {\ltext exists $R(n,\tilde \e)\geq 0 $ such that} 
\ba
|\widehat  G_n(\xi)|\leq (n+1)\kappa^n\tilde \e,\quad\hbox{for }|\xi|>R(n,\tilde \e),
\ea
where $\kappa=\| \nu \|_{L^\infty}\leq 1$. 
Thus
\beq\label{Las uus1}
|\widehat  G_n(\xi)|\leq (m+1)\tilde \e,\quad\hbox{for }|\xi|>R_0=\max_{n\leq m}R (n,\tilde \e),\ n=0,1,2,\dots,m.
\eeq
 As $(e_{-nk} G_n)^{\wedge}(\xi)=\widehat  G_n(\xi-nk)$,
 we see that for any $L>0$ there is $k_0>0$ such that
 if $|k|>k_0$ then $j|k|-L>R_0$ for $1\leq n\leq m$. Then it
 follows from (\ref{Las uus1}) that if $|k|>k_0$, then 
 \ba
 \sup_{|\xi|< L}|\widehat  g_k(\xi)|\leq (m+1)^2\tilde \e.
 \ea
This proves the limit (\ref{Las 9.9}), with the convergence being uniform
for $\xi$ belonging in a compact set. 
 \hfill \proofbox

\begin{proposition}\label{Lemma new 9.4}
{\lltext Assume that $\nu\in B^p_{exp}({B(R_0)})$ with $p>4$ 
 and $\uusiPhi_k(z)$ be the solution
 of (\ref{eq: uusi Psi})-(\ref{eq: asymptotics uusi Psi}).}
Then
\beq\label{Las 9.12 A}
\lim_{k\to \infty}\uusiPhi_k(z)=z\quad \hbox{uniformly for $z\in \C$}.
\eeq
\end{proposition}

\noindent {\bf Proof.} Step 1. We will first show that for all {\lltext $q$ with $4<q<p$ 
 we have
 $\dbar_z \uusiPhi_k(z) \to 0$ weakly
in $X^q(\C)$ as $k\to \infty$.  {\ltext Let $\uusieta \in X^{-q}(\C)$ and  $\e_1>0$.
By Theorem \ref{thm: New 9.1}, there is $C_1>0$ such that $\sup_k \| \dbar \uusiPhi_k\|_{X^q}\leq C_1$.
Since $C^\infty_0(\C)$ is dense in $X^{-q}(\C)$, {\lltext cf.\ \cite[Sec.\ II.10]{KR}}, we can find
a function $\uusieta _0\in C^\infty_0(\C)$ such that $\| \uusieta -\uusieta _0\|_{X^{-q}}\leq \min(1,\e_1/C_1)$.
Then
\beq
\label{Las 9.12}
|\bra \uusieta ,\dbar \uusiPhi_k\cet| \leq |\bra \uusieta _0,\dbar \uusiPhi_k\cet|+\| \uusieta -\uusieta _0\|_{X^{-q}(\C)} \|\dbar \uusiPhi_k \|_{X^{q}(\C)},
\eeq
where the second term on the right hand side is smaller than $\e_1$. 
Moreover, by Lemma \ref{Lemma new 9.3},
we can write $\dbar \uusiPhi_k=h_k+g_k$ {\ltext so that
(\ref{Las 9.7})-(\ref{Las 9.9}) are satisfied for $\e=\e_1(\|\uusieta\|_{X^{-q}}+1)^{-1}$ and some $C_0>0$.} Then 
$|\bra \uusieta _0,h_k\cet |\leq \e_1$.}

Since $\hat \uusieta _0$ is a rapidly decreasing function, $\hat g_k(\xi)$  is uniformly bounded for $\xi\in \C$
and $k\in \C\setminus\{0\}$ by Lemma \ref{Lemma new 9.3}, and $\hat g_k\to 0$ uniformly
in all bounded domains as $k\to \infty$, we see that
\beq\label{Las 9.13}
\bra \uusieta _0,g_k\cet= \bra \hat \uusieta _0,\hat g_k\cet\to 0\quad \hbox{as }k\to \infty
\eeq
Combining these, we see that   $\bra \uusieta _0,\dbar \uusiPhi_k\cet \to 0$  as $k\to \infty$,
and thus $\dbar_z \uusiPhi_k(z) \to 0$ weakly
in $X^q(\C)$ as $k\to \infty$.

Step 2.  Next we show the pointwise convergence 
\beq\label{Las 9.13 B}
\lim_{k\to \infty} \dbar_z \uusiPhi_k(z) =0.
\eeq
For this end, we observe that the function $\uusieta _z(w)=\pi^{-1}(w-z)^{-1}\chi_{B(R_0)}(w)$
satisfies $\uusieta _z\in X^{-q}(\C)$ for $q>1$. Since $ \uusiPhi_k(z)- z= \mathcal O(\frac{1}{z}) $ and
$\dbar \uusiPhi_k$ is supported in $ \overline {B(R_0)}$,
we have 
\beq\label {Las 9.13 B2}
\uusiPhi_k(z)=z-\frac 1\pi  \int_{B(R_0)} (w-z)^{-1}\dbar_w \uusiPhi_k(w)\,dm(w)=z-\bra \uusieta _z,\dbar \uusiPhi_k\cet.
\eeq
{\lltext As $\dbar \uusiPhi_k \to 0$ weakly
in $X^q(\C)$, 
we see (\ref{Las 9.13 B}) holds} for all $z\in \C$. 

Step 3. 
{\ltext By (\ref{modulus of cont.}) and (\ref{Kari 1}) we see that
the family $\{\uusiPhi_k(z)\}_{k\in \C\setminus\{0\}}$ of homeomorphisms has a uniform modulus
of continuity in compact sets.}
Moreover, since 
\ba
 \sup_k \| \dbar \uusiPhi_k\|_{L^1(\C)}\leq  \sup_k \| \dbar \uusiPhi_k\|_{X^q({B(R_0)})}= C_2<\infty
 \ea
{\ltext we obtain by (\ref{Las 9.13 B}) for $|z|>R_0+1$}
\beq\label {Las 9.13 C}
|\uusiPhi_k(z)-z|= |\bra \uusieta _z,\dbar \uusiPhi_k\cet|\leq \frac C{|z|}\| \dbar \uusiPhi_k\|_{L^1(\C)}
\leq \frac {CC_2}{|z|}.
\eeq
Thus, as the functions $\{\uusiPhi_k(z)\}_{k\in \C\setminus\{0\}}$ are uniformly equicontinuous in compact 
sets,
(\ref{Las 9.13 C}) and the pointwise convergence (\ref{Las 9.13 B}) yield
the 
uniform convergence (\ref{Las 9.12 A}).
 \hfill \proofbox

\subsection{Properties of the solutions of the non-linear Beltrami equation}

Let $\a\in \C$, $|\lambda|\leq 1$ and $\mu(z)$ by supported in ${B(R_0)}$, $R_0\geq 1$, and assume 
that $K=K_\mu$ satisfies  (\ref{Las 9.1}) with $q,C_0>0$ and $R_1=1$, c.f. (\ref{distor2 mod}). {\lltext Motivated  
by Lemma \ref{lemma 2}, we consider next the solutions
 $\varphi_k$ of the equation
\beq\label{Las 9.14}
& &\dbar_z \varphi_\lambda(z,k)=-\lambda  \frac {\overline k}k 
\mu(z)e_{-k}(\varphi_\lambda(z,k))\,\overline {\p_z \varphi_\lambda(z,k)},\quad z\in \C,\\
\label{Las 9.14 B}
& & \varphi_\lambda(z,k)=z+O(\frac 1z)\quad \hbox{as }|z|\to \infty.
\eeq
Let  $\psi_\lambda(\,\cdotp,k)=\varphi_\lambda(\,\cdotp,k)^{-1}$ be the inverse function of  $\varphi_\lambda(\,\cdotp,k)$.
A simple computation based on differentiation of
the identity $\psi_\lambda(\varphi_\lambda(z,k),k)=z$ in the $z$ variable shows that then 
\beq\label{Las 9.15}
& &\dbar_z \psi_\lambda(z,k)=-\lambda  \frac {\overline k}k 
\mu(\psi_\lambda(z,k))e_{-k}(z)\,\p_z \psi_\lambda(z,k),\quad z\in \C,\\
\label{Las 9.15 B}
& & \psi_\lambda(z,k)=z+O(\frac 1z),\quad \hbox{as }|z|\to \infty.
\eeq
We consider the equations (\ref{Las 9.14}) and (\ref{Las 9.15})  simultaneously
by defining the sets 
\ba
B_\mu&=&\{(\varphi,\nu)\, ;\ |\nu|\leq |\mu|\hbox{ and
$\varphi:\C\to \C$ is a homeomorphism} 
\\
& &\quad \quad \quad \hbox{ with
$\dbar \varphi=\nu\overline {\p\varphi}$, $\varphi(z)=z+O(z^{-1})$}\}
\ea
and 
\ba
{\mathcal G}_\mu=\{g\in W_{loc}^{1,Q}(\C)\,;\ \dbar g=(\nu\circ\varphi^{-1}) \p g,\ g(z)=z+O(z^{-1}),
\ (\varphi,\nu)\in B_\mu\}.
\ea  %

{\lltext Now  $\exp(\exp(qK_\mu))\in L^1({B(R_0)})$
with some $0<q<\infty$ {\lltext and $|\nu|\leq |\mu|$ a.e.\ Then $K_{\nu}(z)\leq K_\mu(z)$ a.e.\
Let $\dbar \varphi=\nu \overline {\p\varphi}$ in $\C$, $\varphi(z)=z+O(z^{-1})$,
so that $\dbar \varphi=\tilde \nu {\p\varphi}$ with
$|\tilde \nu(z)|=|\nu(z)|$.
Then for $\psi=\varphi^{-1}$ we have $K(z,\psi)=K_\nu(\psi(z))$.
 Thus by Lemma \ref{Lassi lemma 1}  we have 
\beq\label{Las 9.16}
\sup_{g\in {\mathcal G}_\mu} \| \exp(\beta K(\,\cdotp,g))\|_{L^1(B(R))}=
\sup_{(\varphi,\nu)\in B_\mu} \| \exp(\beta K_\nu\circ \varphi^{-1})\|_{L^1(B(R))}<\infty
\eeq
for all $\b>0$ and $R>0$. Using this and Theorem \ref{Stoilow},
we see that  the functions $g\in {\mathcal G}_\mu$
are homeomorphisms. Moreover,
 for $g\in {\mathcal G}_\mu$ the condition $g\in W_{loc}^{1,Q}(\C)$ 
 implies
 that $Dg\in X^{-1}_{loc}(\C)$. Furthermore by  (\ref{Las 9.6}), we have
\beq\label{Las 9.17}
\sup_{(\varphi,\nu)\in B_\mu} \| J_\varphi\|_{X^{1,q}(B(R))}<\infty
\eeq
for all $q>0$.}

\begin{lemma}\label{Lemma new 9.5}
The set ${\mathcal G}_\mu$ is relatively compact in the topology of uniform convergence.
\end{lemma}

\noindent {\bf Proof.} {\ltext 
Let  $( \varphi,\nu)\in B_\mu$ 
and  $\psi=\varphi^{-1}$ and $\overline \p \tilde g=(\nu\circ \varphi^{-1})\p g$, $g(z)=z+O(z^{-1})$. 

As $\mu$ is supported in $B(R_0)$, the function $\varphi$ is analytic outside $\overline B(R_0)$
we see using (\ref{eq: Kari 1}) for function $\varphi$ that for $R>0$ we have
 $\varphi(B(R))\subset B(R+3R_0)$,  $\psi(B(R))\subset B(R+3R_0)$, and that
$\psi$  is analytic outside $\overline {B(4R_0)}$.

Thus   (\ref{eq: Kari 1}) and  the same arguments which we used to prove  the estimate 
(\ref{eq: essential2})  yield that for $R>0$
\beq\label{Las 9.18}
\hspace{-5mm}
\| Q(|Dg|)\|_{L^1(B(R))}&\leq&
 \pi (R+3R_0)^2+
 \int_{B(R)} \exp(qK_{\nu}(\psi(w))-q) \, dm(w)\\
 \nonumber
&\leq&
 \pi (R+3R_0)^2+
 \int_{B(R+3R_0)}  \exp(qK_{\nu}(z)-q)J_\varphi(z) \, dm(z),
 \eeq where $Q(t)=|t|^2/\log (|t|+1)$.}
{\lltext We will next use Young's inequality (\ref{app 9}) with the admissible pair
$(F,G)$ where (c.f.\ \cite[Ch.\ 1.3]{KR})
\beq\label {Las 9.23}
F(t)=e^t-t-1,\quad G(t)=(1+t)\log (1+t)-t.
\eeq
 By Young's inequality,
we have 
\ba 
& &\int_{B(R+3R_0)} \exp(q K_{\nu}(z)-q)J_\varphi(z) \,dm(z)
\\
 && \hspace{-8mm}\leq \int_{B(R+3R_0)} \hspace{-3mm}\exp(\exp(qK_{\nu}(w)-q))\,dm(w)+
 \int_{B(R+3R_0)} \hspace{-3mm}(1+ J_\varphi(w))\log(1+J_\varphi(w))\,dm.
\ea
This,}   (\ref{Las 9.1}), and (\ref{Las 9.18}) show that there is a constant $C(R,\mu)$
such that for $g\in {\mathcal G}_\mu$ 
\beq\label{Las 9.19}
\| Q(|Dg|)\|_{L^1(B(R))}\leq C(R,\mu).
\eeq
{\lltext As $g\in {\mathcal G}_\mu$ are homeomorphism, this and
(\ref{modulus of cont.})  imply that
functions $g\in {\mathcal G}_\mu$ are equicontinuous in compact sets of $\C$.
As $\supp(\nu\circ\psi)\subset B(4R_0)$,
 $g$ is analytic outside the disc $B(4R_0)$ and the inequality (\ref{eq: Kari 1}) yields for $R>0$ and $g\in {\mathcal G}_\mu$,
\ba
g(B(R))\subset B(R+12R_0).
\ea 
These} imply by Arzela-Ascoli theorem that the set $ \{g|_{B(R)}\,;\ g\in {\mathcal G}_\mu\}$
is relatively relatively compact in the topology of uniform convergence for any $R>0$.
Thus by using a diagonalization argument we see that for arbitrary
sequence $g_n\in {\mathcal G}_\mu$, $n=1,2,\dots$ there is a subsequence
$g_{n_j}$ which converges uniformly in all discs $B(R)$, $R>0$.
Finally, by Young's inequality we get using the same notations as in (\ref{Las 9.13 B2}) 
{\ltext for $|z|>4R_0+1$
\beq\label {Las 9.20}
|g_k(z)-z|&=&\left | \frac 1\pi  \int_{B(4R_0)} (w-z)^{-1}\dbar_w g_k(w)\,dm(w)\right|\\
\nonumber &\leq &\frac 1{\pi(|z|-1|)} 
\int_{B(4R_0)} (Q(|\dbar_w g_k(w)|)+G_0(1))\,dm(w)
\eeq 
where $Q(t)$ and $G_0(t)=|t|^2\log (|t|+1)$ form a Young complementary pair (c.f.\ Appendix).}
Thus
$
|g_k(z)-z|\leq  {C_\mu}/({|z|-1})$ for $|z|>4R_0+1.
$
Using this and the uniform convergence of $g_{n_j}$  in all discs $B(R)$, $R>0$,
we see that $g_n$ has a subsequence converging uniformly in $\C$.
\hfill \proofbox 
\medskip

\begin{theorem}\label{Theorem new 9.6}
{\lltext Let $\lambda,k\in \C\setminus\{0\},$ $|\lambda|=1$.}
Assume that $\varphi_\lambda(z,k)$ satisfies (\ref{Las 9.14})-(\ref{Las 9.14 B})
with $\mu$ supported in $\D$ which satisfies (\ref{Las 9.1})  with
$q>0$ and $R_1=1$. Then 
\ba
\lim_{k\to \infty}\varphi_\lambda(z,k)=z
\ea
uniformly in $z\in \C$ and $|\lambda|=1$.
\end{theorem}

 \noindent {\bf Proof.} 
{\ltext   Let $\psi_\lambda(\cdotp,k)$ be the inverse function of $\varphi_\lambda( \,\cdotp,k)$.
{\lltext It is sufficient to show that
\ba
\lim_{k\to \infty}\psi_\lambda(z,k)=z
\ea
uniformly in $z\in \C$ and $|\lambda|=1$.

 Then, $\psi_\lambda( \,\cdotp,k)$ is the solution of (\ref{Las 9.15})-(\ref{Las 9.15 B}). 
Denote $\nu(z)=-\lambda  {\overline k}k^{-1} \mu(z)$ 
and note that
$|\nu(z)|=|\mu(z)|$.} Hence $(\varphi_\lambda( \,\cdotp,k),\nu(\,\cdotp)\,e_{-k}(\,\cdotp)\, )\in B_\mu$ 
and   $\psi_\lambda( \,\cdotp,k)\in {\mathcal G}_\mu$.
Moreover, as $\varphi_\lambda( \,\cdotp,k)$ is homeomorphism in $\C$ and analytic outside 
in $B(1)$, it follows from (\ref{eq: Kari 1}) that 
$\varphi_\lambda( \,\cdotp,k)$ maps the ball $B(1)$  in to $B(3)$ and
moreover, its inverse $\psi_\lambda( \,\cdotp,k)$ maps the ball $B(3)$  in to $B(4)$
and $\C\setminus B(3)$ in to $\C\setminus B(2)$. 
}

{\ltext It follows from Lemma \ref{Lemma new 9.5} that if the claim is not valid, there are sequences 
 $(\lambda_n)_{n=1}^\infty$, $|\lambda_n|=1$ and $(k_n)_{n=1}^\infty$, $k_n\to \infty$ such that
 \ba
\psi_\infty(z)=\lim_{n\to \infty}\psi_{\lambda_n}(z,k_n),
\ea
 where the convergence is uniform $z\in \C$, such that $\psi_{\infty}(z)$ is not equal
 to $z$.} 
Thus, to prove the claim, it is enough to show that such a limit satisfies 
$\psi_{\infty}(z) =z$.
Note that by considering subsequences, we can assume that $\lambda_n\to \lambda$ and 
{\ltext ${\overline k_n}k_n^{-1}\to \beta$ as $n\to \infty$ where $|\lambda|=|\beta|=1$.
Denote next $\nu_0(z)=-\lambda  \beta \mu(z)$}

 Let us consider the solution of
 \beq\label{Mat A1}
& &\dbar_z \Phi_\lambda(z,k)= \nu_0(\psi_\infty(z)) e_{-k}(z)\,{\p_z \Phi_\lambda(z,k)},\\
\label{Mat A2}
& & \Phi_\lambda(z,k)=z+O(\frac 1z)\quad \hbox{as }|z|\to \infty.
\eeq
By Proposition \ref{Lemma new 9.4},  $\Phi_\lambda(z,k)\to z$ as $k\to \infty$ uniformly in $z\in \C$.
Since for every $z\in\C$ the function $\uusieta_z:w\mapsto \chi_{B(4)}(w)(z-w)^{-1}$ is in $X^{-q}(\C)$ for $q>1$,
we obtain using (\ref{Las 9.13 B2})
\beq\nonumber
|\psi_{\lambda_n}(z,k_n)-\Phi_{\lambda}(z,k_n)|\hspace{-2mm}&=&\hspace{-2mm}\frac 1\pi \left | 
 \int_{B(4)}  (w-z)^{-1}\dbar_w(\psi_{\lambda_n}(w,k_n)-\Phi_{\lambda}(w,k_n))\,dm(w)\,\right| \\
 \label{Matti eq A}
 & &\hspace{-2cm}\leq  |\uusieta_z\|_{X^{-q}(B(4))}  \| \dbar (\psi_{\lambda_n}( \,\cdotp ,k_n)-\Phi_{\lambda}( \,\cdotp,k_n))\|_{X^q(B(4))}.
 \eeq
{\ltext  Let us next assume that we can prove that
  \beq\label {Las 9.21}
\lim_{n\to \infty} \| \mu\circ  \psi_{\lambda_n}( \,\cdotp ,k_n)-\mu\circ  \psi_{\infty}( \,\cdotp ,k_n)\|_{L^s(\C)}=0
,\quad \hbox{for some $s>2$.}
\eeq
If this is the case, let $p\in (4q,\infty)$. {\lltext By assumption (\ref{Las 9.1}) and Lemma \ref{Lassi lemma 1}  
there is $N$ such that 
 the Beltrami coefficients of functions $\psi_{\lambda_n}( \,\cdotp ,k_n)$
 are in $B^p_{exp,N}(\D)$ for all $n\in \Z_+$ and $p>4$.} 
By Theorem \ref{thm: New 9.1} 
 and (\ref{Las 9.21}),
 \ba
 \lim_{n\to \infty}\| \dbar(\psi_{\lambda_n}( \,\cdotp ,k_n)-\Phi_{\lambda}( \,\cdotp,k_n))\|_{X^q(\C)}=0.
 \ea
 As $\lim_{n\to \infty}\Phi_{\lambda}(z,k_n)=z$ uniformly in $z\in \C$,
 this and (\ref{Matti eq A}) shows that
$ \psi_\infty(z)=z$.

Thus, to prove the claim it is enough to show (\ref{Las 9.21}).}
First, as $\psi_{\lambda_n}( \,\cdotp ,k_n)\to  \psi_{\infty}( \,\cdotp)$ uniformly
as $n\to\infty$ {\lltext and as $\psi_{\lambda_n}( \,\cdotp ,k_n)$ maps
and $\C\setminus B(3)$ in $\C\setminus B(2)$},
we see using Dominated convergence theorem that the   formula
(\ref {Las 9.21}) is valid when $\mu$ is replaced by a smooth compactly
supported function.  Next, let 
$(F,G)$ be the complementary Young pair given by (\ref{Las 9.23})
and $E_F(B(R))$ be the closure of $L^\infty(B(R))$ in $X_F(B(R))$.
 By \cite[Thm.\ 8.21]{Ad}, 
 the set $C^\infty_0(\D)$ is 
dense in  $E_F(\D)$ with respect to 
the norm of $X_F$. Thus when is $\mu$ a non-smooth Beltrami coefficient satisfying
the assumption (\ref{Las 9.1})} and 
$\e>0$ we can find 
 a smooth function $\theta\in C^\infty_0(\D)$, $\|\theta\|_\infty<2$ such that $\|\mu-\theta\|_F<\e$.
Then, since $|\mu-\theta|$ is supported in $\D$ and bounded by 3, 
{\lltext we have
\beq\nonumber
& &\hspace{-5mm}\| \mu\circ  \psi_{\lambda_n}( \,\cdotp ,k_n)-\theta\circ  \psi_{\lambda_n}( \,\cdotp ,k_n)\|_{L^s(\C)}^s
=\int_{\D} |\mu(z)-\theta(z)|^s J_{g_n}(z)\,dm(z)\\
& &\hspace{-5mm}\leq 3^{s-1}
(\int_{\D}F( |\mu(z)-\theta(z)|)\,dm(z))(  \int_{\D}G(J_{g_n}(z))\,dm(z) )\label {Las 9.22}
\eeq
where $g_n$ is the inverse of the function $\psi_{\lambda_n}( \,\cdotp,k_n)$.
}
Then, 
\ba
  \int_{\D}G(J_{g_n})\,dm \leq C\|J_{g_n}\|_{X^{1,1}(B(2))}\
\ea
and by (\ref{Las 9.17}), $ \|J_{g_n}\|_{X^{1,1}(\D)}$ is uniformly bounded in $n$. 
Using (\ref {Las 9.22})  and (\ref{app 5})
we see that}
 (\ref{Las 9.21})  holds
for all $\mu$ satisfying the assumption (\ref{Las 9.1}) and thus claim of the theorem follows.
%
%
%
%
%
%
\hfill
 \proofbox 
%
%
%
%

\subsection{$\dbar$ equations in $k$ plane}\label{subsec: 4.4}

Let us consider a Beltrami coefficient $\mu\in B^p_{exp}(\D)$ and 
 approximate $\mu$ with functions $\mu_n$ supported in $\D$ for
 which $\lim_{n\to\infty}\mu_n(z)=\mu(z)$ and $\|\mu_n\|_\infty\leq c_n<1$, \mattiHOX{complex
  conjugate added to (\ref{eq: repeated AA}), check citations}
 see e.g.\  
(\ref{truncted mu}). Let 
$f_{\mu}(\, \,\cdotp,k)\in W^{1,Q}_{loc}(\C)$  
be
the solution of
 the equations 
 \beq\label{eq: repeated AA}
& &\dbar_z f_{\mu}(z,k)=\mu(z) \overline {\p_z f_{\mu}(z,k)},\quad\hbox{for a.e. }z\in \C,
\\
& &\label{23b AA}
f_{\mu}(z,k)= e^{ikz}(1+\O_k(\frac 1 z ))\quad\hbox{for }|z|\to \infty,
\eeq
and $f_{\mu_n}(\, \,\cdotp,k)\in W^{1,Q}_{loc}(\C)$  be the solution of
the similar equations 
 to Beltrami coefficients ${\mu_n}$ and $\mu$, cf.\ Lemma \ref{lemma 2}.
 Here $\O_k(h(z))$ means a function of $(z,k)$ that satisfies
$|\O_k(h(z))|\leq C(k)|h(z)|$ for all $z$ with some constant $C(k)$
depending on
$k\in \C$.
Let 
\ba
\uusivarphi_{\mu}(z,k)=(ik)^{-1} \log (f_{\mu}(z,k)), \quad 
 \uusivarphi_{\mu_n}(z,k)=(ik)^{-1} \log (f_n(z,k)),
 \ea c.f.\ (\ref{eq: representation}).
 Then by (\ref{eq: Kari 1}) we have
\beq\label{modulus changed by three2}
|\uusivarphi_{\mu_n}(z,k)|\leq \, |z|+3,\quad |\uusivarphi_{\mu}(z,k)|\leq \,  |z|+3.
\eeq
By the proof of Lemma \ref{lemma 2} we see that
by choosing subsequence of $\mu_n$, $n\in\Z_+$, which we continue
to denote by $\mu_n$, we can assume that  
\beq\label{phi n limit}
\lim_{n\to \infty}\uusivarphi_{\mu_n}(z,k)=\uusivarphi_{\mu}(z,k),\quad& &\hbox{uniformly in $(z,k)\in B(R)\times \{k_0\}$,}\\ \nonumber
& &\hbox{for all $R>0$ and $k_0\in \C$.}
\eeq
Let us write the solutions $f_{\mu_n}$ and $f_{\mu}$ as 
\ba
& &f_{\mu_n}(z,k)=\lambdapoistettu e^{ik\uusivarphi_{\mu_n}(z,k)}=\lambdapoistettu e^{ikz}M_{\mu_n}(z,k),\\
& &f_{\mu}(z,k)=\lambdapoistettu e^{ik\uusivarphi_{\mu}(z,k)}=\lambdapoistettu e^{ikz}M_{\mu}(z,k).
\ea
Similar notations are introduced when $\mu$ is replaced by $-\mu$ etc.
{\ltext Let \ba
& &h^{(+)}_{\mu_n}(z,k)=\frac 12 (f_{\mu_n}(z,k)+f_{-\mu_n}(z,k)),
\quad h^{(-)}_{\mu_n}(z,k)=\frac i2(\overline {f_{\mu_n}(z,k)}-\overline  {f_{-\mu_n}(z,k)}),\\
& &
u^{(1)}_{\mu_n}(z,k)=h^{(+)}_{\mu_n}(z,k)-ih^{(-)}_{\mu_n}(z,k),
\quad u^{(2)}_{\mu_n}(z,k)=-h^{(-)}_{\mu_n}(z,k)+ih^{(+)}_{\mu_n}(z,k).
\ea
Then by (\ref{modulus changed by three2}), $h^{(+)}_{\mu_n}(z,k)$ and
$h^{(-)}_{\mu_n}(z,k)$ are uniformly bounded for $(z,k)\in B(R)\times \C$.
By (\ref{phi n limit}), we can define the pointwise limits
\beq\label{u n limit}
\lim_{n\to \infty} h^{(\pm)}_{\mu_n}(z,k)=h^{(\pm)}_{\mu}(z,k),\quad 
\lim_{n\to \infty} u^{(j)}_{\mu_n}(z,k)=u^{(j)}_{\mu}(z,k),\ \ j=1,2.
\eeq
The above formulae imply  
\beq\label{u12 relations}
\hbox{$u^{(2)}_{\mu}(z,k)=iu^{(1)}_{-\mu}(z,k)$ and  $u^{(1)}_{\mu}(z,k)=-iu^{(2)}_{-\mu}(z,k)$.}
\eeq
Moreover, for \ba
& &\tau_{\mu_n}(k)=\frac 12 (\overline{t_{\mu_n}(k)}-\overline{t_{-\mu_n}(k)}),\quad \tau_{\mu_n}(k)=\frac 12 (\overline{t_{\mu_n}(k)}-\overline{t_{-\mu_n}(k)}),\\
& &
t_{\pm\mu_n}(k)=\frac i {2\pi } \int_{\p \D} M_{\pm\mu_n}(z,k)\,dz,\quad 
t_{\pm\mu}(k)=\frac i{2\pi }  \int_{\p \D} M_{\pm\mu}(z,k)\,dz
\ea
we see using dominated convergence theorem} that 
$\lim_{n\to \infty}t_{\mu_n}(k)=t_{\mu}(k)$ for all $k\in \C$,
and hence 
\beq\label{tau n limit}
\lim_{n\to \infty}\tau_{\mu_n}(k)=\tau_{\mu}(k)\,\quad\hbox{for all $k\in \C$}.
\eeq
Then, as $|\mu_n|\leq c_n<1$ correspond 
to conductivities $\sigma_n$ satisfying $\sigma_n,\sigma_n^{-1}\in L^\infty(\D)$, we have by \cite[formula (8.2)]{AP}
the $\dbar$-equations with respect to the $k$ variables  
\beq\label{dbar tau n}
\dbar_k u^{(j)}_{\mu_n}(z,k)=-i \tau_{\mu_n}(k)\overline{u^{(j)}_{\mu_n}(z,k)},\quad k\in\ \C,\ j=1,2,
\eeq
see also Nachman \cite{N,N1} for different formulation on such equations.
{\ltext For $z\in \C$  functions $u^{(j)}_{\mu_n}(z,\,\cdotp)$, $n\in \Z_+$ are uniformly
bounded,
the limit  (\ref{u n limit}) and the  dominated convergence theorem 
imply that $u^{(j)}_{\mu_n}(z,\,\cdotp)\to u^{(j)}_{\mu}(z,\,\cdotp)$ as $n\to \infty$  in
$L^p(B(R))$ for all $p<\infty$ and $R>0$.
Since functions $|\tau_{\mu_n}(k)|$, $n\in \Z_+$ are uniformly bounded in compact sets, the pointwise limits (\ref{u n limit}), (\ref{tau n limit}) and the equation (\ref{dbar tau n}) 
yield that
\beq\label{dbar tau}
\dbar_k u^{(j)}_{\mu}(z,k)=-i \tau_{\mu}(k)\overline{u^{(j)}_{\mu}(z)},\quad k\in \C,\ j=1,2
\eeq
holds for all $z\in \C$ in sense of distributions and $u^{(j)}_{\mu}(z,\,\cdotp)\in W^{1,p}_{loc}(\C)$
for all $p<\infty$. }
}

\subsection{Proof of uniqueness results for isotropic conductivities}

\medskip

{\bf Proof of Theorem \ref{main 2b isotropic}.} 
Let us consider isotropic conductivities
$\sigma_j$, $j=1,2$.
Due to the above proven results, the proof will go along the lines of Section 8 of \cite{AP},
where $L^\infty$-conductivities are considered,
and its reformulation presented in Section 18 of \cite{AIM}   in a quite straight forward way 
when the changes explained below are made. The key point is
the following proposition.

\begin{proposition}\label{prop. 9.1} Assume that $\mu\in B^p_{exp}(\D)$
and let $f_{\pm \mu}(z,k)$ satisfy 
(\ref{eq: repeated AA})-(\ref{23b AA})
with the Beltrami coefficients $\pm\mu$.
Then
$
f_{\pm \mu}(z,k)=e^{izk}M_{\pm \mu}(z,k),
$ 
where 
\beq\label{strict inequality}
\re \frac {M_{+\mu}(z,k)}{M_{-\mu}(z,k)}>0
\eeq
for every $z,k\in \C$.
\end{proposition}
{\bf Proof.} Let us consider the Beltrami coefficients $\mu_n(z)$, $n\in \Z_+$ 
defined Section \ref{subsec: 4.4} that converge pointwise to $\mu(z)$ and satisfy
$|\mu_n|\leq c_n<1$.
By  Lemma \ref{lemma 1} the functions $M_{\pm \mu_n}(z,k)$ do
not obtain value zero anywhere.
By \cite[Prop.\ 4.3]{AP}, the inequality (\ref{strict inequality})
holds for the functions $M_{\pm \mu_n}(z,k)$. Then, 
$f_{\pm \mu_n}(z,k)\to f_{\pm \mu}(z,k)$ as $n\to \infty$
for all $k,z\in \C$, and thus
see that
\beq\label{non-strict inequality}
\re \frac {M_{+\mu}(z,k)}{M_{-\mu}(z,k)}=\lim_{n\to \infty}\re \frac {M_{+\mu_n}(z,k)}{M_{-\mu_n}(z,k)} \geq 0.
\eeq
To show that the equality does not hold in (\ref{non-strict inequality}), we assume
the opposite. In this case, there are  $z_0$
and $k_0$ such that
\beq\label{it equation}
M_{+\mu}(z_0,k_0)=it M_{-\mu}(z_0,k_0)
\eeq
with some $t\in \R\setminus \{0\}$. Then
\ba
f(z,k_0)=e^{ik_0z}(M_{+\mu}(z,k_0)-it M_{-\mu}(z,k_0))
\ea
is a solution of (\ref{eq: repeated AA}) and satisfies the asymptotics 
\ba
f(z,k_0)=(1-it)e^{ik_0z}(1+\O(\frac 1 z ))\quad\hbox{for }|z|\to \infty.
\ea
By using (\ref{eq: modified anisotropic Beltrami A1}) to write the
equation (\ref{eq: repeated AA}) in the form 
(\ref{eq: modified anisotropic Beltrami}) and 
applying  Lemma \ref{lemma 1}, we see that the solution $f(z,k_0)$ can be written
in the form 
\ba
f(z,k_0)=(1-it)e^{ik_0\varphi(z)}.
\ea
This is in contradiction with the assumption (\ref{it equation}) that implies $f(z_0,k_0)=0$
and  proves (\ref{strict inequality}).
\hfill \proofbox 
\medskip

Let $f_{\pm \mu}(z,k)$ be as in Prop.\  \ref{prop. 9.1} 
{\ltext and use below for the functions defined in (\ref{u n limit}) the 
short hand notation $u^{(1)}_\mu(z,k)=u_1(z,k)$ and $u^{(2)}_\mu(z,k)=u_2(z,k)$.   
Then  $u_1(z,k)$ and $u_2(z,k)$}
are solutions of the equation (\ref{dbar tau}). A direct computation 
shows also that 
\ba
\nabla\, \,\cdotp \sigma\nabla u_1(\, \,\cdotp,k)=0\quad\hbox{and}\quad
\nabla\, \,\cdotp \frac 1\sigma\nabla u_2(\, \,\cdotp,k)=0,
\ea
where $\sigma(z)=(1-\mu(z))/(1+\mu(z))$ is the conductivity corresponding
to $\mu$. Note that the conductivity $1/\sigma(z)=(1+\mu(z))/(1-\mu(z))$ 
 is the conductivity corresponding
to $-\mu$.

Generally, the near field measurements, that is, the Dirichlet-to-Neumann
map $\Lambda_\sigma$ on $\p \Omega$ determines the scattering measurements, in particular the scattered fields outside $\Omega$,
see Nachman \cite{N}. In our setting this means that we can
use Lemma \ref{lem: 3}  and argue, e.g. as in the proof of Proposition 6.1
 in \cite{AP} to see}, that $\Lambda_\sigma$
determines uniquely the solutions $f_{\pm \mu}(z_0,k)$ and $\tau_{\pm\mu}(k)$ for 
$z_0\in\C\setminus \overline \D$ and $k\in \C$. 
  We note that a constructive method based on integral equations on $\p \D$
 to determine $f_{\pm \mu}(z_0,k)$ from $\Lambda_\sigma$ is presented in \cite{AMPPS}.

As $u_j(z,\,\cdotp)$, $j=1,2$ are bounded and non-vanishing functions which satisfy (\ref{dbar tau}), we have
$\dbar u_j(z, \,\cdotp)\in L^\infty_{loc}(\C)$. This implies
that $\p u_j(z, \,\cdotp)\in \hbox{BMO}_{loc}(\C)\subset L^p_{loc}(\C)$ for all $p<\infty$, see e.g.\ \cite[Thm. 4.6.5]{AIM},
and hence $u_j(z, \,\cdotp)\in W^{1,p}_{loc}(\C)$.   

%
%

Let us now consider the isotropic conductivities $\sigma$ and $\tilde \sigma$ in $\Omega=\D$ 
which are equal to 1 near $\p \D$ and satisfy
(\ref{eq: expexp integrability}). 
Assume that $\Lambda_\sigma=\Lambda_{\tilde \sigma}$.
Then, by the above considerations, $\tau_{\pm \mu}(k)=\tau_{\pm\tilde \mu}(k)$ for $k\in\C$.

Let $\mu=(1- \sigma)/(1+\sigma)$ and $\tilde \mu=(1-\tilde \sigma)/(1+\tilde \sigma)$ be the Beltrami coefficients
corresponding to $\sigma$ and $\tilde \sigma$.

By applying Lemma \ref{lemma 2 A} with $k=0$ we 
see that $f_\mu(z,0)=1$ for all $z\in \C$ and hence $u_1(z,0)=1$. 
By Lemma \ref{lemma 1}, the  map $z\mapsto f_\mu(z,k)$ is continuous.
Thus $u_1\in  \mathcal X^p$, $1<p<\infty$, where
$\mathcal X^p$ is the space of functions $v(z,k)$, $(k,z)\in \C^2$ 
for which $v(z, \,\cdotp)\in W^{1,p}_{loc}(\C)$ and
$v(z,\,\cdotp)$ are bounded  for all $z\in \C$
and the function $v(\,\cdotp,k)$ is continuous for all $k\in \C$.
These properties are crucial in the following Lemma
which is a  reformulation of the properties
of the functions $u_1(z,k)$, $z,k\in \C$ proven in \cite{AP} 
for $L^\infty$-conductivities. 

\begin{lemma}\label{lem: asymptotics} 
(i) Functions $u_1(z,k)$ with   $k\neq 0$ have
the $z$-asymptotics
\begin{equation}\label{ass2}
u_1(z,k)=\exp(ikz+v(z;k)),
\end{equation}
where
 $C(k)>0$ is
such that $|v(z,k)|\leq C(k)$  for all $z\in \C$.

(ii) Functions $u_1(z,k)$ have  the $k$-asymptotics 
\begin{equation}\label{ass1}
u_1(z,k)=\exp(ikz+k\varepsilon_\mu (k;z)),\quad k\not=0
\end{equation}
where for  each fixed $z$  we have $\varepsilon_\mu(k;z)\to 0$ as $k\to\infty$.

(iii) Let $1<p<\infty$.
The  $u_1(z,k)$  given in (\ref{u n limit}) is the unique function in $\mathcal X^p$ 
such that  $u_1(z,k)$  non-vanishing,  $u_1(z,0)=1$ for all $z\in \C$,
and satisfies the $\dbar$-equation (\ref{dbar tau}) with the asymptotics 
and (\ref{ass2}) and (\ref{ass1}).
\end{lemma}

{\bf Proof.} 
(i) Let us omit the $(z,k)$ variables in some expressions and denote $u_1(z,k)=u_1$,
$f_\mu(z,k)=f_\mu$ etc.
By definition of $u_1$, 
\begin{equation}\label{sec8:eq84}
u_1
=\frac{1}{2}\left(f_{\mu}+f_{-\mu}+\ol{f_{\mu}}-\ol{f_{-\mu}}\right)
=f_{\mu}\left(1+\frac{f_{\mu}-f_{-\mu}}{f_{\mu}+f_{-\mu}}\right)^{-1}
\left(1+\frac{\ol{f_{\mu}}-\ol{f_{-\mu}}}{f_{\mu}+f_{-\mu}}\right),
\end{equation}
where each factor non-vanishing  by Prop.\ \ref{prop. 9.1}. Thus (\ref{23b AA}) yields (\ref{ass2}).

(ii) Let 
$F_t (z,k)= e^{-it/2}(f_{\mu} (z,k)\cos \frac  t2 + i f_{-\mu} (z,k)\sin \frac t2),$ $t\in\R.$
Then
\ba
& &\ol{\p}_zF_t(z,k) =\mu(z) e^{-it}\ol{\p_zF_t(z,k)},\quad\hbox{for }z\in \C,\\
& &F_t(z,k)=e^{ikz}\left(1+{O}_k \left({z}^{-1}\right)\right) \text{ as }z\to\infty.
\ea
Thus $F_t (z,k) =\exp(k \varphi_\lambda(z,k))$ where  $\lambda=e^{-it}$ and $\varphi_\lambda(z,k)$
solves (\ref{Las 9.14}). Note that $f_{\mu} (z,k) =\exp(k \varphi_{\lambda_0}(z,k))$
where $\lambda_0=1$. Then
\begin{equation}
\label{sec8:eq88}
\frac{f_{\mu}-f_{-\mu}}{f_{\mu}+f_{-\mu}}+e^{it}= 
\frac{2 e^{it}\, F_t }{f_{\mu}+f_{-\mu}} =  \frac{\exp(k \varphi_\lambda(z,k))}{\exp(k \varphi_{\lambda_0}(z,k))}\,  \frac{2 e^{it}}{1 + M_{-\mu}(z,k)/M_{\mu}(z,k)}.
\end{equation}
 By Theorem \ref{Theorem new 9.6} we have  for $z\in \C$ and $k\in \C\setminus\{0\}$
\beq
\label{sec8:eq89}
& &\hspace{2.5cm}e^{-|k|\varepsilon_{1}(k)} \leq |M_{\pm \mu}(z,k)|\leq e^{|k|\varepsilon_{1}(k)},\\
\label{sec8:eq810}
& &\hspace{-1cm}e^{-|k|\varepsilon_{2}(k)}  \leq \inf_{|\lambda|=1} \left|\frac{\exp(k \varphi_\lambda(z,k))}{\exp(k \varphi_{\lambda_0}(z,k))}\right| \leq \sup_{|\lambda|=1} \left|\frac{\exp(k \varphi_\lambda(z,k))}{\exp(k \varphi_{\lambda_0}(z,k))}\right|\leq e^{|k|\varepsilon_{2}(k)}\hspace{1cm}
\eeq
where $\varepsilon_{j}(k)\to 0$ as $k\to\infty$. Since $\re(M_{-\mu}/M_{\mu})>0$,
estimates \eqref{sec8:eq88}, \eqref{sec8:eq89} and \eqref{sec8:eq810}
yield for $z\in \C$, $k\not=0$
\ba
\inf_{t\in \R}\left|\frac{f_{\mu}-f_{-\mu}}{f_{\mu}+f_{-\mu}}+e^{it}\right|\geq
e^{-|k|\varepsilon(k)}\quad\hbox{and}\quad \frac{\,|f_{\mu}-f_{-\mu}|\,}{|f_{\mu}+f_{-\mu}|}\leq
1-e^{-|k|\varepsilon(k)}.
\ea
This and (\ref{sec8:eq84}) yield the $k$-asymptotics (\ref{ass1}).
  
%

(iii) As observed above, the function $u_1(z,k)$  given in (\ref{u n limit})  satisfies the conditions stated in (iii). 

Next, let $u_1(z,k)$ and $\tilde u_1(z,k)$ be two function which satisfy
the assumptions of the claim.
Let us  consider  the logarithms
\ba
\delta_1 (z,k)=\log u_1(z,k),\quad \tilde \delta_1 (z,k)=\log \tilde u_1(z,k),\quad k,z\in \C. 
\ea

As $u_1(z, \,\cdotp)\in W^{1,p}_{loc}(\C)$ for some $p<\infty$  
and $u_1(z,\,\cdotp)$  is bounded and non-vanishing function,
we see that $\delta_1(z, \,\cdotp)\in W^{1,p}_{loc}(\C)$.
As $u_1(z,0)=1$, we have
\beq\label{mu vanishes}
 \delta_1(z,0)=0,\quad\hbox{for }z\in \C.
\eeq
Moreover,
 $z\mapsto \delta_1(z,k)$ is continuous for any 
 $k$.  
Let $k\not =0$ be fixed. Then by (\ref{ass2})
\beq\label{mu asympt} 
\delta_1(z,k)=ikz+v(z,k),\quad z\in \C,
\eeq
where $v(\, \,\cdotp,k)$ is bounded and we see using  elementary homotopy theory \cite{asteteoria}
 that the map $H_{k}:\C\to \C$, $H_{k}(z)= \delta_1(z,k)$ is a surjective.

The function $\tilde\delta_1(z,k)$ has the same above properties as $\delta_1(z,k)$.
Next we want to show that 
$\delta_1(z,k)=\tilde\delta_1(z,k)$
for  all $z\in\C$ and $k\not=0$. As the map $H_{k}:z\mapsto \delta_1(z,k)$ is surjective 
for all $k\not=0$, this follows if we  show that
\beq\label{inclusion}
w\not=z\hbox{ and }k\not=0\ \Rightarrow\ \delta_1(w,k)\not=\tilde\delta_1(z,k).
\eeq
For this end, let $z,w\in \C$, $z\not=w$.
Functions $u_1$ and $\tilde u_1$ satisfy the same equation (\ref{dbar tau}) with the coefficient 
$\tau(k)=\tau_\mu(k)$.
Subtracting these equations from each other we see that 
the difference $g(k;w,z) =\delta_1 (w,k)-\tilde\delta_1 (z,k)$ satisfies
\beq\label{final dbar equation}
& &\dbar_k g(k;w,z)=\gamma(k;w,z)\,g(k;w,z),\quad k\in \C, \\
\nonumber& &\gamma(k;w,z)=
-i\tau(k)\exp(i\,\Im \delta_1(k;w,z)) E(i\,\im g(k;w,z)),
\eeq
where $E(t)=(e^{-t}-1)/t.$ Here, $\gamma(\,\cdotp ;w,z)$ is a locally bounded function.
As $w\not =z$,  the principle of the argument for pseudo-analytic functions, see 
\cite[Prop.\ 3.3]{AP}, the equation (\ref{final dbar equation}),
the boundedness of $\gamma$, and  the asymptotics
$g(k;w,z)=ik(w-z)+k\e(k,w,z)$, where $\e(k,w,z)\to 0$ as $k\to\infty$
imply that   $k\mapsto g(k;w,z)$ 
 vanishes
for one and only one value of $k\in \C$. Thus by (\ref{mu vanishes}),
$g(k;w,z)=0$ implies that $k=0$, and hence (\ref{inclusion}) holds.
Thus 
$\delta_1(z,k)=\tilde\delta_1(z,k)$
and $u_1(z,k)=\tilde u_1(z,k)$ for all  $z\in \C$ and $k\not=0$. \hfill \proofbox \medskip

{\bf Remark  4.1.}
Note that $\tau_{\pm\mu}(k)$ is determined
by $\Lambda_\sigma$. Thus Lemma \ref{lem: asymptotics}  means that  $u_1(z,k)$  can
constructed as the unique  complex curve   $z\mapsto u_1(z,\,\cdotp),$  ${z\in \C}$
in the space of the solutions of 
the $\dbar$-equation (\ref{dbar tau}) which has the properties stated in (iii).
 \medskip

When $u_j(z,k)$ and $\tilde u_j(z,k)$, $j=1,2$ are functions corresponding to $\mu$ and $\tilde \mu$,
the above shows that $u_1(z,k)=\tilde u_1(z,k)$. 
Using $\tau_{-\mu}$ instead of 
$\tau_{\mu}$ and (\ref{u12 relations}), we see by Lemma \ref{lem: asymptotics} that  $u_2(z,k)=\tilde u_2(z,k)$
for all  $z\in \C$,  and $k\not=0$.

%
%
%

Thus $f_{\pm \mu}(z,k)=f_{\pm \tilde \mu}(z,k)$ for all $z\in \C$
and $k\not=0$. {\mmmtext By \cite[Thm.\ 20.4.12]{AIM},
%
the Jacobians 
of $f_{\pm \mu}\in W^{1,Q}_{loc}(\C)$ are non-vanishing
almost everywhere.} Thus we  see using the Beltrami equation
(\ref{eq: repeated AA}) and the fact that $f_{\pm \mu}(z,k)=f_{\pm \tilde \mu}(z,k)$ for all $z\in \C$
and $k\not=0$
that  $\mu=\tilde \mu$ a.e. Hence $\sigma=\tilde \sigma$ a.e.
{\ltext This proves the claim of Theorem \ref{main 2b isotropic}.}
  \hfill \proofbox

\section{Reduction of the inverse problem for an anisotropic conductivity to
the isotropic case}\label{subsec: reduction to scalar case}

Let this section, we assume that the weight function $\A$ satisfies the
almost linear growth condition  (\ref{cond for A function}).
Let $\sigma=\sigma^{jk}\in \Sigma_{\mathcal A}(\C)$ be 
a conductivity  matrix for which
and  $\sigma(z)=1$ for $z$ in $\C\setminus \Omega$
and in some neighborhood of $\p \Omega$.

Let $z_0\in \p \Omega$, and define
\beq\label{EQ 9}
\H_\sigma(z)=\int_{\eta_z}(\Lambda_\sigma (u|_{\p \Omega}))(z')\,ds(z'),
\eeq
where $\eta_z$ is the path (oriented to positive direction)
from $z_0$ to $z$ along $\p \Omega$.
This map is called the $\sigma$-Hilbert transform, and
it can be considered a bounded map
\beq\label{EQ 9b}
\H_\sigma: H^{1/2}(\p \Om)\to H^{1/2}(\p \Om)/\C.
\nonumber
\eeq  
{\mmmtext As shown in
beginning of Subsection \ref{sec: Invariance},
there exists a homeomorphism  $F:\C\to \C$ such that $F(\Omega)=\tilde \Omega$,
$\tilde \sigma =F_*\sigma$ is isotropic (i.e.\ scalar function times identity matrix),
$F$ and $F^{-1}$ are $W^{1,P}$-smooth, and $F(z)=z+O(1/z)$.
 Moreover, $F$ satisfies conditions $\mathcal N$ and 
$\mathcal N^{-1}$,
Also, as $\sigma$ is one near the boundary, we have}
that $F$ and $F^{-1}$ are  $C^\infty$ smooth near the boundary.

{\lltext By definition of $\tilde\sigma= F_*\sigma$, we see that 
\beq
\det(\tilde \sigma(y))=\det(\sigma(F^{-1}(y)))\label{eq: det invaritant}
\eeq
for $y\in \tilde \Omega.$  
Thus under the assumptions of Thm.\ \ref{main 2b trace} where
$\det(\sigma),\det(\sigma)^{-1}\in L^\infty(\Omega)$ we see that
the isotropic conductivity $\tilde \sigma$ satisfies
 $\tilde \sigma,\tilde \sigma^{-1}\in L^\infty(\Omega)$.

Let us next consider the case when assumptions of Theorem \ref{main 2b both}
are valid and we have  $\A(t)=pt-p$, $p>1$. 
Then, as $F$ satisfies the condition $\mathcal N$,
the area formula gives
\beq\label{eq: reduction1}
I_1&= &\int_{\tilde \Omega}\exp(\exp(q(\tilde\sigma(y)+\frac 1{\tilde\sigma(y)}))dm(y)\\
&=&\int_{\Omega}\exp(\exp (q(\det(\sigma(x))^{1/2} +\frac 1{\det(\sigma(x))^{1/2}})))J_F(x)\,dm(x).
\nonumber
\eeq

In the case when $\A(t)=pt-p$, $p>1$
\cite[Thm.\  1.1]{AGRS}
implies that 
 $J_F\log^\beta(e+ J_F)\in L^1(\Omega)$ for $0<\beta<p$.
     {\lltext  Then, Young's inequality (\ref{app 9}) with 
the admissible pair (\ref{Las 9.23}) implies that
 \beq\label{eq: reduction2A}
& &\int_{\Omega}\exp(\exp (q(\det(\sigma(x)) ^{\frac 12}+\frac 1{\det(\sigma(x))^{\frac 12}}  ))) J_F(x)\,dm(x)\\
& &\hspace{-.8cm}\leq\left(\int_{\Omega}\exp(\exp(\exp (q(\det(\sigma)^{\frac 12} +\frac 1{\det(\sigma)^{\frac 12}}))))\,dm\right) 
\left(\int_{\Omega}(1+J_F)\log (1+J_F)\,dm \right)\nonumber
\eeq
and if conductivity $\sigma$ satisfies (\ref{eq: exp^3 integrability}),
we see that $I_1$ is finite for some $q>0$.

Thus under assumptions of Theorem \ref{main 2b both}  we see
that $I_1$ is finite 
for the isotropic conductivity $\tilde \sigma$.}

}

Let $\rho=F|_{\p \Omega}$. 
It follows from Lemma \ref{lem: change of coord.} and (\ref {eq: push forward of Lambda}) that  $\rho_*\Lambda_\sigma=\Lambda_{\tilde \sigma}$.
Then, $ \H_{\tilde \sigma }h= \H_\sigma( h\circ \rho^{-1})$
for all $h\in H^{1/2}(\p \tilde \Om)$.

Next we seek for a function
$G_{ \Omega}(z,k)$, $z\in \C\setminus  \Om$, $k\in \C$ that satisfies 
\beq\label{EQ 10 B}
& &\overline \p_z G_{ \Omega}(z,k)=0 \quad \hbox{for }z\in \C\setminus\overline { \Om},\\
\label{EQ 10 Bb}
& &G_{ \Omega}(z,k)=e^{ikz}(1+\O_k(\frac 1{z})),\quad\hbox{as }z\to \infty,\\
\label{EQ 10 Bc}
& &\im G_{ \Omega}(\,\cdotp,k)|_{\p   \Om}=\H_{  \sigma}
(\re G_{ \Omega}(\,\cdotp,k)|_{ \p    \Om}).
\eeq
To study it, we consider a similar function $G_{\tilde \Omega}(\,\cdotp,k):
\C\setminus  {\tilde \Omega}\to \C$
corresponding
to the scalar conductivity $\tilde \sigma,$
which  satisfies
in  the domain $\C\setminus \overline {\tilde \Omega}$ 
the  equations (\ref {EQ 10 B})-(\ref {EQ 10 Bb})  
and the boundary condition
$\im G_{\tilde \Omega}(\,\cdotp,k)= \H_{\tilde \sigma }
(\re G_{\tilde \Omega}(\,\cdotp,k))$ on ${z\in \p  {\tilde \Omega}}$.
Below, let  $\tilde \mu=(1-\tilde \sigma)/(1+\tilde \sigma)$  be the Beltrami coefficient corresponding
to the conductivity $\tilde \sigma$. 
%

\begin{lemma}\label{lem: 3} Assume that $\sigma\in \Sigma_{\mathcal A}(\Omega)$ is one near $\p \Omega$. 
Then for all  $k\in \C$

(i) For $k\in \C$ 
and   $z\in \C\setminus \overline{\tilde \Omega}$ we have
$G_{\tilde\Omega}(z,k)=W(z,k)$ where  $W(\,\cdotp,k)\in W^{1,P}_{loc}(\C)$
is the a unique solution of 
\beq\label{EQ 11}
& &\overline \p_z W(z,k)=\tilde \mu(z)\overline {\p_z W(z,k)}, \quad \hbox{for }z\in \C,\\
\label{EQ 11 b}
& &W(z,k)=e^{ikz}(1+\O_k(\frac 1 z)),\quad\hbox{as }z\to \infty.
\eeq

(ii) The equations (\ref{EQ 10 B})-(\ref{EQ 10 Bc}) have a unique  solution $G_\Om(\,\cdotp,k)
\in C^\infty(\C\setminus \Om)$
and  $
G_\Om(z,k)=G_{\tilde\Omega}(F(z),k)$ for $z\in \C\setminus \Om$.

\end{lemma}

{\bf Proof.} 
The definition of the Hilbert transform 
$ \H_{\tilde \sigma }$ implies that
 any solution $G_{\tilde\Omega}(z,k)$ of (\ref{EQ 10 B})-(\ref{EQ 10 Bc}) can be extended to a solution $W(z,k)$
of (\ref{EQ 11}). On other hand, the restriction of the solution
$W(z,k)$ of  (\ref{EQ 11})-(\ref{EQ 11 b}) satisfies (\ref{EQ 10 B})-(\ref{EQ 10 Bc}). 
The equations 
(\ref{EQ 11})-(\ref{EQ 11 b})
have a unique solution
by Theorem \ref{thm:CGO}. {\lltext As the solution $W(\,\cdotp ,k)$ is analytic in $\C\setminus
\supp(\tilde \sigma)$, the claim (i) follows.

The claim (ii) follows immediately 
as $F:\C\setminus \overline \Omega\to
\C\setminus \overline {\tilde \Omega}$ is conformal, $F(z)=z+\O(1/z)$, and
  $ \H_{\tilde \sigma }h= \H_\sigma( h\circ \rho)$
for all $h\in H^{1/2}(\p \tilde \Om)$.
\hfill \proofbox \medskip

\begin{lemma}\label{lem: 4} 
Assume that $\Omega$ is given and that  
 $\sigma\in \Sigma_{\mathcal A}(\Omega)$ is one near $\p \Omega$. 
{\lltext Then 
the Dirichlet-to-Neumann form
$Q_\sigma$ determines the values of the restriction
$F|_{\C\setminus \Om}$,
the boundary $\p \tilde \Omega$, and the Dirichlet-to-Neumann map
$\Lambda_{\tilde \sigma}$ of} the isotropic conductivity $\tilde \sigma=F_*\sigma$ on 
$\tilde \Omega$.
\end{lemma}

{\bf Proof.} When $\sigma$ is identity near $\p\Omega$,
the Dirichlet-to-Neumann form
$Q_\sigma$ determines the
 Dirichlet-to-Neumann map
$\Lambda_\sigma$.
By  Lemma \ref{lemma 2}  we \
have 
$
W(z,k)=\exp(ik\varphi(z,k))
$
where by Theorem \ref{Theorem new 9.6}
\beq\label{EQ 18}
\lim_{k\to \infty} \sup_{z\in \C}|\varphi(z,k)-z|=0.
\eeq
As   $G(z,k)=W(F(z),k)$ we have
\beq\label{EQ 19}
\lim_{k\to \infty} \frac{\log G(z,k)}{ik}=
\lim_{k\to \infty} \varphi(F(z),k)=F(z).
\eeq

By Lemma \ref{lem: 3} $F(z,k)$ can be constructed for any $z\in \C\setminus\Omega$
by solving the equations
(\ref{EQ 10 B})-(\ref{EQ 11 b}). Thus the restriction of $F$ to $ \C\setminus\Omega$
is determined by the values
of limit (\ref{EQ 19}). As
$\tilde \Omega=\C\setminus F(\C\setminus \Omega)$
and $\Lambda_{\tilde \sigma}=(F|_{\p\Omega})_*\Lambda_\sigma$,
this proves the claim. \hfill \proofbox 
\medskip

{\lltext 
Above we saw that if the assumptions of 
Theorem \ref{main 2b both} for $\sigma$ are satisfied then for the isotropic conductivity 
$\tilde\sigma=F_*\sigma $  we have  $\tilde \sigma,\tilde \sigma^{-1}\in L^\infty(\tilde \Omega)$. 
Also,  under the assumptions of 
Theorem \ref{main 2b both} for $\sigma$ the integral $I_1$ in 
(\ref{eq: reduction1}) is finite for some $q>0$. Thus Theorems \ref{main 2b both} and  \ref{main 2b trace}
follow by 
 Theorem \ref{main 2b isotropic} and Lemma \ref{lem: 4}.}
 \hfill \proofbox 
\medskip

\section*{Appendix A: Orlicz spaces}

For the proofs of the facts discussed in this appendix we 
refer to \cite{Ad,KR}.
 
Let $F,G:[0,\infty) \to [0,\infty)$ be bijective convex functions. The pair $(F,G)$
is called  a Young complementary pair if 
\ba
F'(t)=f(t),\quad 
G'(t)=g(t),\quad g=f^{-1}.
\ea
In the following we will consider also extensions of
these functions defined by $F,G:\C \to [0,\infty)$ by setting
$F(t)=F(|t|)$ and $G(t)=G(|t|)$. {\ltext By \cite[Sec.\ I.7.4]{KR},
there are examples of such pairs for which $F(t)=\frac 1p t^p \log^a t$
and $G(t)=\frac 1q t^q \log^{-a} t$ where $p,q\in (1,\infty),$ $p^{-1}+q^{-1}=1$
and $a\in \R$.}
We define that  
 $u:D\to \C$, $D\subset \R^2$ is in an Orlicz class $K_F(D)$ if 
 \beq\label{app 1}
 \int_\D F(|u(x)|)\, dm(x)<\infty.
 \eeq
 The Orlicz space $X_F(D)$ is the smallest vector space
 containing the convex set $K_F(D)$.

For a
 Young complementary pair $(F,G)$ one can define for $u\in X_F(D)$ the norm
  \beq\label{app 2}
 \| u\|_F=\sup\{ \int_D |u(x)v(x)|\, dm(x)\ ;\ \int_D G(u(x))\, dm(x)\leq 1\}.
 \eeq
  There is also a Luxenburg norm
 \beq\label{app 3}
 \| u\|_{(F)}=\inf \{ t>0\  ;\ \int_D F(\frac {u(x)}t)\, dm(x)\leq 1\}
 \eeq
which is equivalent to the norm $\| u\|_F$ and
one always has 
\beq\label{app 4}
 \| u\|_{(F)}\leq  \| u\|_{F}\leq 2\| u\|_{(F)}.
 \eeq
By  \cite[Thm.\ 8.10]{Ad}, $L_X(D)$ is a Banach space
with respect to the norm $\| u\|_{(F)}$.
Moreover, it holds that (see \cite[Thm.\ II.9.5 and  II.10.5]{KR}),
\beq
\label{app 5}
& & \| u\|_{(F)}\leq 1\quad \Rightarrow\quad \int_D F(u(x))\, dm(x)\leq \| u\|_{F},\\
\label{app 6}
& &\| u\|_{(F)}\geq 1\quad \Rightarrow\quad \int_D F(u(x))\, dm(x) \geq  \| u\|_{(F)}.
 \eeq
We also recall the
Young's inequality 
\cite[Thm.\ II.9.3]{KR},
$
uv\leq F(u)+G(v)$ for $u,v\geq 0$
which implies
\beq
\label{app 9}
\left |\int_D u(x)v(x)\,dm(x)\right |\leq \| u\|_{F}\| u\|_{G}.
\eeq
The set $K_F(D)$ is a vector space when $F$
satisfies the $\Delta_2$ condition, that is, there is $k>1$ such that
 $F(2t)\leq kF(t)$ for all $t\in \R_+$, see \cite[Lem.\ 8.8]{Ad}. 
 In this case $X_F(D)=K_F(D)$. 

We will use functions
\ba
M_{p,q}(t)=|t|^p (\log(1+|t|))^{q},\quad 1\leq p<\infty,\ q\in \R
\ea
and use for $F(t)=M_{p,q}(t)$ the notations $X_F(D)=X^{p,q}(D)$
and  $\| u\|_F=\| u\|_{X^{p,q}(D)}$. 
For $p=2$ we denote 
$M_{2,q}(t)=M_{q}(t)$ and $
X^{2,q}(D)=X^{q}(D).$
Note that if $D$ is bounded, $1<p<\infty$ and $0<\e<p-1$
then 
\ba
L^{p+\e}(D)\subset X^{p,q}(D)\subset L^{p-\e}(D).
\ea
Finally, we note that the dual space of $X^q(D)$ is
 $X^{-q}(D)$ and 
\beq
\label{app 10}
\left |\int_D u(x){v(x)}\,dm(x)\right |\leq \| u\|_{X^q(D)}\| v\|_{X^{-q}(D)}.
\eeq
\bibliographystyle{amsalpha}

\end{document}